\documentclass[a4paper]{amsart}
\usepackage[utf8]{inputenc}
\DeclareUnicodeCharacter{03C4}{$\tau$}


\usepackage{afterpage}
\usepackage{hyperref}
\usepackage{amsmath,amssymb,amsthm}
\usepackage{microtype}
\usepackage{bbm}
\usepackage{graphicx}
\usepackage{rotating}
\usepackage{longtable}
\usepackage{booktabs}
\usepackage{tabularx}
\usepackage{fancyvrb}
\usepackage{tabu}

\usepackage{spectralsequences}
\usepackage{tikz-cd}
\usepackage[T1]{fontenc}
\usepackage{cleveref}

\usepackage[url=false,isbn=false,citestyle=authoryear,style=alphabetic]{biblatex}
\newtheorem{thm}{Theorem}[section]
\crefname{thm}{theorem}{theorems}
\newtheorem{lemma}[thm]{Lemma}
\newtheorem{cor}[thm]{Corollary}

\theoremstyle{definition}
\newtheorem{remark}[thm]{Remark}
\newtheorem{defi}[thm]{Definition}
\newtheorem{notation}[thm]{Notation}
\newtheorem{eg}[thm]{Example}

\newcommand\kk{{\mathbbm{k}}}
\newcommand\A{{\mathcal{A}}}
\newcommand\B{{\mathcal{B}}}
\newcommand\C{{\mathcal{C}}}
\newcommand\D{{\mathcal{D}}}
\newcommand\E{{\mathbb{E}}}
\newcommand\M{{\mathcal{M}}}
\renewcommand\H{{\mathcal{H}}}
\newcommand\Z{{\mathbb{Z}}}
\newcommand\F{{\mathbb{F}}}
\newcommand\HFp{{\F_p}}
\newcommand\kappabar{{\bar{\kappa}}}

\newcommand\Sp{{\mathrm{Sp}}}
\newcommand\Syn{{\mathrm{Syn}}}
\newcommand\Cat{{\mathrm{Cat}}}
\newcommand\Free{{\mathrm{Free}}}
\renewcommand\S{{\mathbb{S}}}
\newcommand\op{{\mathrm{op}}}
\newcommand\gr{{\mathrm{gr}}}
\newcommand\Stbl{{\Cat^{\mathrm{ex}}}}

\DeclareMathOperator\Sq{Sq}
\DeclareMathOperator\Tot{Tot}
\DeclareMathOperator\Hom{Hom}
\DeclareMathOperator\End{End}
\DeclareMathOperator\Aut{Aut}
\DeclareMathOperator\Pic{Pic}
\DeclareMathOperator\Ext{Ext}
\DeclareMathOperator\Fun{Fun}
\DeclareMathOperator\Mod{Mod}
\DeclareMathOperator\dgMod{dgMod}
\DeclareMathOperator\Comod{Comod}

\DeclareMathOperator\Alg{Alg}

\DeclareMathOperator\Ch{Ch}
\DeclareMathOperator*\colim{colim}

\newcommand\iq{{\text{?`}}}
\newcommand\mot{{m}}

\newcommand\prompt[1]{\textcolor{gray!70!black}{#1}}
\catcode`\$=13
\DefineVerbatimEnvironment{shell}{Verbatim}{%
  frame=lines,%
  codes={\catcode`\$=13},%
  defineactive=\def${\prompt{\$}},%
  commandchars=\\\{\}%
}
\catcode`\$=3

\newcommand\fakeqed{\pushQED{\qed}\popQED}
\newcommand\sidedeg{
    \begin{aligned}
      |\mu_0| &= 0, & d(\mu_0) &= p, \\
      |\tau| &= 1, & d(\tau) &= 0.
    \end{aligned}
}

\title{The \texorpdfstring{$E_3$}{E3} page of the Adams spectral sequence}
\author{Dexter Edralin Chua}

\begin{document}

\begin{abstract}
In \cite{baues-book}, Baues computed the secondary Steenrod algebra, the algebra of all secondary cohomology operations. Together with Jibladze, they showed that this gives an algorithm that computes all Adams $d_2$ differentials for the sphere \cite{baues-e3}.

The goal of this paper is to reinterpret their results in the language of synthetic spectra in order to achieve stronger computational results. Using this, we obtain an algorithm that computes hidden extensions on the $E_3$ page that jump by one filtration, in addition to the $d_2$ differentials of Baues--Jibladze. We then implement and run this algorithm for the sphere up to the 140\textsuperscript{th} stem.

Combined with a generalized version of the Leibniz rule, these hidden extensions allow us to compute many longer differentials with ease. In particular, we resolve all remaining unknown $d_2$, $d_3$, $d_4$ and $d_5$ differentials of the sphere up to the 95\textsuperscript{th} stem.

\end{abstract}
\maketitle

\tableofcontents

\section{Introduction}
One of the most fundamental problems in homotopy theory is the computation of stable homotopy groups of spheres. While simple to define, they have proved to be extremely difficult to compute.

The standard way to compute stable homotopy groups is the Adams spectral sequence \cite{structure-applications}, which seeks to compute the stable homotopy groups of a finite spectrum $X$ from its cohomology
\[
  H^* X = \bigoplus_k \pi_{-k} F(X, \HFp).
\]
The crucial observation is that $H^*X$ is not just a group, but supports the action of cohomology operations. To capture this action, we define the algebra of all stable cohomology operations
\[
  \A = \bigoplus_k \pi_{-k} \End(\F_p).
\]
This algebra is known as the Steenrod algebra, and can be described explicitly in terms of generators and relations. Then $H^*X$ is naturally a module over $\A$, and the Adams spectral sequence takes the form
\[
  E^{s, t}_2 = \Ext_\A^{s, t}(H^*X, \F_p) \Rightarrow \pi_{t - s} X^\wedge_p.
\]

Unfortunately, even when $X$ is the sphere, this spectral sequence is highly non-trivial --- the $E_2$ page does not admit a simple description, and the differentials are hard to compute.

In practice, the first problem does not present a huge obstacle. Using a computer, one can iteratively construct a minimal free $\A$-resolution of $H^*X$ in a fairly efficient manner. This not only lets us read off the $\Ext$ groups; it also lets us compute the composition product of $\Ext$. Since the Adams spectral sequence is multiplicative, this lets us apply the Leibniz rule effectively, which massively simplifies the work involved in computing differentials. Similarly, we can compute Massey products and apply Moss' convergence theorem \cite{moss} to obtain new differentials.

Equipped with a (practically) full description of the Adams $E_2$ page, much subsequent work was focused on developing techniques to compute differentials by hand. In \cite{baues-e3}, Baues and Jibladze tackled the Adams spectral sequence from a different approach --- they discovered an algorithm that computes all $d_2$ differentials in the Adams spectral sequence of the sphere, thereby obtaining a computer description of the $E_3$ page.

The key insight of their algorithm is that just as the Adams $E_2$ page is controlled by the Steenrod algebra, the Adams $E_3$ page is controlled by the secondary Steenrod algebra, the algebra of all secondary cohomology operations.

Recall that secondary cohomology operations are defined by relations between cohomology operations. For example, the relation $\beta \beta = 0$ gives rise to a secondary cohomology operation $\Phi$, which is defined on elements $x \in H^* X$ such that $\beta x = 0$. To construct the action, represent $x$ as a map $\tilde{x} \colon X \to \HFp$. We then have a sequence
\[
  \begin{tikzcd}
    X \ar[r, "\tilde{x}"] & \HFp \ar[r, "\beta"] & \Sigma \HFp \ar[r, "\beta"] & \Sigma^2 \HFp
  \end{tikzcd}
\]
where any successive composition is trivial. The secondary cohomology operation is then defined as the Toda bracket
\[
  \Phi x = \langle \beta, \beta, \tilde{x} \rangle.
\]
Just as $\beta$ detects $p$, this cohomology operation $\Phi$ detects $p^2$. For example, it acts non-trivially on the cohomology of $\Z/p^2$.

Thus, to construct the secondary Steenrod algebra, we need to know not only all cohomology operations, but also homotopies between them. This suggests the definition
\[
  \A^{(2)} = \bigoplus_k \tau_{[0, 1]} \Sigma^k \End(\F_p).
\]
We similarly define the secondary cohomology functor by
\[
  \H^{(2)} X = \bigoplus_k \tau_{[0, 1]} \Sigma^k F(X, \F_p).
\]
One then sees that the action of $\A^{(2)}$ on $\H^{(2)}X$ lets us recover all secondary cohomology operations acting on $H^*X$.

While the Steenrod algebra is an actual algebra, $\A^{(2)}$ is \emph{a priori} only a graded $\E_1$-ring. Nevertheless, in \cite{baues-book}, Baues showed that $\A^{(2)}$ is in fact a differential graded algebra over $\Z/p^2$, and explicitly computed this differential graded algebra.

Equipped with this computation, Baues and Jibladze showed that we can compute all Adams $d_2$ differentials of a spectrum $X$ algorithmically given $\H^{(2)} X$. Together with a computation of $\H^{(2)}\S$, they were able to implement an algorithm to compute all $d_2$ differentials for the sphere. (Unfortunately, their implementation only managed to reach $t = 40$, and this theory was considered impractical by most of the broader computational homotopy theory community. This is, in fact, not true, as our implementation shows.)

The main thesis of this paper is that knowledge of the secondary Steenrod algebra in fact gives us ``full control'' of the Adams $E_3$ page, not just the $E_3$ page as a group. While the $E_3$ page of the Adams spectral sequence inherits a multiplication from the $E_2$ page, this is not the full picture; the ``$E_3$ page product'' ought to know about products up to one filtration higher. For example, it should be able to detect hidden extensions that jump by one filtration, as well as relations of the form
\[
  \nu^3 = \eta^2 \sigma + \eta \epsilon.
\]

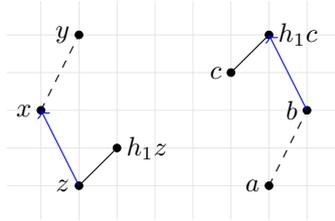
\begin{figure}[ht]
  \centering
  \begin{tikzpicture}[scale = 0.5]
    \draw [opacity = 0.1] (-1.9, -0.9) grid (6.9, 4.9);

    \node [left] at (-1, 2) {$x$};
    \node [left] at (0, 4) {$y$};
    \node [left] at (0, 0) {$z$};
    \node [right] at (1, 1) {$h_1 z$};

    \draw (0, 0) -- (1, 1);
    \draw [dashed] (-1, 2) -- (0, 4);

    \foreach \x/\y in {-1/2, 0/4, 0/0, 1/1} {
      \draw [fill] (\x, \y) circle (0.1);
    }
    \draw [blue, ->] (0, 0) -- (-1, 2);

    \node [left] at (4, 3) {$c$};
    \node [right] at (5, 4) {$h_1 c$};
    \node [left] at (5, 0) {$a$};
    \node [left] at (6, 2) {$b$};

    \draw (4, 3) -- (5, 4);
    \draw [dashed] (5, 0) -- (6, 2);

    \foreach \x/\y in {4/3, 5/4, 5/0, 6/2} {
      \draw [fill] (\x, \y) circle (0.1);
    }

    \draw [blue, ->] (6, 2) -- (5, 4);
  \end{tikzpicture}
  \caption{An example Adams chart with hidden extensions}\label{figure:fake-adams}
\end{figure}

This knowledge is extremely useful for computing the Adams spectral sequence. Consider a hypothetical Adams chart as in \Cref{figure:fake-adams}.

In this diagram, there are hidden $\eta$ extensions from $x$ to $y$ and from $a$ to $b$. Using a generalized version of the Leibniz rule, we can deduce that
\[
  d_3(h_1 z) = \text{``}\eta x\text{''} = y.
\]
Similarly, we can divide the differential $d_2(b) = h_1 c$ along $\eta$ to learn that $d_3(a) = c$. Crucially, this lets us relate differentials of different lengths, and in particular differentials on pages beyond the $E_2$ page.

\subsection*{Main results}
The paper has three main results, and is divided into three parts accordingly.

\subsubsection*{\ref{part:n-ary}\quad The \texorpdfstring{$n$}{n}-ary Steenrod algebra}
Our first result is to formalize the relationship between the secondary Steenrod algebra and the Adams $E_3$ page using the language of synthetic spectra \cite{synthetic} \cite[Appendix A]{manifold-synthetic}. Recall that for any Adams type spectrum $E$, the category $\Syn_E$ of $E$-based synthetic spectra is a symmetric monoidal category that interpolates between $\Comod_{E_*E}$ and $\Sp$. Specifically, there is an endomorphism $\tau$ of the unit $\S$ such that
\[
  \Mod_{C\tau} \cong \Comod_{E_*E},\quad \Mod_{\tau^{-1}\S} \cong \Sp.
\]
Further, there is a synthetic analogue functor $\nu \colon \Sp \to \Syn_E$ such that
\[
  C\tau \otimes \nu X \cong E_* X,\quad \tau^{-1} \nu X \cong X
\]
under the respective isomorphisms, and the $\tau$-Bockstein spectral sequence
\[
  \pi_{*, *} C\tau \otimes \nu X = \Ext_{E_*E} (E_*, E_* X) \Rightarrow \pi_{*, *} \tau^{-1} \nu X = \pi_* X
\]
of $\nu X$ is exactly the $E$-based Adams spectral sequence for $X$, at least up to a sign.

In the language of synthetic spectra, $\Mod_{C\tau}$ is the category controlling the Adams $E_2$ page. Similarly, $\Mod_{C\tau^n}$ fully captures information about the Adams $E_{n + 1}$ page. In particular, the ``$E_{n + 1}$ page product'' we alluded to is simply the composition product in $\Mod_{C\tau^n}$.

Using this, we reinterpret and extend Baues and Jibladze's result as

\begin{thm}
  Define the $n$-ary Steenrod algebra by
  \[
    \A^{(n)} = \bigoplus_k \tau_{[0, n)} \Sigma^k \End(\HFp).
  \]
  and let $E = \HFp$. Then there is a cocontinuous functor
  \[
    \H^{(n)} \colon \Mod_{C\tau^n} \to \Mod_{\A^{(n)}}^\op
  \]
  that is fully faithful when restricted to the full stable subcategory generated by objects of the form $C\tau^n \otimes \nu X$ where $X$ is a finite type spectrum.
\end{thm}

We will prove this in \Cref{part:n-ary} of the paper, as well as compatibility results between different $n$'s. Of course, we are mostly interested in the $n = 2$ case; the higher $n$ result is not practically useful without a computation of the $n$-ary Steenrod algebra itself.

\subsubsection*{\ref{part:secondary}\quad Computing \texorpdfstring{$E_3$}{E3} page data}
In \Cref{part:secondary}, we specialize to the case $n = 2$, where the secondary Steenrod algebra $\A^{(2)}$ was explicitly computed by Baues. Since $\A^{(2)}$ is a differential graded algebra over $\Z/p^2$, standard homological algebra gives us a model category presentation of $\Mod_{\A^{(2)}}$, which we can use to perform computations in $\Mod_{C\tau^2}$.

After describing the explicit algorithms to perform these computations, we implement them at the prime $2$ and compute the following data up to the $140$\textsuperscript{th} stem:
\begin{enumerate}
  \item all $d_2$ differentials;
  \item all $\Mod_{C\tau^2}$ products with $E_3$ page indecomposables up to the 39\textsuperscript{th} stem; and
  \item select $\Mod_{C\tau^2}$ Massey products, including the Adams periodicity operator.
\end{enumerate}

The primary purpose of this part is to document the details of the algorithm itself, and is largely aimed towards an audience interested in implementing the algorithm. The reader is encouraged to skip this part entirely if they are only interested in the mathematical underpinnings of the algorithm (and are satisfied with ``we have a model category presentation so we can compute anything''). Those who are interested in using the results to perform Adams spectral sequence calculations should read \Cref{section:data}, which explains how to retrieve and interpret our generated data.

\subsubsection*{\ref{part:computation}\quad Computing Adams differentials}
Finally, in \Cref{part:computation}, we use the computer generated data to compute new Adams differentials. We first formally define our notion of a hidden extension and prove a generalized version of the Leibniz rule. Using this, we resolve various unknown differentials in \cite{more-stable-stems}. In particular, we resolve all remaining unknown $d_2$, $d_3$, $d_4$ and $d_5$ differentials up to the 95\textsuperscript{th} stem.

\subsection*{Acknowledgements}

I would like to thank Christian Nassau, Dan Isaksen, Haynes Miller, John Rognes, Martin Frankland, Mike Hopkins, Piotr Pstr\k{a}gowski, Robert Bruner and Robert Burklund for helpful discussions related to this paper.

\part{The \texorpdfstring{$n$}{n}-ary Steenrod algebra}\label{part:n-ary}
\section{Overview}
The goal of this \namecref{part:n-ary} is to construct the comparison functor
\[
  \H^{(n)} \colon \Mod_{C\tau^n} \to \Mod_{\A^{(n)}}^\op
\]
and show that it is an equivalence on finite type objects.

We begin by introducing the category $\Mod_{C\tau^n}$ in \Cref{section:ctaun}. After proving basic categorical properties of the category, we move on to study duals in this category. The goal is to show that despite not being dualizable, the object $C\tau^m \in \Mod_{C\tau^n}$ for $m < n$ still behaves as if it were dualizable under many circumstances. For example, the natural map to its double dual is an equivalence.

We next introduce the $n$-ary Steenrod algebra in \Cref{section:steenrod}. After constructing $\A^{(n)}$ as a graded $\E_1$-ring, we show that it in fact strictifies to an algebra over a suitable quotient of the sphere. In the case $n = 2$, this recovers Baues' result that $\A^{(2)}$ is a differential graded algebra over $\Z/p^2$. We end by computing the secondary $\A(0)$ as a primer to the full secondary Steenrod algebra introduced in \Cref{section:nassau}.

In \Cref{section:comparison}, we construct the comparison functor $\H^{(n)}$, show that it is an equivalence on finite type objects, and prove naturality in $n$.

Our main result implies that the composition product in $\Mod_{C\tau^n}$ can be computed as the composition product in $\Mod_{\A^{(n)}}$. However, this is not quite true for the composition of bigraded mapping groups $[\Sigma^{a, b} X, Y]$; they only agree up to a sign! This is the infamous discrepancy between the product in the Adams $E_2$ page and the product in $\Ext$ \cite[p. 196]{structure-applications}. We will understand this in \Cref{section:bigraded} by carefully keeping track of the coherence data defining locally bigraded categories. Note that this is important even when $p = 2$, since we are now working over $\Z/4$, not $\Z/2$!

\subsection*{Conventions}
\begin{notation}
  If $\mathcal{C}$ is a spectrally enriched category and $X, Y \in \mathcal{C}$, we use $\mathcal{C}(X, Y)$ to denote the mapping space and $F_{\mathcal{C}}(X, Y)$ for the mapping spectrum. Thus, $\mathcal{C}(X, Y) = \Omega^\infty F_{\mathcal{C}}(X, Y)$. We will write $F_\Sp$ as $F$.

  If $\mathcal{C}$ is presentably symmetric monoidal, we write $\underline{\mathcal{C}}(X, Y)$ for the internal Hom.
\end{notation}

\begin{notation}
  We always use $DX$ to mean the weak dual of $X$ in the appropriate category. That is, $DX = \underline{\C}(X, \mathbf{1})$.
\end{notation}

\begin{notation}
  We write $\nu_n \colon \Sp \to \Mod_{C\tau^n}$ for the functor $X \mapsto C\tau^n \otimes \nu X$.
\end{notation}

\begin{notation}
  We define bigraded suspension in $\Syn_E$ by
  \[
    (\Sigma^{a, b} X)(P) = \Sigma^{-b} X(\Sigma^{-a - b} P).
  \]
  In particular, categorical suspension is $\Sigma^{1, -1}$, while $\nu \Sigma = \Sigma^{1, 0} \nu$. This is chosen to be compatible with the Adams spectral sequence, and differs from \cite{synthetic}.
\end{notation}

\section{Modules over \texorpdfstring{$C\tau^n$}{Ctaun}}\label{section:ctaun}
The goal of this \namecref{section:ctaun} is to construct and understand the category $\Mod_{C\tau^n}$.

Recall that the Adams $d_m$ differential is given by the connecting homomorphism of
\[
  \Sigma^{0, -m + 1} C\tau \to C\tau^m \to C\tau^{m - 1},
\]
which can be lifted to $\Mod_{C\tau^n}$ if $m \leq n$. In our future arguments, we would like to manipulate $C\tau^m$ as if it were dualizable. However, if $m < n$, then this is not true. Nevertheless, $C\tau^m$ is ``finite enough'' that in the situations of interest, it behaves as if it were dualizable. Specifically, we shall show that
\begin{thm}\label{thm:almost-compact}
  Let $X \in \Mod_{C\tau^n}$ be such that the underlying object in $\Syn_E$ is dualizable. Then
  \begin{enumerate}
    \item If $Y \in \Sp$ is a filtered colimit of objects in $\Sp_E^{fp}$, then the map
      \[
        DX \otimes \nu_n Y \to \underline{\Mod_{C\tau^n}}(X, \nu_n Y)
      \]
      is an equivalence.
    \item The map $X \to DDX$ is an equivalence.
  \end{enumerate}
\end{thm}
We will prove these in \Cref{thm:almost-dualizable,thm:double-dual} respectively.

\subsection{The category \texorpdfstring{$\Mod_{C\tau^n}$}{ModCtaun}}
To define $\Mod_{C\tau^n}$ at all, we have to construct $C\tau^n$ as an $\E_\infty$-ring, which does not follow from the definition as the cofiber as $\tau^n$. To make it a ring, we need an alternative description of $C\tau^n$.

Recall that $\Syn_E$ comes with a natural $t$-structure compatible with the symmetric monoidal structure \cite[Propositions 2.16, 2.29]{synthetic}. By \cite[Lemma 4.29]{synthetic}, we can write
\[
  C\tau^n = \tau_{<n} \S.
\]
This immediately gives
\begin{cor}\pushQED{\qed}
  There is a sequence of $\E_\infty$-rings
  \[
    \S \to \cdots \to C\tau^n \to \cdots \to C\tau^3 \to C\tau^2 \to C\tau.\qedhere
  \]
\end{cor}

This allows us to define the categories $\Mod_{C\tau^n}$, which are symmetric monoidal. As expected, these come with a natural $t$-structure.

\begin{lemma}
  Let $(\Mod_{C\tau^n})_{\geq 0}$ and $(\Mod_{C\tau^n})_{\leq 0}$ be the full subcategories of $\Mod_{C\tau^n}$ consisting of modules whose underlying object in $\Syn_E$ is connective and co-connective respectively. Then these form a right-complete $t$-structure compatible with filtered colimits and the symmetric monoidal structure.
\end{lemma}

\begin{proof}
  By \cite[Proposition 1.4.4.11]{ha}, there is a $t$-structure whose connective part is generated by $\{\nu_n P\}_{P \in \Sp_E^{fp}}$, and standard arguments (e.g.\ \cite[Lemma 5.3.2.12.3]{ha}) show that the connective part is $(\Mod_{C\tau^n})_{\geq 0}$. It follows from the adjunction that the co-connective objects are those whose underlying object is co-connective.

  The right-completeness and compatibility with filtered colimits follow from the same properties of the $t$-structure on $\Syn_E$. Compatibility with the symmetric monoidal structure follows from the bar construction model of the tensor product.
\end{proof}

Our original motivation was to study the cofiber sequence
\[
  \Sigma^{0, -m+1} C\tau \to C\tau^m \to C\tau^{m - 1}
\]
whose connecting map is the Adams differential. This cofiber sequence is easy to construct in $\Syn_E$, and we can lift it uniquely to one in $\Mod_{C\tau^n}$ by virtue of

\begin{lemma}\label{lemma:unique-lift}
  For any $k \in \Z$, any diagram in $(\Syn_E)_{[k, k + n)}$ lifts uniquely to a diagram in $\Mod_{C\tau^n}$.
\end{lemma}

\begin{proof}
  It suffices to prove this for the $k = 0$ case, and we have to show that $(\Mod_{\mathbf{1}_{<n}})_{[0, n)} \to (\Syn_E)_{[0, n)}$ is an equivalence. This follows from the general fact that given a compatible localization functor on a symmetric monoidal category, the category of local objects is equivalent to the category of local modules over the localization of the unit.
\end{proof}

\begin{cor}\label{cor:lift-cofib}
  For $k \leq m \leq n$ and $X \in \Sp$, the object $C\tau^m \otimes \nu X$ has a unique $C\tau^n$-structure. Further, the cofiber sequence
  \[
    \Sigma^{0, -k} C\tau^{m - k} \otimes \nu X \to C\tau^m \otimes \nu X \to C\tau^k \otimes \nu X
  \]
  has a unique lift to $\Mod_{C\tau^n}$.\fakeqed
\end{cor}

\subsection{Almost dualizable objects}
When $m < n$, the object $C\tau^m \in \Mod_{C\tau^n}$ is not dualizable. However, it does have the redeeming quality of being \emph{almost} dualizable.
\begin{defi}
  Let $\D$ be a presentably stable symmetric monoidal $\infty$-category with a compatible $t$-structure. We say $X \in \D$ is \emph{almost dualizable} if we can write $X = \colim X_k$ where
  \begin{enumerate}
    \item $X_k$ is dualizable; and
    \item $X_k \to X$ is a $k$-equivalence (i.e.\ it is an equivalence after $\tau_{\leq k}$).
  \end{enumerate}
\end{defi}

\begin{remark}
  In favorable circumstances, one can show that this agrees with the notion of almost compactness of \cite[Definition 7.2.4.8]{ha}.
\end{remark}

\begin{eg}
  A spectrum is almost dualizable iff it is finite type, i.e.\ it is bounded below and $H_*(X; \Z)$ is finite dimensional in each degree.
\end{eg}

\begin{lemma}
  Let $X \in \Mod_{C\tau^n}$ be almost dualizable and $Y \in \Sp$. Suppose $Y$ can be written as a filtered colimit of objects in $\Sp_E^{fp}$. Then the natural map
  \[
    DX \otimes \nu_n Y \to \underline{\Mod_{C\tau^n}}(X, \nu_n Y)
  \]
  is an equivalence.
\end{lemma}

\begin{proof}
  First observe that if $Y$ is in fact in $\Sp_E^{fp}$, then $\nu_n Y$ is dualizable, and the result holds unconditionally for all $X$.

  Let $X = \colim X_k$ as in the definition of almost dualizable. Let $X^k$ be the cofiber of $X_k \to X$. Then $X^k$ is $k$-connected. By \cite[Lemma 4.29]{synthetic}, we know $\nu_n Y$ is $n$-coconnected.

  We can write our map as
  \[
    DX \otimes \nu_n Y \to \lim \underline{\Mod_{C\tau^n}}(X_k, \nu_n Y) = \lim (DX_k \otimes \nu_n Y),
  \]
  whose fiber is $\lim (DX^k \otimes \nu_n Y)$.

  By right-completeness, it suffices to show that $DX^k \otimes \nu_n Y$ is $(n - k)$-coconnected. Since the $t$-structure is compatible with filtered colimits and $\nu_n$ preserves filtered colimits, we may assume $Y \in \Sp_E^{fp}$. Then $DX^k \otimes \nu_n Y = \underline{\Mod_{C\tau^n}}(DX^k, \nu_n Y)$, and the result follows.
\end{proof}

\begin{thm}\label{thm:almost-dualizable}
  Let $X \in \Mod_{C\tau^n}$. If the underlying object of $X$ is dualizable, then $X$ is almost dualizable.
\end{thm}

\begin{proof}
  By shifting $X$, we may assume that $X$ is connective. Let $X_\bullet$ be the bar construction on $X$ as a $C\tau^n$-module. Then we can write
  \[
    X = \colim X_\bullet = \colim_m \left(\colim_{\Delta^\op_{< m}} X_\bullet\right).
  \]
  Since $X_\bullet = (C\tau^n)^\bullet \otimes X$ is free on a dualizable object, it is dualizable. Further, when $k > m$, the cofiber of
  \[
    \colim_{\Delta^\op_{< m}} X_\bullet \to \colim_{\Delta^\op_{< k}} X_\bullet
  \]
  is $m$-connected \cite[Proposition 1.2.4.5.4]{ha}, so $\colim_{\Delta^\op_{< m}}X_\bullet \to X$ is an $m$-equivalence.
\end{proof}

\subsection{Weak duals in \texorpdfstring{$\Mod_{C\tau^n}$}{Ctaun}}
Finally, we compute the weak dual of $C\tau^m \in \Mod_{C\tau^n}$, and show that the natural map $C\tau^m \to DDC\tau^m$ is an equivalence. We begin by computing the (strong) dual in $\Syn_E$.

\begin{lemma}
  In $\Syn_E$, the dual of $\tau \colon \Sigma^{0, -1} \S \to \S$ is $\Sigma^{0, 1} \tau \colon \S \to \Sigma^{0, 1} \S$.

  Thus, we have
  \[
    D C\tau^n = \Sigma^{-1, n + 1} C\tau^n.
  \]
\end{lemma}

\begin{proof}
  For the first part, the map $\tau$ is constructed by starting with the following diagram in $\Sp_E^{fp}$:
  \[
    \begin{tikzcd}
      \S^{-1}  \ar[r] \ar[d] & * \ar[d] \\
      * \ar[r] & \S
    \end{tikzcd}
  \]
  applying $\nu$ to get
  \[
    \begin{tikzcd}
      \S^{-1, 0} \ar[r] \ar[d] & * \ar[d] \\
      * \ar[r] & \S
    \end{tikzcd}
  \]
  and then taking the induced map $\S^{0, -1} = \Sigma \S^{-1, 0} \to \S$. Since $\nu$ is symmetric monoidal and the dual of the first diagram is a suspension of the original diagram, the result follows.

  The second part follows immediately from the first.
\end{proof}

If $X \in \Mod_{C\tau^n}$, then its weak dual is defined to be $\underline{\Mod_{C\tau^n}}(X, C\tau^n)$. The key to understanding the weak dual is the observation that $C\tau^n$ is not just a free $C\tau^n$-module, but a cofree one as well.

Let $F \dashv U \colon \Mod_{C\tau^n} \rightleftharpoons \Syn_E$ be the free-forgetful adjunction. Since $U$ preserves all colimits, it has a right adjoint $C$.
\begin{lemma}
  We have
  \[
    C \S = \Sigma^{-1, n + 1} C\tau^n.
  \]
\end{lemma}

\begin{proof}
  Since $UC$ is right adjoint to $UF$, we find that
  \[
    UC(X) = DC\tau^n \otimes X.
  \]
  So the result follows from \Cref{cor:lift-cofib}.
\end{proof}

\begin{remark}
  In fact, one can show that $C \cong \Sigma^{-1, n + 1} F$.
\end{remark}

\begin{cor}
  There is a natural equivalence of functors
  \[
    UD \cong \Sigma^{1, -n - 1} DU.\qedhere
  \]
\end{cor}

\begin{proof}
  This follows from the more general relation
  \[
    U \underline{\Mod_{C\tau^n}}(X, CY) = \underline{\Syn_E}(UX, Y),
  \]
  whose proof is formal abstract nonsense using the projection formula
  \[
    U(FZ \otimes X) \cong Z \otimes UX.\qedhere
  \]
\end{proof}

\begin{cor}
  In $\Mod_{C\tau^n}$, if $m \leq n$, then $DC\tau^m \cong \Sigma^{0, m - n} C\tau^m$.\fakeqed
\end{cor}

\begin{thm}\label{thm:double-dual}
  Let $X \in \Mod_{C\tau^n}$ be such that $UX$ is dualizable. Then the natural map $X \to DD X$ is an equivalence.
\end{thm}

\begin{proof}
  Since the equivalence $UD = \Sigma^{1, -n - 1}DU$ preserves the co-evaluation map, this follows from the conservativity of $U$.
\end{proof}

\section{The \texorpdfstring{$n$}{n}-ary Steenrod algebra}\label{section:steenrod}
\subsection{Constructing the \texorpdfstring{$n$}{n}-ary Steenrod algebra}\label{section:n-ary-algebra}
Let $n \in [0, \infty]$. Informally, we can define the $n$-ary Steenrod algebra as
\[
  \A^{(n)} = \bigoplus_k \tau_{[0, n)} \Sigma^k \End(\HFp).
\]

While this is easy to write down as a graded spectrum, the $\E_1$-ring structure requires performing a categorical dance.
\begin{defi}
  The category of graded spectra is $\Sp^\gr = \Sp^\Z$, where $\Z$ is viewed as a discrete abelian group. This is a symmetric monoidal category under Day convolution.

  We give this a $t$-structure by declaring the connective part to be $\Sp_{\geq 0}^\Z$.
\end{defi}

\begin{defi}
  We define the bigraded spheres $\S^{a, b} \in \Sp^\gr$ to be $\S^a$ in degree $b$ and $0$ elsewhere. We define $[k] \colon \Sp^\gr \to \Sp^\gr$ to be $\S^{0, k} \otimes (-)$.
\end{defi}

The first step in constructing $\A^{(n)}$ as an $\E_1$-ring in $\Sp^\Z$ is to construct the $\E_1$-ring that is $\End(\HFp)$ in every degree. This follows from the functoriality of Day convolution.

\begin{lemma}[{\cite[Corollary 3.8]{day-functorial}}]
  Let $f \colon \C \to \C'$ be symmetric monoidal and $\mathcal{D}$ a presentably symmetric monoidal category. Then
  \begin{enumerate}
    \item $f^* \colon \mathcal{D}^{\C'} \to \mathcal{D}^\C$ is lax symmetric monoidal; and
    \item $f_* \colon \mathcal{D}^{\C} \to \mathcal{D}^{\C'}$ is symmetric monoidal,
  \end{enumerate}
  where $f^*$ is the restriction functor and $f_*$ is the left adjoint to $f^*$.\fakeqed
\end{lemma}

In our case, we have symmetric monoidal functors
\[
  \{0\} \overset\iota\hookrightarrow \Z \overset{\Delta}\twoheadrightarrow \{0\}
\]
which result in four functors between $\Sp^\gr$ and $\Sp$:
\begin{itemize}
  \item $\iota^* X_\bullet = X_0$.
  \item $(\iota_* X)_0 = X$ and vanishes in non-zero degrees.

  \item $(\Delta^* X)_n = X$ for all $n$.
  \item $\Delta_* X_\bullet = \bigoplus_n X_n$.
\end{itemize}

Then $\Delta^* \End(\HFp)$ is the $\E_1$-ring that is $\End(\HFp)$ in every degree. We next need to apply the shift $\Sigma^k$.

\begin{lemma}
  There is a cocontinuous $\E_1$-monoidal functor $\Phi \colon \Sp^\gr \to \Sp^\gr$ that sends $\{X_k\}$ to $\{\Sigma^k X_k\}$.
\end{lemma}

\begin{proof}
  By \cite[Proposition 4.8.1.10]{ha}, such a functor is equivalent to an $\E_1$-monoidal functor $\Z \to \Sp^\gr$, which we choose to send $k$ to $\S^{k, k}$. Then $\Phi$ is the unique cocontinuous functor extending this, and must be of the given form.
\end{proof}

\begin{defi}
  We define $\Phi^{(n)} \colon \Sp \to \Sp^\gr$ by $\Phi^{(n)} = \tau_{[0, n)} \Phi \Delta^*$.
\end{defi}
This is a lax $\E_1$-monoidal functor and $\Phi^{(n)} \Sigma X = (\Phi^{(n)} X)[-1]$.

\begin{defi}
  The $n$-ary Steenrod algebra is the $\E_1$-ring in $\Sp^\gr$ given by
  \[
    \A^{(n)} = \Phi^{(n)} \End(\HFp).
  \]
  The $n$-ary cohomology functor
  \[
    \H^{(n)} \colon \Sp \to \Mod_{\A^{(n)}}^\op
  \]
  is given by
  \[
    \H^{(n)}(X) = \Phi^{(n)} F(X, \HFp).
  \]
\end{defi}

The following lemmas are immediate from definition:
\begin{lemma}
  $\H^{(n)}$ sends sums to products and $\H^{(n)}(\Sigma X) = \H^{(n)}(X)[1]$.\fakeqed
\end{lemma}

\begin{lemma}\label{lemma:hn-homotopy}\pushQED{\qed}
  Let $\tau \in \pi_{1, 1} \A^{(n)}$ be the identity map in $\Sigma \End(\HFp)$. Then there is an isomorphism of rings
  \[
    \pi_{*, *} \A^{(n)} = \A [\tau] / \tau^n,
  \]
  where $\A$ is the ordinary Steenrod algebra. Further, as an $\pi_{*, *} \A^{(n)}$-module, we have
  \[
    \pi_{*, *} \H^{(n)} (X) = H^*(X) [\tau] / \tau^n.
  \]
  Thus, we have
  \[
    \H^{(n - 1)} (X) = \A^{(n - 1)} \otimes_{\A^{(n)}} \H^{(n)} (X).\qedhere
  \]
\end{lemma}

\begin{remark}
  $\A^{(\infty)}$ is a \emph{shift algebra} in the sense of \cite{abstract-goerss-hopkins}, and $\H^{(\infty)} (X)$ is a periodic $\A^{(\infty)}$-module.

  This lets us apply the results of \cite[Section 4]{abstract-goerss-hopkins}. In particular, $\H^{(n)}(X)$ is a potential $(n - 1)$-stage, and there is an obstruction theory for the space of possible values of $\H^{(n)}(X)$ given $H^*(X)$. In many cases of interest (e.g.\ the sphere), this space is connected, so any $\A^{(n)}$-module with the right homotopy groups must be $\H^{(n)}(X)$.
\end{remark}

\begin{remark}
  It follows from the descriptions of the homotopy groups that if $A \to B \to C$ induces a short exact sequence on cohomology, then $\H^{(n)} (A) \to \H^{(n)} (B) \to \H^{(n)} (C)$ is a cofiber sequence.
\end{remark}

In order to connect $\Mod_{\A^{(n)}}$ to $\Mod_{C\tau^n}$, we will need an alternative description of $\Mod_{\A^{(n)}}$. Recall the $P_\Sigma$ construction from \cite[Definition 5.5.8.8]{htt}, and write $P_\Sigma^\Sp(\C)$ for the stabilization of $P_\Sigma(\C)$. By \cite[Remark C.1.5.9]{sag}, $P_\Sigma^\Sp(\C)$ is the full subcategory of $\Fun(\C^\op, \Sp)$ of (finite) product-preserving functors.

\begin{defi}
  We let $\Free_{\A^{(n)}}$ be the full subcategory of $\Mod_{\A^{(n)}}$ consisting of finite direct sums of modules of the form $\A^{(n)}[k]$.
\end{defi}

\begin{lemma}
  There is an equivalence of categories
  \[
    P_\Sigma^\Sp(\Free_{\A^{(n)}}) \cong \Mod_{\A^{(n)}}
  \]
  with inverse given by the spectral Yoneda embedding.
\end{lemma}

\begin{proof}
  The inclusion of $\Free_{\A^{(n)}}$ gives a cocontinuous map
  \[
    F\colon P_\Sigma(\Free_{\A^{(n)}}) \to \Mod_{\A^{(n)}}.
  \]
  We claim this is fully faithful with essential image given by $(\Mod_{\A^{(n)}})_{\geq 0}$. By \cite[Proposition 5.5.8.22]{htt}, we need to show that the inclusion of $\Free_{\A^{(n)}}$ is fully faithful with image given by compact projective generators. The first part is clear and the second follows from \cite[Corollary 7.1.4.14]{ha}.

  Since $\Mod_{\A^{(n)}}$ is the stabilization of its connective part, the stabilization $\tilde{F}$ of $F$ is an equivalence of categories.

  Let $G$ and $\tilde{G}$ be the right adjoints to $F$ and $\tilde{F}$ respectively. By combining \cite[Corollary 5.2.6.5, Proposition 5.5.8.10]{htt}, we know $G$ is given by the Yoneda embedding.

  For the spectral version, note that we must have $\Omega^\infty \tilde{G} = G$, since both are right adjoints to $F = \tilde{F} \Sigma^\infty_+$. By \cite[Corollary 1.4.2.23]{ha}, any two left exact functors $\Mod_{\A^{(n)}} \to P_\Sigma(\Free_{\A^{(n)}})$ that agree after applying $\Omega^\infty$ must in fact agree. So $\tilde{G}$ must be the spectral Yoneda embedding.
\end{proof}

\begin{remark}
  Given a presheaf $X \colon \Free_{\A^{(n)}}^\op \to \Sp$, the underlying $\A^{(n)}$-module of $X$ is given by
  \[
    X_k = X(\A^{(n)} [k]).
  \]
  If $\alpha \in \pi_{0, \ell} \A^{(n)} = \A_\ell$, then its action on $X$ is given by applying the presheaf to
  \[
    \alpha \colon \A^{(n)}[\ell] \to \A^{(n)}.
  \]
  This is all standard; the less-obvious part is the action of $\tau$. Informally, we expect this to be given by $X$ acting on $\tau \colon \Sigma \A^{(n)}[1] \to \A^{(n)}$. However, $\Sigma \A^{(n)}[1]$ is not an object in $\Free_{\A^{(n)}}$. Nevertheless, the map $\tau$ is represented by a commutative diagram
  \[
    \begin{tikzcd}
      \A^{(n)}[1] \ar[r] \ar[d] & * \ar[d] \\
      * \ar[r] & \A^{(n)}.
    \end{tikzcd}
  \]
  Applying $X$ to this diagram then gives a map $\Sigma X \to X[-1]$, which is the action of $\tau$.
\end{remark}

The category $\Free_{\A^{(n)}}$ in turn admits a more direct definition in terms of Eilenberg--Maclane spectra.
\begin{defi}
  Let $\M_{\HFp}$ be the full subcategory of $\Sp$ consisting of finite sums of shifts of $\HFp$.
\end{defi}

\begin{lemma}\label{cor:free-an}
  The $n$-ary cohomology functor gives an equivalence of categories
  \[
    h_n \M_{\HFp} \overset\sim\to \Free_{\A^{(n)}}^\op.
  \]
  Under the isomorphism $\Mod_{\A^{(n)}} \cong P_\Sigma^\Sp(h_n \M_{\HFp}^\op)$, the $n$-ary cohomology $\H^{(n)} X$ of a spectrum $X$ corresponds to the presheaf
  \[
    M \mapsto \tau_{[0, n)} F(X, M).
  \]
\end{lemma}

This is an immediate consequence of the following more general lemma:
\begin{lemma}
  Let $X \in \Sp$ and $M \in \M_{\HFp}$. Then the natural map
  \[
    \H^{(n)}\colon \Sp(X, M) \to \Mod_{\A^{(n)}}(\H^{(n)}(M), \H^{(n)}(X))
  \]
  is $n$-truncation.
\end{lemma}

\begin{proof}
  Since $\H^{(n)}$ preserves shifts and direct sums, we may assume $M = \A^{(n)}$. Then the right-hand side is
  \begin{multline*}
    \Mod_{\A^{(n)}}(\A^{(n)}, \H^{(n)}(X)) = \Sp^\gr(\iota_* \S, \H^{(n)}(X)) = \Sp(\S, \iota^* \H^{(n)}(X)) \\
    = \Sp(\S, \tau_{[0, n)} F(X, \HFp)) = \tau_{\leq n} \Sp(X, \HFp).
  \end{multline*}
\end{proof}

We end with a lemma on the naturality of this isomorphism.
\begin{lemma}\label{lemma:psigma-nat}
  Let $m < n$. Under the isomorphism $\Mod_{\A^{(n)}} \cong P_\Sigma^\Sp(h_n \M_{\HFp}^\op)$, the forgetful functor $\Mod_{\A^{(m)}} \to \Mod_{\A^{(n)}}$ corresponds to restriction along $h_n \M_\HFp^\op \to h_m \M_\HFp^\op$.
\end{lemma}

\begin{proof}
  It suffices to show that they have the same left adjoint. By construction, the left adjoint to restriction along $h_n \M_\HFp^\op \to h_m \M_\HFp^\op$ is the unique stable cocontinuous functor $P_\Sigma^\Sp(h_n \M_\HFp^\op) \to P_\Sigma^\Sp(h_m \M_\HFp^\op)$ that extends the map $h_n \M_\HFp^\op \to h_m \M_\HFp^\op$. Since $\A^{(m)} \otimes_{\A^{(n)}}(-)$ also fits this description, we are done.
\end{proof}

\begin{remark}
  The equivalence $\Mod_{\A^{(n)}} \cong P_\Sigma^\Sp(h_n \M_\HFp^\op)$ lets us directly construct the category of modules over $\A^{(n)}$ without constructing $\A^{(n)}$ itself. A version of this was studied by \cite{mapping-algebra} using the language of model categories. They were then able to prove directly that it encodes information about the Adams $E_{n + 1}$ page.
\end{remark}

\subsection{Strictifying the \texorpdfstring{$n$}{n}-ary Steenrod algebra}
\emph{A priori}, the algebra $\A^{(n)}$ is a ring over $\S$. In the $n = 1$ case, we know it is in fact a ring over $\F_p$, which is much easier to work with. In the $n = 2$ case, Baues \cite[Section 5]{baues-book} showed that $\A^{(2)}$ is a ring over $\Z/p^2$, which also allows us to employ homological algebra machinery to perform computations. In general, $\A^{(n)}$ is a ring over a suitable truncation of the sphere.

\begin{defi}
  For $E$ an Adams type homology theory and $n \geq 1$, define the truncation functor
  \[
    \tau^E_{<n} = (-)^E_{<n} \colon \Sp_{\geq 0} \to \Sp_{\geq 0}
  \]
  by
  \[
    \tau^E_{<n} X = X^E_{<n} = F_{\Syn_E}(\S, C\tau^n \otimes \nu X).
  \]
\end{defi}

\begin{lemma}
  The functor $X \mapsto X^E_{<n}$ is lax symmetric monoidal and natural in $n$. Moreover, there is a natural transformation of lax symmetric monoidal functors $X \to X^E_{<n}$. On homotopy groups, this kills elements whose image in the $E$-based Adams spectral sequence has $t$-coordinate at least $n$.\fakeqed
\end{lemma}

\begin{remark}
  The definition of $\tau^E_{<n}$ makes sense for non-connective spectra as well, but the effect on negative homotopy groups is more subtle.
\end{remark}

\begin{eg}
  $\S^{\HFp}_{<1} = \F_p$, $\S^{\HFp}_{<2} = \Z/p^2$ and $\pi_* \S^{\F_2}_{<3} = \Z/8 [\eta] / \eta^2$.
\end{eg}

By construction, we know that
\begin{lemma}
  $F(C\tau^n, -) \colon \Mod_{C\tau^n} \to \Sp$ lifts to a lax symmetric monoidal functor $\Mod_{C\tau^n} \to \Mod_{\S^{\HFp}_{<n}}$.\fakeqed
\end{lemma}

\begin{thm}
  $\A^{(n)}$ lifts to an $\E_1$-algebra in $\Mod_{\S^{\HFp}_{<n}}^\Z$.
\end{thm}

\begin{proof}
  By \cite[Construction C.17]{r-motivic}, there is a symmetric monoidal functor $\Z \to \Syn_{\HFp}$ that sends $n$ to $\S^{0, -n}$. Call this functor $\S^{0, -\bullet}$.

  Consider the composite
  \[
    \Sp \overset{\nu}\longrightarrow \Syn_{\HFp} \overset{\Delta^*}\longrightarrow \Syn_{\HFp}^\Z \overset{\otimes \S^{0, -\bullet}}\longrightarrow \Syn_{\HFp}^\Z \overset{C\tau^n \otimes }\longrightarrow \Mod_{C\tau^n}^\Z \overset{F(C\tau^n, -)}\longrightarrow \Sp^\Z\overset{\Phi}\longrightarrow \Sp^\Z.
  \]
  All the functors are lax symmetric monoidal, and the last two maps naturally lift to $\Mod_{\S_{<n}^\HFp}^\Z$. So it suffices to show that $\A^{(n)}$ is the image of $\End(\HFp)$ under this map.

  Let this composite be $F_{\nu \otimes C\tau^n}$. We consider the variations $F_Y$ and $F_\nu$ defined as follows:
  \begin{itemize}
    \item $F_\nu$ is obtained by replacing the fourth and fifth maps of $F_{\nu \otimes C\tau^n}$ with $\Syn_{\HFp}^\Z \overset{F(\S, -)}\longrightarrow \Sp^\Z$.
    \item $F_Y$ is obtained by replacing the first map of $F_\nu$ with $Y\colon \Sp \to \Syn_{\HFp}$ (recall that $Y = \tau^{-1} \nu$ is the spectral Yoneda embedding).
  \end{itemize}
  We then have natural transformations
  \[
    \begin{tikzcd}
      F_\nu \ar[r] \ar[d] & F_{\nu \otimes C\tau^n} \\
      F_Y
    \end{tikzcd}
  \]
  of lax symmetric monoidal functors. It suffices to show that
  \begin{enumerate}
    \item $F_Y \cong \Phi \Delta$;
    \item $F_\nu \End(\HFp) = \tau_{\geq 0} F_Y \End(\HFp)$; and
    \item $F_{\nu \otimes C\tau^n} \End(\HFp) = \tau_{< n} F_\nu \End(\HFp)$.
  \end{enumerate}
  The last two can be checked on homotopy groups. To prove the first, note that on $\tau$-invertible spectra, $F(\S, -)$ is canonically equivalent to $\tau^{-1} \colon \Syn_{\HFp} \to \Sp$ as a lax symmetric monoidal functor. Further, $\S^{0, -\bullet}$ is constructed so that after $\tau$-inversion, it is the symmetric monoidal functor that is constantly the unit. So we are done.
\end{proof}

\subsection{The secondary \texorpdfstring{$\A(0)$}{A(0)}}\label{section:a0}
Before we move on, it is prudent to give some intuition for what $\A^{(n)}$ looks like. In \Cref{section:nassau}, we are going to give a full description of $\A^{(2)}$. However, this description is fairly complex and it is easy to get lost in the details. To provide a simpler example, we instead look at the secondary $\A(0)$, defined by
\[
  \A(0)^{(2)} = \bigoplus_{k \in \Z} \tau_{[0, 1]} \Sigma^{k} \End_\Z(\F_p).
\]

To compute this, we use the following explicit presentation of $\F_p \in \Mod_\Z$:
\[
  \F_p =
  \left(
  \begin{tikzcd}
    \Z \{x_1\} \ar[d, "p"] \\ \Z \{x_0\}
  \end{tikzcd}
  \right).
\]
Then $\End_\Z(\F_p)$ is given by
\[
  \End_\Z(\F_p) = \F_p \otimes \F_p^* = \left(
    \begin{tikzcd}
      \Z \{x_1 \otimes x_0^*\} \ar[d, "{(p, -p)}"] \\
      \Z\{x_0 \otimes x_0^*, x_1 \otimes x_1^*\} \ar[d, "{(p, p)}"] \\
      \Z\{x_0 \otimes x_1^*\}.
    \end{tikzcd}
  \right).
\]
As is well-known, $\pi_* \End_Z(\F_p) = \F_p\{1, \beta\}$, with explicit representatives given by
\[
  1 = x_0 \otimes x_0^* - x_1 \otimes x_1^*,\quad \beta = x_0 \otimes x_1^*.
\]

By definition, $\A(0)^{(2)}$ is the sum of truncations
\[
  \begin{tikzcd}[column sep=0.5em, row sep=tiny]
    \color{gray} k = 0 & \color{gray} k = 1  & \color{gray} k = 2 \\
    \Z\{x_1 \otimes x_0^*\} \ar[dd, "p"] & \Z\{x_1 \otimes x_1^*\} \oplus \F_p\{x_0 \otimes x_0^* - x_1 \otimes x_1^*\}\ar[dd, "{(p, 0)}"] & \F_p\{x_0 \otimes x_1^*\} \ar[dd] \\
    \vphantom{x} \\
    \Z\{x_0 \otimes x_0^* - x_1 \otimes x_1^*\} & \Z\{x_0 \otimes x_1^*\} & 0
  \end{tikzcd}
\]
In $\A(0)^{(2)}$, we let $1$ and $\beta$ be the corresponding classes in cohomological degree $0$, and define the following classes in cohomological degree $1$:
\[
  \mu_0 = x_1 \otimes x_0^*,\quad \tau = x_0 \otimes x_0^* - x_1 \otimes x_1^*.
\]
Thus, $\mu_0$ is the null-homotopy of $p$, while $\tau$ detects the copy of $1$ in cohomological degree $1$. We can then write $\A(0)^{(2)}$ as the chain complex
\[
  \A(0)^{(2)} = \left(
    \begin{tikzcd}
      \Z\{\mu_0, \mu_0 \beta\} \oplus \F_p\{\tau, \tau \beta\} \ar[d, "d"] \\
      \Z \{1, \beta\}
    \end{tikzcd}
  \right),\quad\sidedeg
\]
As for the algebra structure, $\tau$ acts centrally, while we have the crucial relation
\[
  \beta \mu_0 = \mu_0 \beta + \tau.
\]
This relation encodes the fact that $\beta$ detects $p$.

We now have a full description of $\A(0)^{(2)}$. However, this is a differential graded algebra over $\Z$, instead of the promised $\Z/p^2$. To remedy this, observe that as a chain complex, $\A(0)^{(2)}$ is in fact equivalent to one over $\F_p$ --- we can simply quotient out the $\mu_0$ factors and end up with $\A(0)$ in cohomological degrees $0$ and $1$. However, this quotienting does not respect the algebra relation $\beta \mu_0 = \mu_0 \beta + \tau$. Nevertheless, since $p \tau = 0$, we can quotient out $p \mu_0$, and get our final presentation
\[
  \A(0)^{(2)} = \left(
    \begin{tikzcd}
      \F_p\{\mu_0, \mu_0 \beta\} \oplus \F_p\{\tau, \tau \beta\} \ar[d, "d"] \\
      \Z/p^2 \{1, \beta\}
    \end{tikzcd}
  \right),\quad\sidedeg
\]

Equipped with a presentation of $\A(0)^{(2)}$, we can now compute the secondary cohomology of various $\Z$-modules. The simplest $\Z$-module is, of course, $\Z$ itself. Tracing through the definitions gives the following presentation of the secondary cohomology of $\Z$:
\[
  \kk = \left(\begin{tikzcd}
      \F_p\{\mu_0\} \oplus \F_p\{\tau\} \ar[d, "d"] \\
      \Z/p^2
  \end{tikzcd}\right),\quad\sidedeg
\]
The only non-trivial $\A(0)^{(2)}$ action is given by
\[
  \beta \mu_0 = \tau.
\]
Alternatively, this can be described as $\A(0)^{(2)} / (\A(0)^{(2)} \beta)$.

\begin{remark}
  While there is a ring map $\kk \to \A(0)^{(2)}$, this does not map to the literal center of $\A(0)^{(2)}$. Instead, there is a chain homotopy between the left and right multiplication maps.
\end{remark}

More interestingly, we can look at two-cell complexes. The first example we can look at is $\Z/p$. Since this has cells in degrees $0$ and $1$, as a chain complex, we have
\[
  \H^{(2)} (\Z/p) = \kk \{a\} \oplus \kk \{b\},\quad |a| = 0,\quad |b| = 1.
\]
However, there is a non-trivial $\A(0)^{(2)}$ action given by $\beta a = b$. This distinguishes it from $\H^{(2)}(\Z \oplus \Z[1])$, which has the same underlying chain complex but with $\beta a = 0$. Of course, this difference already manifests itself on the level of ordinary cohomology, without having to go to the secondary level.

On the other hand, $\Z/p^2$ and $\Z \oplus \Z[1]$ \emph{do} have the same ordinary cohomology groups. We can compute that $\H^{(2)} (\Z/p^2)$ is again $\kk\{a\} \oplus \kk \{b\}$, but now the $\A(0)^{(2)}$ action is given by
\[
  \beta a = p b.
\]
Since $pb$ is null-homotopic, this is not visible on the level of ordinary cohomology. Instead, this is detected by the secondary cohomology operation associated to the equation $\beta \beta = 0$. Indeed, we can compute
\[
  \langle \beta, \beta, a\rangle = 0 \cdot a + \beta \cdot \mu_0 b = \tau b
\]
with no indeterminacy.

\section{The comparison functor}\label{section:comparison}
Set $E = \HFp$. In this \namecref{section:comparison}, we will construct the comparison functor
\[
  \H^{(n)} \colon \Mod_{C\tau^n} \to \Mod_{\A^{(n)}}^\op
\]
and show that it has the desired properties. The constructions will work when $n = \infty$ as well, in which case we set $\Mod_{C\tau^n}$ to be $\widehat{\Syn_\HFp}$, the category of hypercomplete synthetic spectra, and $C\tau^n \otimes (-)$ is the hypercompletion functor.

\subsection{Constructing the comparison functor}
The comparison functor will be a natural extension of the $n$-ary cohomology functor $\H^{(n)} \colon \Sp \to \Mod_{\A^{(n)}}^\op$ along $\nu_n \colon \Sp \to \Mod_{C\tau^n}$. To construct this functor, we use the isomorphism
\[
  \Mod_{\A^{(n)}} \cong P_\Sigma^\Sp(h_n \M_{\HFp}^\op).
\]
By \Cref{cor:free-an}, the $n$-ary cohomology functor can then be described by
\[
  \H^{(n)}(X)(M) = \tau_{[0, n)} F(X, M).
\]

To extend $\H^{(n)}$ along the map $\nu_n$, we need to express $\tau_{[0, n)} F(X, M)$ in terms of $\nu_n X$. This follows from the following lemma:
\begin{lemma}\label{lemma:truncation}
  Let $M$ be a homotopy ${\HFp}$-module. Then $\nu M = \tau_{\geq 0} F(-, M)$. Thus, for any $X \in \Sp$, we have
  \[
    F_{\Syn_{\HFp}}(\nu X, \nu M) = \tau_{\geq 0} F(X, M).
  \]
  Further, the functor $C\tau^n \otimes(-)$ exhibits $F_{\Mod_{C\tau^n}}(\nu_n X, \nu_n M)$ as the $n$-truncation of $F_{\Syn_{\HFp}}(\nu X, \nu M)$.
\end{lemma}

\begin{proof}
  By construction, $\nu M$ is the sheafification of $\tau_{\geq 0} F(-, M)$. Thus, we have to show that $\tau_{\geq 0} F(-, M)$ is already a sheaf. By \cite[Theorem 2.8]{synthetic}, we have to show that if $A \to B \to C$ is a cofiber sequence in $\Sp_{\HFp}^{fp}$ with the second map being an $(\HFp)_*$-surjection, then
  \[
    \tau_{\geq 0} F(C, M) \to \tau_{\geq 0} F(B, M) \to \tau_{\geq 0} F(A, M)
  \]
  is a fiber sequence.\footnote{\cite[Theorem 2.8]{synthetic} states this for presheaves of spaces instead of spectra. However, the proof reduces it to the case of spectra and proves it for spectra.} Since this is a fiber sequence before applying $\tau_{\geq 0}$, it suffices to show that
  \[
    [B, M] \to [A, M] \to 0
  \]
  is exact. Since ${\HFp}$ is Adams type, this is given by
  \[
    \Hom_{\F_p}((\HFp)_* B, M_*) \to \Hom_{\F_p}((\HFp)_* A, M_*) \to 0.
  \]
  Since $(\HFp)_* A \to (\HFp)_* B$ splits, the result follows.

  To prove the second part, note that $\nu$ preserves filtered colimits, so it suffices to prove this when $X$ is finite, which follows from Yoneda's lemma.

  The last part follows from the construction of $\tau$ (for $n = \infty$, use that $\nu M$ is already hypercomplete).
\end{proof}

\begin{cor}
  The composite
  \[
    \begin{tikzcd}[column sep=large]
      \M_{\HFp} \ar[r, hook] & \Sp \ar[r, "\nu"] & \Syn_{\HFp} \ar[r, "C\tau^n \otimes (-)"] & \Mod_{C\tau^n}
    \end{tikzcd}
  \]
  identifies the image with $h_n \M_{\HFp}$.
\end{cor}

This allows us to define the comparison functor as follows:
\begin{defi}
  We define $\H^{(n)}\colon \Mod_{C\tau^n} \to P_\Sigma^\Sp(h_n \M_\HFp^\op)^\op \cong \Mod_{\A^{(n)}}^\op$ by
  \[
    \H^{(n)}(X)(M) = F_{\Mod_{C\tau^n}}(X, \nu_n M).
  \]
\end{defi}
It is easy to see that $\H^{(n)}$ preserves the two suspension functors (but see \Cref{section:bigraded} for crucial details), and a little diagram chase shows that $\H^{(n)}$ sends $\tau$ to $\tau$.

For the rest of the \namecref{section:comparison}, we will use $\H^{(n)}$ to refer to this extended functor on $\Mod_{C\tau^n}$ instead of the $n$-ary cohomology functor.

\begin{remark}
  A useful property of the comparison functor $\H^{(n)}$ is that it is cocontinuous, unlike the $n$-ary cohomology functor. The trade-off is that $\nu_n$ is, of course, not cocontinuous.
\end{remark}

\subsection{Fully faithfulness of the comparison functor}
We shall show that $\H^{(n)}$ is fully faithful when restricted to the full subcategory of finite type objects.

Recall that a spectrum $X$ is finite type if it is bounded below and $H_*(X; \Z)$ is finite dimensional in each degree. In other words, it is a bounded below spectrum built with finitely many cells in each degree. We let $\Sp^{ft} \subseteq \Sp$ be the full subcategory of finite type spectra. By K\"unneth's formula, this is closed under tensor products.

\begin{defi}
  Let $\Mod_{C\tau^n}^{ft} \subseteq \Mod_{C\tau^n}$ be the full stable subcategory generated by $\{\nu_n P\}_{P \in \Sp^{ft}}$.
\end{defi}

\begin{thm}\label{thm:y-ff}
  $\H^{(n)}$ restricts to a fully faithful functor $\Mod_{C\tau^n}^{ft} \to \Mod_{\A^{(n)}}^\op$. In fact, for any $X \in \Mod_{C\tau^n}$ and $Y \in \Mod_{C\tau^n}^{ft}$, the map
  \[
    \H^{(n)} \colon \Mod_{C\tau^n}(X, Y) \to \Mod_{\A^{(n)}}(\H^{(n)}Y, \H^{(n)} X)
  \]
  is an equivalence.
\end{thm}

The proof requires some auxiliary lemmas, which we will prove after proving the main theorem.

\begin{proof}
  Let $\C \subseteq \Mod_{C\tau^n}$ be the full subcategory of $\Mod_{C\tau^n}$ consisting of spectra $Y$ such that for any $X \in \Mod_{C\tau^n}$, the map
  \[
    F_{\Mod_{C\tau^n}}(X, Y) \to F_{\Mod_{\A^{(n)}}}(\H^{(n)}Y, \H^{(n)} X)
  \]
  is an equivalence. Then $\C$ is stable and is closed under limits preserved by $\H^{(n)}$ (that is, limits that are sent to colimits in $\Mod_{\A^{(n)}}$).

  By the spectral Yoneda lemma, $\M_{\HFp} \subseteq \C$.

  Next, we show that if $P \in \Sp^{ft}$, then $\nu_n (\HFp \otimes P) \in \C$. Indeed, we can write
  \[
    \HFp \otimes P = \bigoplus \Sigma^{k_i} {\HFp}
  \]
  where $k_i \to \infty$. By \Cref{cor:sum-prod}, we have
  \[
    \nu_n (\HFp \otimes P) = \bigoplus \nu_n \Sigma^{k_i} \HFp = \prod \nu_n \Sigma^{k_i} \HFp.
  \]
  Since $\H^{(n)}$ is cocontinuous, it preserves direct sums. Thus, we have
  \[
    \H^{(n)} \prod \nu_n \Sigma^{k_i} \HFp = \H^{(n)} \bigoplus \nu_n \Sigma^{k_i} \HFp = \prod \H^{(n)} \nu_n \Sigma^{k_i} \HFp = \bigoplus \H^{(n)} \nu_n \Sigma^{k_i} \HFp.
  \]
  So the direct product is preserved by $\H^{(n)}$, and $\nu_n (\HFp \otimes P) \in \C$.

  Finally, if $P \in \Sp^{ft}$, then by \Cref{lemma:ass-converge}, its Adams spectral sequence converges. That is, we have
  \[
    \nu_n P = \varprojlim \nu_n (CB^\bullet(\HFp) \otimes P).
  \]
  By \Cref{lemma:y-tot}, this limit is preserved by $\H^{(n)}$. So $\nu_n P \in \C$.
\end{proof}

We now prove the various lemmas used in the proof.
\begin{lemma}
  Let $X$ be $k$-connective. Then $\pi_{a, b} \nu X = \pi_{a, b} \nu_n X = 0$ when $a < k$.
\end{lemma}

\begin{proof}
  We first prove the $\nu X$ version. By \cite[Theorem 4.58]{synthetic}, this is true for $b \leq 0$. For $b > 0$, the long exact sequence from \cite[Proposition 4.57]{synthetic} gives
  \[
    \Ext^{b - 1, a + b}_{\A_*}(\F_p, H_* X) \to \pi_{a, b - 1} \nu X \to \pi_{a, b} \nu X \to \Ext_{\A_*}^{b, a + b}(\F_p, H_* X).
  \]
  So the result follows from the vertical vanishing line of $\Ext$.

  As for the $\nu_n$ version, the $n < \infty$ case follows from the cofiber sequence
  \[
    \begin{tikzcd}
      \Sigma^{0, -n} \nu X \ar[r, "\tau^n"] & \nu X \ar[r] & C\tau^n \otimes \nu X.
    \end{tikzcd}
  \]
  When $n = \infty$, since $X$ is connective, we have $\nu_n X = \nu X^\wedge_p$. Since $X^\wedge_p$ is also $k$-connective, the result follows.
\end{proof}

\begin{cor}\label{cor:sum-prod}
  Let $\{X_i\}_{i \in \mathbb{N}}$ be a sequence of spectra such that $X_i$ is $k_i$-connective and $k_i \to \infty$. Then
  \[
    \bigoplus X_i = \prod X_i,\quad \nu_n \left(\bigoplus X_i\right) = \bigoplus \nu_n X_i = \prod \nu_n X_i.
  \]
\end{cor}

\begin{proof}
  The first part is standard. The equality $\nu_n \left(\bigoplus X_i\right) = \bigoplus \nu_n X_i$ follows from $\nu_n$ preserving finite coproducts and filtered colimits, hence infinite coproducts. To show that
  \[
    \bigoplus \nu_n X_i = \prod \nu_n X_i,
  \]
  we use the fact that $\Mod_{C\tau^n}$ is generated by shifts of $\{\nu_n P\}_{P \in \Sp_\HFp^{fp}}$ under colimits. So it suffices to show that
  \[
    \left[ \nu_n P, \bigoplus \nu_n X_i\right]^{*, *} =
    \left[ \nu_n P, \prod \nu_n X_i\right]^{*, *}
  \]
  Since $\nu_n P$ is compact, this is equivalent to showing that
  \[
    \bigoplus \left[ \nu_n P, \nu_n X_i\right]^{*, *} =
    \prod \left[ \nu_n P, \nu_n X_i\right]^{*, *}.
  \]
  Since $\nu_n P$ is dualizable and $\nu_n DP \otimes \nu_n X_i = \nu_n (DP \otimes X_i)$, we may assume that $P = \S$. So we have to show that
  \[
    \bigoplus \pi_{*, *} \nu_n X_i = \prod \pi_{*, *} \nu_n X_i
  \]
  This follows from the previous vanishing line.
\end{proof}

\begin{lemma}\label{lemma:ass-converge}
  If $X$ is any bounded below spectrum, then $\nu_n X$ is $\nu_n \HFp$-nilpotent complete in $\Mod_{C\tau^n}$. That is,
  \[
    \nu_n P \cong \varprojlim \nu_n (CB^\bullet({\HFp}) \otimes P).
  \]
\end{lemma}

\begin{proof}
  If $n < \infty$, this is \cite[Lemma A.12]{manifold-synthetic}, since limits in $\Mod_{C\tau^n}$ are computed in $\Syn_{\HFp}$.

  If $n = \infty$, then by \cite[Propositions 5.4, 5.6]{synthetic}, we have $\nu_\infty X = \nu X_{\HFp}$. Since $X$ is bounded below,  by \cite[Theorem 6.6]{localization-spectra}, we know that $X_{\HFp} = X_{\HFp}^{\wedge}$, the ${\HFp}$-nilpotent completion of $X$. By \cite[Proposition A.11]{manifold-synthetic}, we know that $\nu X_{\HFp}^{\wedge}$ is $\nu {\HFp} = \nu_\infty {\HFp}$-nilpotent complete.
\end{proof}

\begin{lemma}\label{lemma:y-tot}
  Let $P \in \Sp^{ft}$. Then $\H^{(n)}$ preserves the limit
  \[
    \nu_n P \overset\sim\to \varprojlim \nu_n (CB^\bullet({\HFp}) \otimes P).
  \]
\end{lemma}

\begin{proof}
  Sifted colimits in $P_\Sigma^\Sp(h_n \M_{\HFp}^\op)$ are evaluated pointwise, so we evaluate both sides on $\nu_n M \in h_n \M_{\HFp}$. The left-hand side is
  \[
    (\H^{(n)} \nu_n P)(\nu_n M) = \tau_{[0, n)} F(P, M),
  \]
  while right-hand side is given by
  \[
    \begin{aligned}
      \varinjlim F_{\Mod_{C\tau^n}}(\nu_n CB^\bullet({\HFp}) \otimes P, \nu_n M) &= \varinjlim \tau_{[0, n)} F(CB^\bullet({\HFp}) \otimes P, M) \\
                                                                                  &= \varinjlim \tau_{[0, n)} F(P, F(CB^\bullet({\HFp}), M)).
    \end{aligned}
  \]
  Since $M$ is an ${\HFp}$-module, the augmented simplicial object $F(CB^\bullet({\HFp}), M) \to M$ has extra degeneracies. So we are done.
\end{proof}

\subsection{Naturality of the comparison functor}\label{section:comp-nat}
Our ultimate goal is to use the comparison functor to compute the Adams differential, which is the long exact sequence associated to the cofiber sequence
\[
  \Sigma^{0, -k} C\tau^{m - k} \to C\tau^m \to C\tau^k.\tag{$\dagger$}
\]
More precisely, we want to look at the long exact sequence induced by applying the functor $[X, (-)\otimes_{C\tau^n} Y]_{\Mod_{C\tau^n}}^{*, *}$ to the cofiber sequence, where $X, Y \in \Mod_{C\tau^n}^{ft}$.

Since $C\tau^m \otimes_{C\tau^n} Y$ is not in $\Mod_{C\tau^n}^{ft}$ when $m < n$, we cannot simply apply \Cref{thm:y-ff} to translate this to the world of $\A^{(n)}$-modules. Nevertheless, \Cref{thm:almost-compact} tells us we can instead apply $[D(-) \otimes_{C\tau^n} X, Y]^{*, *}_{\Mod_{C\tau^n}}$ to the sequence ($\dagger$) to obtain the same result.

Thus, we are motivated to compute $\H^{(n)}(DC\tau^m \otimes_{C\tau^n} X)$ in terms of $\H^{(n)}(X)$.

\begin{thm}\label{lemma:shift-ctau}
  Let $m < n$. Then there is a natural transformation of $\A^{(n)}$-modules
  \[
    \eta\colon \A^{(m)} \otimes_{\A^{(n)}} \H^{(n)} X \to \H^{(n)} (DC\tau^m \otimes_{C\tau^n} X)
  \]
  which is an equivalence on the stable subcategory generated by $\{\nu_n Y\}_{Y \in \Sp}$. Moreover, when $X$ is of the form $\Sigma^a \nu_n Y$, the cofiber sequence induced by \textup{($\dagger$)} corresponds to the cofiber sequence induced by
  \[
    \begin{tikzcd}
      \Sigma^{k, k} \A^{(m - k)} \ar[r, "\tau^k"] & \A^{(m)} \ar[r] & \A^{(k)}.
    \end{tikzcd}
  \]
\end{thm}

\begin{remark}
  We expect the compatibility property to hold unconditionally. However, a proof eludes us.
\end{remark}

The first part naturally breaks into two lemmas.

\begin{lemma}\label{lemma:hn-nat-dual}
  Let $m \leq n$. Then there is a natural equivalence of $\A^{(n)}$-modules
  \[
    \H^{(n)} (DC\tau^m \otimes_{C\tau^n} X) \cong \H^{(m)} (C\tau^m \otimes_{C\tau^n} X).
  \]
\end{lemma}
Note that on the left-hand side, we are using the tensor product in $\Mod_{C\tau^n}$, whereas on the right, we are using the base change functor $\Mod_{C\tau^n} \to \Mod_{C\tau^m}$.

\begin{proof}
  By \Cref{lemma:psigma-nat}, we can write the right-hand side as the presheaf
  \begin{align*}
    \H^{(m)}(C\tau^m \otimes_{C\tau^n} X)(\nu_n M) &= F_{\Mod_{C\tau^m}} (C\tau^m \otimes_{C\tau^n} X, C\tau^m \otimes_{C\tau^n} \nu_n M) \\
                                                   &= F_{\Mod_{C\tau^n}} (X, C\tau^m \otimes_{C\tau^n} \nu_n M) \\
                                                   &= F_{\Mod_{C\tau^n}} (DC\tau^m \otimes_{C\tau^n} X, \nu_n M) \\
                                                   &= \H^{(n)} (DC\tau^m \otimes_{C\tau^n} X)(\nu_n M),
  \end{align*}
  where the third equality uses \Cref{thm:almost-compact}.
\end{proof}

\begin{lemma}\label{lemma:hn-nat}
  There is a natural transformation of $\A^{(m)}$-modules
   \[
     \A^{(m)} \otimes_{\A^{(n)}} \H^{(n)} X \to \H^{(m)} (C\tau^m \otimes_{C\tau^n} X)
   \]
   that is an equivalence on the stable subcategory generated by $\{\nu_n Y\}_{Y \in \Sp}$.
\end{lemma}

\begin{proof}
  Taking the dual of $C\tau^n \to C\tau^m$ gives a map $DC\tau^m \to DC\tau^n = C\tau^n$. Since $\H^{(n)}$ is contravariant, this gives a map of $\A^{(n)}$-modules
  \[
    \H^{(n)} X \to \H^{(n)}(DC\tau^m \otimes_{C\tau^n} X) \cong \H^{(m)} (C\tau^m \otimes_{C\tau^n} X).
   \]
  The desired natural transformation is then the adjoint to this map.

  One then observes that this is an equivalence when $X = \nu_n Y$, where both sides are the $m$-ary cohomology of $Y$.
\end{proof}

\begin{proof}[Proof of \Cref{lemma:shift-ctau}]
  The first part follows from \Cref{lemma:hn-nat-dual,lemma:hn-nat}. As for the second part, tracing through the proof shows that the reduction map $C\tau^m \to C\tau^k$ always corresponds to the natural projection $\A^{(m)} \to \A^{(k)}$. The map $\tau^k \colon \Sigma^{0, -k} C\tau^{m - k} \to C\tau^m$ requires more work.

  For brevity, we drop the subscripts in the tensor products. Then we have a commutative diagram
  \[
    \begin{tikzcd}
      \Sigma^{k, k} \A^{(m - k)} \otimes \H^{(n)} X \ar[r, "\tau^k"] \ar[d, dashed] & \A^{(m)} \otimes \H^{(n)} X \ar[r] \ar[d, "\eta"] & \A^{(k)} \otimes \H^{(n)} X \ar[d, "\eta"] \\
      \H^{(n)}(\Sigma^{0, k} DC\tau^{m - k} \otimes X) \ar[r, "\tau^k"] & \H^{(n)}(DC\tau^m \otimes X) \ar[r] & \H^{(n)}(DC\tau^k \otimes X)
    \end{tikzcd}
  \]
  where the dashed vertical arrow is induced by the universal property of a cofiber sequence. Our goal is to show that the dashed vertical arrow is in fact $\eta$ when $X = \nu_n Y$.

  In this case, $\eta$ is an equivalence, and the leftmost column is the $k$-connective cover of the middle column. Thus, there is a unique choice of the dashed arrow that makes the left-hand square commutes. So it suffices to show that $\eta$ also makes the left-hand square commute.\footnote{We are trying to show that selecting the dashed map to be $\eta$ gives a map of cofiber sequences, which is \emph{a priori} stronger than showing that the two squares commute. In this special case, our argument shows that the latter is in fact sufficient.}

  The trick is that we know $\H^{(n)}(D(-) \otimes X)$ sends the map
  \[
    \tau^k \colon \Sigma^{0, -k} C\tau^m \to C\tau^m
  \]
  to
  \[
    \tau^k \colon \Sigma^{k, k} \A^{(m)} \to \A^{(m)}.
  \]
  The maps labelled $\tau^k$ in the diagram above are related to these $\tau^k$ maps by the restriction maps $C\tau^m \to C\tau^{m - k}$ and $\A^{(m)} \to \A^{(m - k)}$, which $\H^{(n)}$ is also known to preserve. So in the diagram
  \[
    \begin{tikzcd}[column sep=1.4em]
      \Sigma^{k, k} \A^{(m)} \otimes \H^{(n)} X \ar[r] \ar[d, "\eta"] & \Sigma^{k, k} \A^{(m - k)} \otimes \H^{(n)} X \ar[r, "\tau^k"] \ar[d, "\eta"] & \A^{(m)} \otimes \H^{(n)} X \ar[d, "\eta"] \\
      \H^{(n)}(\Sigma^{0, k} DC\tau^m \otimes X) \ar[r] & \H^{(n)}(\Sigma^{0, k} DC\tau^{m - k} \otimes X) \ar[r, "\tau^k"] & \H^{(n)}(DC\tau^m \otimes X),
    \end{tikzcd}
  \]
  we know both the large rectangle and the left square commute. Moreover, the left-hand square exhibits the middle column as the $(m-1)$-truncation of the leftmost column, and the rightmost column is $(m - 1)$-truncated. So the right-hand square must commute as well, and we are done.
\end{proof}

\section{Locally bigraded categories}\label{section:bigraded}
Famously, the product of the Adams $E_2$ page differs from the (usual) product of the $\Ext$ groups by a sign \cite[p. 196]{structure-applications}. The goal of this \namecref{section:bigraded} is to explain where this sign comes from. Even at the prime $2$, the sign is now important, since the secondary Steenrod algebra is an algebra over $\Z/4$, not $\Z/2$.

The main issue at hand is that $\Mod_{C\tau^n}$ and $\Mod_{\A^{(n)}}$ have \emph{two} suspension functors $\Sigma^{1, 0}$ and $\Sigma^{0, 1}$. To define the composition product, we need to choose natural equivalences $\Sigma^{a, b} \Sigma^{a', b'} \cong \Sigma^{a + a', b + b'}$ in a suitably coherent fashion. While the map $\H^{(n)} \colon \Mod_{C\tau^n} \to \Mod_{\A^{(n)}}$ preserves each suspension functor individually, it does \emph{not} preserve this coherence data.

In this \namecref{section:bigraded}, our goal is to develop a framework to keep track of these coherence data. In \Cref{section:locally-graded}, we warm up on the case where there is only one suspension functor, which is relatively straightforward. In \Cref{section:locally-bigraded}, we follow the template of \Cref{section:locally-graded} to study the bigraded case. In general, it is difficult to show that a functor preserves the coherence data. However, we will show that this is automatic if one of the suspensions is the categorical suspension and the ``obvious'' coherence data is used.

In \Cref{section:sign}, we explain how these choices affect sign rules in the presence of a symmetric monoidal structure. This motivates us to impose a non-obvious choice of coherence data on $\Mod_{C\tau^n}$, which $\H^{(n)}$ then fails to preserve.

\subsection{Locally graded categories}\label{section:locally-graded}

\begin{defi}
  A locally graded category is a category $\C$ with an automorphism $[1] \colon \C \to \C$.
\end{defi}

\begin{eg}
  The category $\Sp$ of spectra is a locally graded category with automorphism given by $X[1] = \Sigma X$.
\end{eg}

\begin{eg}
  The category $\mathrm{Ab}_*$ of graded abelian groups is a locally graded category with automorphism given by $(X[1])_n = X_{n - 1}$.
\end{eg}

The structure of a locally graded category gives rise to graded mapping spaces. Let $\C$ be a locally graded category and $X, Y \in \C$. We can then define
\[
  [X, Y]^t = [X[t], Y],
\]
where $[t]$ is the $t$-fold composition of $[1]$ (using the inverse if negative).

The graded mapping spaces inherit a composition operation
\[
  [X, Y]^t \times [Y, Z]^s \to [X, Z]^{t + s}.
\]
To define this, given $f\in [X, Y]^t$ and $g \in [Y, Z]^s$, we shift $f$ to get a map
\[
  f[s] \colon X[t][s] \cong X[t + s] \to Y[s],
\]
and then compose with $g \colon Y[s] \to Z$ to get a map $X[t + s] \to Z$. Crucially, this involves identifying $X[t][s] \cong X[t + s]$. This is easy, since both are given by iterating the functor $[1]$ $(t + s)$ many times. One might have to be a bit careful when $t$ or $s$ is negative, but it turns out not to be a problem.

Nevertheless, it is worth keeping track of these identifications ``properly'', which will be crucial in the bigraded case. To do so, we define locally graded categories in an ``unbiased'' way. That is, we provide functors $[t] \colon \C \to \C$ for every $t \in \Z$, together with a coherent choice of equivalences $[t] \circ [s] \cong [t + s]$. In other words, we want an $\E_1$-map $\Z \to \Aut(\C)$.

\begin{defi}
  The category of locally graded categories is $\Cat^{B\Z}$.
\end{defi}

To reconcile the two definitions, let $S^1$ be the simplicial set given by identifying the endpoints of $\Delta^1$. Then there is an inclusion map $S^1 \hookrightarrow B\Z$ selecting $1 \in \Z$. One can check that the induced map $\Cat^{B\Z} \to \Cat^{S^1}$ is fully faithful with essential image given by those where $S^1$ selects an automorphism of the category. This then recovers our original definition of a locally graded category. Further, this lets us describe a morphism of locally graded categories as a functor $F \colon \C \to \D$ together with a natural equivalence $F[1]_\C \cong [1]_\D F$.

\subsection{Locally bigraded categories}\label{section:locally-bigraded}
A locally bigraded category is one where there are two compatible shift operators. There is now no obvious biased definition, so we head straight to the unbiased one, and then reverse-engineer the biased one afterwards.

\begin{defi}
  The category of locally bigraded categories is $\Cat^{B\Z \times B\Z}$. Given a locally bigraded category, we write the action of $(a, b) \in \Z \times \Z$ as $\Sigma^{a, b}$.
\end{defi}
We then have bigraded mapping spaces
\[
  [X, Y]^{a, b} = [\Sigma^{a, b} X, Y].
\]

As in the single-graded case, we have a fully faithful embedding $\Cat^{B\Z \times B\Z} \to \Cat^{S^1 \times S^1}$ whose essential image is given by the elements where $S^1 \times S^1$ selects automorphisms.

From this, we see that a local bigrading is given by two automorphisms $\Sigma^{1, 0}$ and $\Sigma^{0, 1}$ together with an equivalence
\[
  \Sigma^{1, 0} \Sigma^{0, 1} \simeq \Sigma^{0, 1} \Sigma^{1, 0},
\]
which we call the swapping homotopy.

Given this data, we define $\Sigma^{a, b} = (\Sigma^{1, 0})^a (\Sigma^{0, 1})^b$. Then the identifications $\Sigma^{a, b} \Sigma^{a', b'} \cong \Sigma^{a + a', b + b'}$ are given by
\begin{align*}
  \Sigma^{a, b} \Sigma^{a', b'} &\cong \Sigma^{a, 0} \Sigma^{0, b} \Sigma^{a', 0} \Sigma^{0, b'} \\
                                &\cong \Sigma^{a, 0} \Sigma^{a', 0} \Sigma^{0, b} \Sigma^{0, b'} \\
                                &\cong \Sigma^{a + a', 0} \Sigma^{0, b + b'} \\
                                &\cong \Sigma^{a + a', b + b'}
\end{align*}
In this chain, the second identification applies the swap map $a'b$ many times, and the rest are by definition.

Informally, a morphism of locally bigraded categories is a functor that commutes with the two shifts and preserves the swapping homotopy. In practice, while it is easy to check that a functor is compatible with the shifts, it is difficult to show that it preserves the swapping homotopy --- we have to write down a 3-morphism to show that a certain cube commutes.

Since we need the identification $\Sigma^{a, b} \Sigma^{a', b'} \cong \Sigma^{a + a', b + b'}$ to define composition of bigraded mapping spaces, a functor that fails to preserve this identification will fail to preserve compositions between bigraded mapping spaces. Indeed, this is the source of mismatch between the product in $\Ext$ and the product in the Adams $E_2$ page.

Fortunately for us, in all cases of interest, the bigrading is of a special form --- one of the shifts is given by categorical suspension. This can be chosen functorially, which will relieve much of our pains.

To state this formally, let $\Stbl \subseteq \Cat$ be the category of stable $\infty$-categories and exact functors.
\begin{lemma}\label{lemma:canonical-suspension}
  The projection $(\Stbl)^\Z \to \Stbl$ has a section $\Sigma \colon \Stbl \to (\Stbl)^\Z$ that selects the categorical suspension functor of each stable $\infty$-category.
\end{lemma}
This argument is due to Piotr Pstr\k{a}gowski.

\begin{proof}
  We have to produce an automorphism of $1 \in \Fun(\Stbl, \Stbl)$. Under the Grothendieck construction, the identity functor is classified by the coCartesian fibration ${}_{\Sp^\omega \backslash} \Stbl \to \Stbl$, where $\Sp^\omega$ is the category of finite spectra. The desired automorphism is then given by precomposition with $\Sigma\colon \Sp^\omega \to \Sp^\omega$.
\end{proof}

Given any stable category $\C$ and an automorphism $\Sigma^{0, 1}$, there is a local bigrading where $\Sigma^{1, 0}$ is the categorical suspension and the swapping homotopy is the natural transformation witnessing the exactness of $\Sigma^{0, 1}$. This construction can be made functorial as follows:
\begin{defi}
  Let $\Sigma \colon \Stbl \to (\Stbl)^{B\Z}$ be the suspension functor of \Cref{lemma:canonical-suspension}. Applying $(-)^{B\Z}$ to this gives a functor
  \[
    \Sigma^{B\Z}\colon  (\Stbl)^{B\Z} \to (\Stbl)^{B\Z \times B\Z}.
  \]
  If $\Sigma^{0, 1} \colon \C \to \C$ is an automorphism of a stable category, we call the image under $\Sigma^{B\Z}$ the canonical local bigrading generated by $\Sigma^{0, 1}$.
\end{defi}

The key point is that if $F \colon \C \to \D$ is a morphism between locally graded stable categories, then it is automatically a functor between the canonical locally bigraded categories. This absolves the need to consider $3$-morphisms.

\begin{eg}\label{eg:moda-bigrade}
  Let $A$ be a graded algebra over $\Z$ and $\Mod_A$ be the $\infty$-category of graded modules over $A$. This has a shift functor $[1]\colon \Mod_A \to \Mod_A$ given by shifting the internal grading.

  Recall that $\Mod_A$ is presented by the category $\Ch(A)$ of chain complexes over $A$. Then the categorical suspension functor $\Ch(A)$ is given by shifting cohomological degrees, while the internal shift $[1]$ is given by shifting internal degrees. As functors between $1$-categories, these commute on the nose, and this gives the canonical bigrading.
\end{eg}

If we give both $\Mod_{C\tau^n}$ and $\Mod_{\A^{(n)}}$ the canonical local bigrading, then $\H^{(n)}$ will be a morphism of locally bigraded categories, and everything will be nice. \Cref{eg:moda-bigrade} suggests we should indeed give $\Mod_{\A^{(n)}}$ the canonical local bigrading, since this is what we get when computing with the model structure. In the next \namecref{section:sign}, we will explain why we should \emph{not} give $\Mod_{C\tau^n}$ the canonical local bigrading.

\subsection{Sign rules}\label{section:sign}
Often, the local bigrading comes from a symmetric monoidal structure. Let $\C$ be a symmetric monoidal category, and choose $\S^{1,0}, \S^{0, 1} \in \Pic(\C)$. We can then define bigraded spheres
\[
  \S^{a, b} = (\S^{1, 0})^{\otimes a} \otimes (\S^{0, 1})^{\otimes b},
\]
and thus bigraded suspension functors
\[
  \Sigma^{a, b} = \S^{a, b} \otimes (-).
\]
To formally define a local bigrading, we choose the two shift maps to be $\Sigma^{1, 0}$ and $\Sigma^{0, 1}$, and choose the swapping homotopy to be the one induced by the symmetric monoidal structure. We call this the symmetric monoidal bigrading.

\begin{lemma}
  Suppose $\S^{1, 0} = \Sigma \mathbf{1}_\C$. Then the symmetric monoidal bigrading agrees with the canonical local bigrading generated by $\Sigma^{0, 1}$.
\end{lemma}

\begin{proof}
  We have to show that for any $X, Y \in \C$, the following diagram commutes naturally:
  \[
    \begin{tikzcd}
      \Sigma (X \otimes Y) \ar[r] \ar[d] & X \otimes \Sigma Y \ar[d] \\
      \Sigma \mathbf{1} \otimes X \otimes Y \ar[r, "\sigma \otimes Y"] & X \otimes \Sigma \mathbf{1} \otimes Y.
    \end{tikzcd}
  \]
  Then taking $X = \S^{0, 1}$, the top map is the canonical bigrading, while the bottom map is the symmetric monoidal bigrading.

  To show this, we show that the two diagonal compositions $\Sigma(X \otimes Y) \to X \otimes \Sigma \mathbf{1} \otimes Y$ are both equal to a third map
  \[
    g \colon \Sigma(X \otimes Y) \to \Sigma(X \otimes \mathbf{1} \otimes Y) \to X \otimes \Sigma \mathbf{1} \otimes Y.
  \]

  For the composite through the bottom-left, consider the diagram
  \[
    \begin{tikzcd}
      X \otimes Y \ar[r] \ar[rd] & \mathbf{1} \otimes X \otimes Y \ar[r] \ar[d, "\sigma \otimes Y"] & * \ar[r] \ar[d, equals] & \Sigma \mathbf{1} \otimes X \otimes Y \ar[d, "\sigma \otimes Y"] \\
      & X \otimes \mathbf{1} \otimes Y \ar[r] & * \ar[r] & X \otimes \Sigma \mathbf{1} \otimes Y.
    \end{tikzcd}
  \]
  Here the left triangle commutes canonically by the definition of a symmetric monoidal category, while the rest of the diagram is a map of cofiber sequences obtained by applying the natural transformation $(-) \otimes X \otimes Y \to X \otimes (-) \otimes Y$ to the cofiber sequence $\mathbf{1} \to * \to \Sigma \mathbf{1}$.

  This diagram gives two commutative diagrams of the form
  \[
    \begin{tikzcd}
      X \otimes Y \ar[d] \ar[r] & * \ar[d] \\
      * \ar[r] & X \otimes \Sigma \mathbf{1} \otimes Y,
    \end{tikzcd}
  \]
  one via the top cofiber sequence and the other via the bottom one, which correspond to two maps $\Sigma (X \otimes Y) \to X \otimes \Sigma \mathbf{1} \otimes Y$. The one via the top sequence is the bottom-left composite, while the one via the bottom sequence is the map $g$. Since the diagram of cofiber sequences commutes, it follows that these two maps agree.

  The top-right composite follows from a similar argument. Start with $\mathbf{1} \to * \to \Sigma \mathbf{1}$ and tensor with $Y$ on the right to get the commutative diagram of cofiber sequences
  \[
    \begin{tikzcd}
      Y \ar[d] \ar[r] & * \ar[d, equals] \ar[r] & \Sigma Y \ar[d] \\
      \mathbf{1} \otimes Y \ar[r] & * \ar[r] & \Sigma \mathbf{1} \otimes Y
    \end{tikzcd}
  \]
  Next we tensor this whole diagram with $X$ on the left to get
  \[
    \begin{tikzcd}
      X \otimes Y \ar[d] \ar[r] & * \ar[d, equals] \ar[r] & X \otimes \Sigma Y \ar[d] \\
      X \otimes \mathbf{1} \otimes Y \ar[r] & * \ar[r] & X \otimes \Sigma \mathbf{1} \otimes Y
    \end{tikzcd}
  \]
  Then the map through the top sequence is the top-right composite, while the one via the bottom sequence is $g$.
\end{proof}

\begin{remark}
  The proof that the diagram commutes is, of course, entirely formal. In fact, it does not use that the tensor product preserves colimits; it only involves the colimit comparison map. Once one decides to prove the result in this generality, there is only one possible proof to write down.
\end{remark}

One checks that
\begin{lemma}
  Suppose the composite
  \[
    \S^{2, 0} \cong \S^{1, 0} \otimes \S^{1, 0} \overset{\sigma}\to \S^{1, 0} \otimes \S^{1, 0} \cong \S^{2, 0}
  \]
  is multiplication by $\alpha \in \End(\S^{0, 0})$ and the corresponding one for $\S^{0, 1}$ is $\beta \in \End(\S^{0, 0})$. If we use the symmetric monoidal structure to identify $\S^{a + a', b + b'} \cong \S^{a, b} \otimes \S^{a', b'}$, then the composite
  \[
    \S^{a + a', b + b'} \cong \S^{a,b} \otimes \S^{a', b'} \overset{\sigma}\to \S^{a', b'} \otimes \S^{a, b} \cong \S^{a + a', b + b'}
  \]
  is multiplication by $\alpha^{aa'} \beta^{bb'}$.\fakeqed
\end{lemma}

Note that when we identify $\S^{a + a', b + b'} \cong \S^{a, b} \otimes \S^{a', b'}$, we have to move $\S^{a', 0}$ over $\S^{0, b}$, but the swap map $\sigma$ immediately moves it back. When determining sign rules of bigraded homotopy groups, the first move uses the homotopy from the definition of the bigrading, and the second uses the symmetric monoidal structure. For the symmetric monoidal bigrading, these agree, so they cancel out. If the two homotopies differ by $(-1)$, then we pick up an extra sign of $(-1)^{ab' + a'b}$.

\begin{lemma}
  For $\Syn_E$, the multiples for $\S^{1, -1}$ and $\S^{1, 0}$ are both $-1$.
\end{lemma}

\begin{proof}
  The former is a general property of categorical suspension. The latter follows from the fact that $\nu \colon \Sp_E^{fp} \to \Syn_E$ is symmetric monoidal and $\nu \S^1 = \S^{1, 0}$.
\end{proof}

Under the canonical bigrading, we get a sign rule of $(-1)^{aa' + a'b + ab'}$, which is bizarre; a more natural sign rule is $(-1)^{aa'}$, which depends only on the stem and not the filtration. For example, under the sign rule of $(-1)^{aa' + a'b + ab'}$, both $h_0$ and $\tau$ multiplications anti-commute with elements in odd stems. To fix the sign rule, we insert a sign:

\begin{defi}[{\cite[Remark 4.10]{synthetic}}]
  Viewing $\Syn_E$ as a category of sheaves over $\Sp_E^{fp}$, the Adams bigrading on $\Syn_E$ is generated by
  \[
    (\Sigma^{1, 0} X)(P) = X(\Sigma^{-1} P),\quad (\Sigma^{1, -1}X)(P) = \Sigma X(P),
  \]
  where the swap map is given by $(-1)$ times the canonical bigrading.
\end{defi}

This then results in a sign rule of $(-1)^{aa'}$.

If we used the canonical bigrading on $\Syn_E$ (hence $\Mod_{C\tau^n}$), then since the functor $\H^{(n)} \colon \Mod_{C\tau^n} \to \Mod_{\A^{(n)}}^\op$ is exact, it is a map of locally bigraded categories, hence preserves the bigraded composition product. Since we decide to use the Adams bigrading on $\Syn_E$ instead, under the functor $\H^{(n)}$, composition products differ by a sign of $(-1)^{(-b')(a + b)} = (-1)^{s' t}$.

\part{Computing \texorpdfstring{$E_3$}{E3} page data}\label{part:secondary}
\section{Overview}
Following Baues, we write $\A = \A^{(1)}$ and $\B = \A^{(2)}$. The objective of this \namecref{part:secondary} is to understand how to do computations in $\Mod_{C\tau^2}$ via the fully faithful embedding
\[
  \Mod_{C\tau^2}^{ft} \to \Mod_{\B}^\op.
\]
These computations will then be used to compute the Adams spectral sequence in \Cref{part:computation}.

We start with \Cref{section:modules}, where we seek to understand $\Mod_{\B}$ via a model category presentation. After describing the differential graded algebra $\B$ and studying some $\B$-modules in depth, we learn how to construct a cofibrant replacement of $\H^{(2)}X$ by lifting a free $\A$-resolution of $H^*X$.

In \Cref{section:e3}, we take this cofibrant replacement and use it to compute the data we sought, namely $d_2$ differentials, $\Mod_{C\tau^2}$ products and $\Mod_{C\tau^2}$ Massey products. The $d_2$ differentials are computed as the obstruction to lifting an $\Ext$ class over $\A$ to an $\Ext$ class over $\B$, while the $\Mod_{C\tau^2}$ products and Massey products are obtained by lifting chain maps and chain homotopies over $\A$ to ones over $\B$. Note that our algorithm to compute $d_2$ differentials ends up being identical to that of \cite{baues-e3}, but our proofs are independent (apart from the computation of $\B$ itself).

Finally, in \Cref{section:data}, we discuss our implementation of the algorithm. We will give an overview of the resulting dataset and provide instructions to reproducing the data. We then discuss the performance characteristics of our implementation to understand how the run time grows with the stem.

\section{The category of secondary Steenrod modules}\label{section:modules}
\subsection{The secondary Steenrod algebra}\label{section:nassau}
In this \namecref{section:nassau}, we explicitly describe the secondary Steenrod algebra as a differential graded algebra at the prime $2$. This was originally computed by Baues in \cite{baues-book}. To make use of his computations, we need to reconcile our definition with his.

\begin{thm}
  $\B$ is equivalent to the secondary Steenrod algebra of \cite{baues-book}.
\end{thm}

\begin{proof}
  By Morita theory, $\B$ is uniquely characterized by the fact that
  \[
    \Free_{\B} \cong h_2 \M_{\HFp}^\op.
  \]
  This was shown for Baues' secondary Steenrod algebra in \cite[Theorem 5.5.6]{baues-book} and ours in \Cref{cor:free-an}.
\end{proof}

Since $\B$ is a differential graded algebra, it admits multiple presentations as chain complexes. Baues' original presentation was large and unwieldy. In \cite{nassau-secondary}, Nassau discovered a smaller and simpler presentation of the secondary Steenrod algebra, which is what we shall regurgitate here.

Recall that the homotopy groups of $\B$ are given by
\[
  \pi_* \B = \A[\tau] / \tau^2.
\]
In particular, they are concentrated in cohomological degrees $0$ and $1$. Thus, we can represent $\B$ as a $2$-term chain complex
\[
  \B = \left(\begin{tikzcd}
    \B_1 \ar[d, "d^\B"] \\
    \B_0
  \end{tikzcd}\right).
\]
The structure of being a differential graded algebra means $\B_0$ is a ring, $\B_1$ is a $\B_0$-$\B_0$-bimodule, and $d^\B$ is a bimodule homomorphism such that
\[
  (d^\B a)b = a(d^\B b) \text{ for all }a, b \in \B_1.
\]

One should think of $\B_0$ as an enlargement of $\A$ so that certain products are not literally zero. For example, if we want the secondary cohomology operation $\langle \beta, \beta, -\rangle$ to ever be non-zero, the product $\beta \beta$ cannot vanish in $\B_0$; it must be killed by a non-trivial homotopy in $\B_1$. The $\B_0$-$\B_0$-bimodule structure on $\B_1$ then encodes various Massey product information.

By definition, $\B_1$ and $\B_0$ fit in a long exact sequence
\[
  \begin{tikzcd}
    0 \ar[r] & \pi_1 \B = \A\{\tau\} \ar[r] & \B_1 \ar[r, "d^\B"] & \B_0 \ar[r, "\pi^\B"] \ar[r] & \pi_0 \B = \A \ar[r] & 0.
  \end{tikzcd}
\]
Our model of $\B$ admits the following crucial simplifying property:
\begin{lemma}\pushQED{\qed}
  The long exact sequence splits as
  \[
    \B_1 \cong \ker \pi^B \oplus \A\{\tau\},\quad |\tau| = 1.
  \]
  Further, this splitting is compatible with the right $\B_0$-action and the left $\ker \pi^\B$-action.\qedhere
\end{lemma}

Under this splitting, the left action of $\B_0$ on $\B_1$ necessarily takes the form
\[
  a \cdot (r, p) = (ar, A(\pi^\B(a), r) + \pi^\B(a) p)
\]
for some function
\[
  A \colon \A \otimes \ker \pi^\B \to \A\{\tau\}.
\]
One should think of this function $A$ as carrying the ``Massey product information'' in $\B$. For example, if $a, b, c \in \A$ are such that $ab = bc = 0$ and $\tilde{b}, \tilde{c} \in \B_0$ are lifts of $b, c$ respectively, then $A(a, \tilde{b}\tilde{c}) \in \langle a, b, c\rangle$.

\begin{eg}
  In the secondary $\A(0)$ that we computed in \Cref{section:a0}, we had $\B_0 = \Z/4\{1, \beta\}$ and
  \[
    A(\beta, p) = \tau.
  \]
\end{eg}

\begin{remark}
  All the non-trivial information in $\B$ is contained in the function $A$. When choosing $\B_0$, we are mostly just trying to fatten $\A$ enough to make room for the non-trivial $A$. In Baues' original model, it was simply taken to be the free $\Z/4$-algebra on $\{\Sq^n\}_{n > 0}$.
\end{remark}

\begin{remark}
  One can show that if $r \in \ker \pi^\B$ is in the center of $\B_0$, then $A(-, r)$ is a derivation on $\A$. When $r = 2$, the derivation $A(-, 2)$ is the Kirstensen derivation that sends $\Sq^n \to \Sq^{n - 1}$.
\end{remark}

In the rest of the \namecref{section:nassau}, we will describe the ring $\B_0$ and the function $A$. We start with $\B_0$, which is in fact a Hopf algebra. As in the ordinary case, its dual admits a nice ``geometric'' description --- it is the Hopf algebra representing power series of the form
\[
  f(x) = \sum_{k \geq 0} \xi_k x^{2^k} + \sum_{0 \leq k < l} 2 \xi_{k, l} x^{2^k + 2^l}
\]
under composition mod $4$. This gives a natural inclusion $\A_* \hookrightarrow (\B_0)_*$, whose dual defines our projection $\pi^\B \colon \B_0 \to \pi_0 \B$.

Explicitly, the Hopf algebra $(\B_0)_*$ is given by
\[
  (\B_0)_* = \Z/4[\xi_k, 2 \xi_{k, l} \mid 0 \leq k < l, \xi_0 = 1]
\]
with the coproduct
\[
  \begin{aligned}
    \Delta \xi_n &= \sum_{i + j = n} \xi_i^{2^j} \otimes \xi_j + 2 \sum_{0 \leq k < l} \xi_{n - 1 - k}^{2^k} \xi_{n - 1 - l}^{2^l} \otimes \xi_{k, l}\\
    \Delta \xi_{n, m} &= \xi_{n, m} \otimes 1 + \sum_{k \geq 0} \xi_{n - k}^{2^k} \xi_{m - k}^{2^k} \otimes \xi_{k + 1} \\
                      &\hphantom{= \xi_{n, m} \otimes 1}+ \sum_{0 \leq k < l} (\xi_{n - k}^{2^k} \xi_{m - l}^{2^l} + \xi_{m - k}^{2^k} \xi_{n - l}^{2^l}) \otimes \xi_{k, l}.
  \end{aligned}
\]
That is, $(\B_0)_*$ is the sub-Hopf algebra of $\Z/4[\xi_k, \xi_{k, l}]$ generated by the $\xi_k$ and $2\xi_{k, l}$. The ring $\B_0$ is then given by $\Hom((\B_0)_*, \Z/4)$. It is generated by the following elements:
\begin{defi}
  Define $\Sq(R)$ and $Y_{k, \ell}$ to be dual to $\xi^R$ and $2\xi_{k, l}$ under the monomial basis.\footnote{Our indexing differs slightly from Nassau's.}
\end{defi}

It is easy to check that
\begin{lemma}\pushQED{\qed}
  $Y_{k, \ell} \Sq(R)$ is dual to $\xi^R \xi_{k, \ell}$ under the monomial basis. Further,
  \[
    \pi^\B(Y_{k, \ell}) = 0\text{ and } \pi^\B(\Sq(R)) = \Sq(R).\qedhere
  \]
\end{lemma}

\begin{lemma}\pushQED{\qed}\label{lemma:y-prod-zero}
  \[
    Y_{*, *} Y_{*, *} = 2 Y_{*, *} = 0.\qedhere
  \]
\end{lemma}
\begin{remark}
  Since $2Y_{*, *} = 0$, we prefer to think of the $\Sq(R)$ in $Y_{k, \ell} \Sq(R)$ as an element of $\A$ instead of $\B_0$. Similarly, there is a left action of $\A$ on the $Y_{*, *} \Sq(*)$ terms.
\end{remark}

To describe the rest of the multiplication, we let $\daleth: \A_* \otimes \A \to \A$ be the contraction operator. In the Milnor basis, we have
\[
  \daleth(\xi^R, \Sq(S)) = \Sq(S - R),
\]
where $\Sq(S - R)$ is zero if any entry is negative.

\begin{lemma}\pushQED{\qed}
  If $a \in \A$, then
  \[
    a Y_{k, l} = \sum_{i, j \geq 0} Y_{k + i, l + j} \daleth(\xi_i^{2^k} \xi_j^{2^l}, a),
  \]
  where we set
  \[
    Y_{k, l} =
    \begin{cases}
      Y_{l, k} & k > l,\\
      2 \Sq(\Delta_{k + 1}) & l = k.
    \end{cases}\qedhere
  \]
\end{lemma}
Here $\Delta_k$ is the sequence that is $1$ in the $\xi_k$ position and $0$ elsewhere.

It remains to determine the multiplication between the $\Sq(R)$. Recall the following definition in the multiplication of $\A$ under the Milnor basis:

\begin{defi}
  Let $X = (x_{ij})$ be a matrix indexed on the non-negative integers. Define
  \begin{gather*}
    r_i(X) = \sum_j 2^j x_{ij},\quad s_j(X) = \sum_i x_{ij},\quad t_n(X) = \sum_{i + j = n} x_{ij}, \\
    R(X) = (r_1(X), r_2(X), \ldots ),\quad S(X) = (s_1(X), \ldots ),\quad T(X) = (t_1(X), \ldots), \\
    b(X) = \frac{\prod t_n!}{\prod x_{ij}!} = \prod_n \binom{t_n}{x_{n0}\; \cdots \; x_{0n}} \in \Z.
  \end{gather*}
\end{defi}

\begin{thm}[{\cite[Theorem 4b]{milnor-steenrod}}]\pushQED{\qed}
  We have
  \[
    \Sq(R) \Sq(S) = \sum_{\substack{R(X) = R\\ S(X) = S}} b(X) \Sq(T(X)).\qedhere
  \]
\end{thm}

Dualizing the secondary coproduct formula gives
\begin{thm}\pushQED{\qed}
  \begin{multline*}
    \Sq(R) \Sq(S) = \sum_{k \geq 0} \sum_{0 \leq m < n} Y_{m + k, n + k} \daleth(\xi_m^{2^k} \xi_n^{2^k}, \Sq(R)) \daleth(\xi_{k + 1}, \Sq(S))\\
    + \sum_{\substack{R(X) = R\\ S(X) = S}} b(X) \Sq(T(X)).
  \end{multline*}\popQED
\end{thm}

This completes the description of $\B_0$. It remains to describe the function $A$.
\begin{lemma}\pushQED\qed
  We have
  \[
    \begin{aligned}
      A(a, 2) &= \daleth(\xi_1, a),\\
      A(a, Y_{k, \ell}) &= \sum_{i, j \geq 0} Z_{k + i, l + j} \daleth(\xi_i^{2^k} \xi_j^{2^{\ell}}, a),\\
      A(a, r \Sq(R)) &= A(a, r) \Sq(R),
    \end{aligned}
  \]
  where
  \[
    Z_{k, \ell} =
    \begin{cases}
      0 & k < \ell,\\
      \Sq(\Delta_k + \Delta_\ell) & k \geq \ell.
    \end{cases}\qedhere
  \]
\end{lemma}

\subsection{Periodic \texorpdfstring{$\B$}{B}-modules}\label{section:periodic}
Equipped with a description of $\B$, we can now describe the category $\Mod_\B$. This admits the expected model category presentation.
\begin{thm}
  Let $A$ be a $\Z$-graded differential graded algebra and $\dgMod_A$ the $1$-category of differential graded modules over $A$. Then there is a model structure on $\dgMod_A$ where
  \begin{enumerate}
    \item fibrations are epimorphisms;
    \item weak equivalences are homology isomorphisms; and
    \item if $j \colon M \to N$ is a cofibration of graded $\Z$-chain complexes, then $A \otimes j\colon A \otimes M \to A \otimes N$ is a cofibration.
  \end{enumerate}
  Further, this model category presents $\Mod_A$, where we view $A$ as a graded $\E_1$-ring in $\Mod_\Z$.
\end{thm}

\begin{proof}
  If $A$ is cofibrant as a chain complex over $\Z$, then this is \cite[Theorem 4.3.3.17]{ha}. Otherwise, let $A'$ be a cofibrant replacement of $A$ in $\Alg(\operatorname{Ch}(\Z))$. By \cite[Theorem 3.1]{schwede-shipley}, we know $A'$ is cofibrant as an $\Z$-chain complex. By \cite[Corollary 3.4]{dgm}, the base change adjunction $\dgMod_A \rightleftharpoons \dgMod_{A'}$ is a Quillen equivalence. So we are done.
\end{proof}

The goal of this \namecref{section:periodic} is to understand $\B$-modules of the form $\H^{(2)}X$. Instead of computing $\H^{(2)} X$ directly, our strategy is to start with $H^*X$ and use obstruction theory to classify $\B$-modules that look like potential candidates for $\H^{(2)} X$.

Recall from \Cref{lemma:hn-homotopy} that $\H^{(2)}X$ and $H^*X$ are related by the equations
\[
  \pi_0 \H^{(2)} X = H^*X,\quad \pi_* \H^{(2)}X = (\pi_0 \H^{(2)}X) [\tau] / \tau^2.
\]

\begin{defi}
  A $\B$-module $M$ is periodic if
  \[
    \pi_* M = \pi_0 M[\tau] /\tau^2
  \]
  as a $\pi_* \B$-module. We say $M$ is a lift of the $\A$-module $\pi_0 M$.
\end{defi}
We should think of the category of periodic $\B$-modules as the secondary version of the heart of $\Mod_{\A}$. In particular, it is a $2$-category.

\begin{thm}
  Let $\bar{M}$ be an $\A$-module. Then the obstruction to lifting $\bar{M}$ to a periodic $\B$-module lies in $\Ext^{3, 1}_\A(\bar{M}, \bar{M})$.

  If the obstruction vanishes, then the set of lifts is a torsor over $\Ext^{2, 1}_\A(\bar{M}, \bar{M})$.
\end{thm}

\begin{proof}
  This follows from (the proof of) \cite[Theorem 4.9]{abstract-goerss-hopkins}.
\end{proof}

This obstruction theory has a natural interpretation. One can think of the ordinary cohomology groups $H^*X$ as encoding how one builds $X$ out of spheres, except we only remember maps up to filtration $1$. The secondary cohomology group then seeks to remember these attaching maps up to filtration $2$. If the attaching maps of $\bar{M}$ supported $d_2$ differentials, then it would be impossible to lift to a periodic $\B$-module, and this obstruction is captured in $\Ext^{3, 1}_\A(\bar{M}, \bar{M})$. If these obstructions vanish, then the set of ways to lift the filtration 1 maps to include filtration 2 information is then a torsor over $\Ext^{2, 1}_\A(\bar{M}, \bar{M})$.

In our case, if $\bar{M}$ is the cohomology of a spectrum, then we know a lift exists, namely the secondary cohomology of said spectrum. For many spectra of interest, the group $\Ext^{2, 1}_\A(\bar{M}, \bar{M})$ is trivial, so there is exactly one lift. Thus, any lift one can write down will work. When the group is non-trivial, one can show that if two lifts differ by $\chi \in \Ext^{2, 1}_\A(\bar{M}, \bar{M})$, then the $d_2$ differentials of the lifts differ by multiplication-by-$\chi$ (e.g.\ this follows from inspecting the $d_2$ algorithm we present later). To find the right lift, one will have to manually compute a small number of $d_2$'s.

In \Cref{section:resolution}, we will describe an explicit procedure to lift an $\A$-module by lifting its minimal free resolution. In this \namecref{section:periodic}, we will instead focus on understanding a few key examples of periodic $\B$-modules.

\begin{notation}
  Let $M$ be a periodic $\B$-module. Then $M$ can be represented by a $2$-term chain complex. That is, it is zero outside of cohomological degrees $0$ and $1$. We will write this chain complex as
  \[
    M = \left(\begin{tikzcd}
      M_1 \ar[d, "d^M"] \\
      M_0
    \end{tikzcd}\right).
  \]
  We write $\pi^M\colon M_0 \to \pi_0 M$ for the natural projection map.
\end{notation}

When manipulating periodic $\B$-modules, we often make use of the following property, which we have already seen for $\B$ itself:
\begin{lemma}
  Let $M$ be a periodic $\B$-module. If $a \in \B_1$ and $m \in M_1$, then
  \[
    (d^\B a) m = a (d^M m).
  \]
\end{lemma}

\begin{proof}
  $am = 0$ for degree reasons, and apply the Leibniz rule.
\end{proof}

We start with the simplest secondary Steenrod module, namely the secondary cohomology of the sphere. This is the natural analogue of $k = \F_p$.
\begin{defi}
  We define $\kk = \H^{(2)}(\S)$.
\end{defi}

\begin{lemma}
  We have
  \[
    \kk = \left(\begin{tikzcd}
        \Z/p\{\mu_0\} \oplus \Z/p\{\tau\} \ar[d, "d"] \\
        \Z/p^2
    \end{tikzcd}\right),\quad\sidedeg
  \]
  The $\B$ action is given by
  \[
    \tilde{\beta} \mu_0 = \tau,
  \]
  where $\tilde\beta \in \B_0$ is any representative of $\beta \in \pi_0 \B = \A$ (as usual, $\beta = \Sq^1$ when $p = 2$). For degree reasons, there are no other possible non-trivial actions.
\end{lemma}
Note that $\mu_0$ is the null-homotopy of $p$, and the action encodes the fact that $\beta$ detects $p$.

\begin{proof}
  Since $\Ext^{2, 1}_\A(k, k) = 0$, \emph{any} periodic $\B$-module lifting $k$ must be a model of $\kk$. Thus, the only work is to check that we described a valid $\B$-module structure, which is straightforward.
\end{proof}

The other family of periodic $\B$-modules we are interested in is the cohomology of Eilenberg--Maclane spectra, i.e.\ free $\B$-modules.

\begin{defi}
  A free $\B$-module is a direct sum of internal degree shifts of $\B$.
\end{defi}
Note that we consider the choice of generators part of the structure of a free $\B$-module.

Recall that $\B_1$ admits a splitting
\[
  \B_1 = \ker \pi^\B \oplus \pi_0 \B\{\tau\}.
\]
Thus, every free module $M$ also comes with a standard\footnote{We shall use ``standard'' to refer something that results from a concrete but non-canonical choice we have made.} splitting
\[
  M_1 \cong \ker \pi^M \oplus \pi_0 M\{\tau\}.
\]
We refer to these two components as the $\ker \pi$ component and the $\tau$ component respectively. Of course, we also have such a splitting in the case of $\kk$ by our explicit description.

We now turn to homomorphisms between periodic $\B$-modules. Since we are going to take free resolutions of periodic $\B$-modules, we are only interested in $\B$-module homomorphisms out of free modules. One immediately sees that

\begin{lemma}
  Let $M$ be a free $\B$-module and $N$ any $\B$-module. Then the natural map
  \[
    [M, N]_{\B} \to [\pi_0 M, \pi_0 N]_{\A}
  \]
  is a bijection.\fakeqed
\end{lemma}

In other words, given a map $\pi_0 M \to \pi_0 N$, there is a unique lift to a map $M \to N$ up to homotopy. While we do not get a well-defined chain map $M \to N$, we can fix some choices once and for all. We fix sections
\[
  \begin{aligned}
    \sigma^\B \colon \A = \pi_0 \B &\to \B_0\\
    \sigma^\kk \colon k = \pi_0 \kk &\to \kk_0.
  \end{aligned}
\]
of $\pi^\B$ and $\pi^\kk$ as functions between sets. These naturally extend to sections for free modules as well.

With these choices, if we have a map $f\colon \pi_0 M \to \pi_0 N$ where $N$ is either free or $\kk$, then we get a standard lift to a chain map $\tilde{f}\colon M \to N$, which we can depict as
\[
  \begin{tikzcd}
    M_1 \ar[d, "d^M"] \ar[r, "f_1"] & N_1 \ar[d, "d^N"] \\
    M_0 \ar[r, "f_0"] & N_0.
  \end{tikzcd}
\]

Given chain maps $f, g$, a homotopy between them is a $\B_0$-module map $h\colon M_0 \to N_1$ such that
\[
  f_0 - g_0 = d^N h,
\]
i.e.\ a lift of the difference along $d^N$.\footnote{The definition of a chain map requires $f_1 - g_1 = h d^M$ as well. However, if $M$ is free, then this is automatic. Indeed, any element in $M_1$ is of the form $am$, where $a \in \B_1$ and $m \in M_0$. Then
\[
  hd^M (am) = h(d^\B(a) m) = d^\B(a) h(m) = a d^N(h(m)) = a(f_0 - g_0)(m) = (f_1 - g_1)(am).
\]
Here the last equality uses the fact that $f$ and $g$ are maps of $\B$-modules, while the third equality uses the crucial relation $a (dm) = (da) m$. This will be a common theme in our future manipulations, where the $(-)_1$ version of the equation we have to satisfy follows formally from the $(-)_0$ version when the source is free. The proof is largely similar and we will not comment further.} In practice, to specify a homotopy, we don't have to specify all of $h$. Since $d^N$ is injective when restricted to the $\ker \pi$ component of $N_1$, this component of $h$ must be equal to $f_0 - g_0$ itself. The freedom in choosing the homotopy lies in the $\tau$ component. Thus, given our choices, we can thus identify a homotopy with a map $M_0 \to \pi_0 N\{\tau\}$, which necessarily factors through $\pi^M$ to give a map $\pi_0 M \to \pi_0 N\{\tau\}$.

\begin{remark}
  In general, the space of homotopies is canonically a torsor over $\Hom_\A(\pi_0 M, \pi_0 N\{\tau\})$. After making all our choices, we have found a basepoint for this space, namely the homotopy with trivial $\tau$ component. We can then identify the space of homotopies with $\Hom_\A(\pi_0 M, \pi_0 N\{\tau\})$ itself.
\end{remark}

\subsection{Free resolutions}\label{section:resolution}
We are now ready to construct a free resolution of a periodic $\B$-module $M$, which will give us a cofibrant replacement of $M$ in $\dgMod_\B$.

We start by taking a free resolution $\bar{P}^\bullet \to \pi_0 M$ of $\pi_0 M$. As usual, each $\bar{P}^{(s)}$ is a free $\A$-module with a fixed choice of generators.

The previous \namecref{section:periodic} gives us a lift of this to a sequence of free $\B$-modules
\[
  \begin{tikzcd}
    \cdots \ar[r] & P^{(3)} \ar[r, "\partial^{(3)}"] & P^{(2)} \ar[r, "\partial^{(2)}"] & P^{(1)} \ar[r, "\partial^{(1)}"] & P^{(0)} \ar[r, "\epsilon"] & M
  \end{tikzcd}
\]
such that the composites of successive maps are homotopic to zero. This alone does not let us assemble this into a cofibrant replacement in $\dgMod_\B$. What we need is a suitable choice of null-homotopy of each composite $\partial^{(k - 1)} \partial^{(k)}$.

\begin{defi}[{\cite[Definition 3.1]{baues-e3}}]
  A secondary chain complex is a sequence
  \[
    \begin{tikzcd}
      \cdots \ar[r] & P^{(3)} \ar[r, "\partial^{(3)}"] & P^{(2)} \ar[r, "\partial^{(2)}"] & P^{(1)} \ar[r, "\partial^{(1)}"] & P^{(0)}
    \end{tikzcd}
  \]
  of periodic $\B$-modules together with specified null-homotopies of $\partial^{(k - 1)} \partial^{(k)}$, such that all three-fold Massey products $\langle \partial^{(k - 2)}, \partial^{(k - 1)}, \partial^{(k)} \rangle$ vanish.
\end{defi}

Writing each module out as a 2-term chain complex itself, we can expand this to a diagram
\[
  \begin{tikzcd}[row sep = large, column sep = large]
    \cdots \ar[r] \ar[d] & P^{(3)}_1 \ar[r, "\partial^{(3)}_1"] \ar[d, "d^{(3)}"] & P^{(2)}_1 \ar[r, "\partial^{(2)}_1"] \ar[d, "d^{(2)}"] & P^{(1)}_1 \ar[r, "\partial^{(1)}_1"] \ar[d, "d^{(1)}"] & P^{(0)}_1 \ar[d, "d^{(0)}"]\\
    \cdots \ar[r] \ar[urr, gray!50!black, pos=0.4] & P^{(3)}_0 \ar[r, "\partial^{(3)}_0"] \ar[urr, "h^{(3)}", gray!50!black, pos=0.4] & P^{(2)}_0 \ar[r, "\partial^{(2)}_0"] \ar[urr, "h^{(2)}", gray!50!black, pos=0.4] & P^{(1)}_0 \ar[r, "\partial^{(1)}_0"] & P^{(0)}_0.
  \end{tikzcd}
\]
The condition of being a secondary chain complex is then
\[
  \begin{aligned}
    d^{(s - 1)} \partial_1^{(s)} &= \partial_0^{(s)} d^{(s)},\\
    \partial_0^{(s - 1)} \partial_0^{(s)} = d^{(s - 2)} h^{(s)}, &\quad\;\, \partial_1^{(s - 1)} \partial_1^{(s)} = h^{(s)} d^{(s)},\\
    h^{(s - 1)} \partial^{(s)}_0 &= \partial^{(s - 2)}_1 h^{(s)}.
  \end{aligned}
\]
These say, respectively, that
\begin{itemize}
  \item Each $\partial^{(s)}$ is a chain map;
  \item $h^{(s)}$ is a null-homotopy of $\partial^{(s - 1)} \partial^{(s)}$; and
  \item The bracket $\langle \partial^{(s - 2)}, \partial^{(s - 1)}, \partial^{(s)}\rangle$ vanishes.
\end{itemize}
Note that when the source is free, the equations $d\partial_1 = \partial_0 d$ and $\partial_1 \partial_1 = hd$ are implied by the others by $\B$-linearity.

Having chosen such homotopies, we can now define

\begin{defi}
  Let $P^{\bullet}$ be a secondary chain complex. Then the total chain complex $\Tot(P^{\bullet})$ is the chain complex\footnote{
    There are many possible choices of sign when forming the total chain complex. The choice of sign here is motivated by two concerns:
    \begin{itemize}
      \item The inclusion map $P^{(0)} \hookrightarrow \Tot(P^\bullet)$ should be given by the obvious inclusion. This precludes the last differential from being $\begin{pmatrix} \partial_0^{(1)} & -d^{(0)}\end{pmatrix}$, which is a more common sign convention
      \item Applying $\A \otimes_{\B}(-)$ to the total chain complex should yield the same complex as $\A \otimes_{\B}(-)$ applied to the secondary chain complex. This requires the top-left entry to be $\partial_0^{(s)}$ instead of $-\partial_0^{(s)}$.
    \end{itemize}
    Of course, these choices are immaterial, but we believe our choice of signs makes it slightly easier to reason about various factors.
  }
  \[
    \setlength\arraycolsep{1pt}
    \begin{tikzcd}[column sep = 4.7em, ampersand replacement=\&]
      P_0^{(3)} \oplus P_1^{(2)} \ar[r, "{\begin{pmatrix}\partial^{(3)}_0 & d^{(2)} \\ -h^{(3)} & -\partial^{(2)}_1 \end{pmatrix}}"] \& P_0^{(2)} \oplus P_1^{(1)} \ar[r, "{\begin{pmatrix}\partial^{(2)}_0 & d^{(1)} \\ -h^{(2)} & -\partial^{(1)}_1 \end{pmatrix}}"] \& P_0^{(1)} \oplus P^{(0)}_1 \ar[r, "{\begin{pmatrix} \partial_0^{(1)} & d^{(0)} \end{pmatrix}}"] \& P_0^{(0)}.
    \end{tikzcd}
  \]
\end{defi}
One readily checks that the conditions of being a secondary chain complex translate to the totalization being a chain complex, and that the natural $\B$-module structure gives it a structure of a differential graded module over $\B$.

\begin{thm}
  If each $P^{(k)}$ is cofibrant, then so is $\Tot(P^{\bullet})$.
\end{thm}

\begin{proof}
  Let $F_s \Tot(P^\bullet)$ be the subcomplex consisting of the $P^{(k)}$ terms with $k \leq s$. Then $\Tot(P^{\bullet}) = \colim F_s \Tot(P^{\bullet})$. So it suffices to show that $F_{s - 1} \Tot(P^{\bullet}) \to F_s \Tot(P^{\bullet})$ is a cofibration. But this fits in the pushout diagram
  \[
    \begin{tikzcd}
      S^{s - 1} \otimes P^{(s)} \ar[r] \ar[d] & F_{s - 1} \Tot(P^{\bullet}) \ar[d] \\
      D^s \otimes P^{(s)} \ar[r] & F_s \Tot(P^{\bullet})
    \end{tikzcd}
  \]
  and the left-hand map is a cofibration. So we are done.
\end{proof}

Filtering by cohomological degree gives a spectral sequence

\begin{lemma}\pushQED{\qed}
  If $P^{\bullet}$ is a secondary chain complex, then there is a spectral sequence
  \[
    E^{p, q}_2 = H^p(\pi_q P^{\bullet}) \implies \pi_{p + q} \Tot(P^\bullet).\qedhere
  \]
\end{lemma}
In particular,
\begin{cor}
  If $\bar{P}^\bullet \to \pi_0 M$ is a free resolution and $P^\bullet \to M$ a lift to a secondary chain complex, then the induced map $\Tot(P^\bullet) \to M$ is a weak equivalence.\fakeqed
\end{cor}
This leaves the question of how one can construct the secondary chain complex in the first place. We start with the sequence
\[
  \begin{tikzcd}
    \cdots \ar[r] & P^{(3)} \ar[r, "\partial^{(3)}"] & P^{(2)} \ar[r, "\partial^{(2)}"] & P^{(1)} \ar[r, "\partial^{(1)}"] & P^{(0)} \ar[r, "\epsilon"] & M
  \end{tikzcd},
\]
and our goal is to choose homotopies inductively to satisfy
\[
  \partial^{(s - 2)}_1 h^{(s)} = h^{(s - 1)} \partial^{(s)}_0 .
\]
To simplify notation, we assume $M_1$ admits a splitting into a $\ker \pi$ component and a $\tau$ component as well; the argument goes through with slight modifications in the general case.

The first potential non-zero homotopy is $h^{(1)}$, which we can choose arbitrarily, since the equation it has to satisfy takes values in the zero group. Inductively, assume we have made valid choices of $h^{(k)}$ for $k < s$.

As previously discussed, the $\ker \pi$ component of $h^{(s)}$ is forced to be $\partial_0^{(s - 1)} \partial_0^{(s)}$, and we have the freedom to choose the $\tau$ component, which we call $h^{(s)}_\tau$.

Let $g$ be a generator of $P^{(s)}$, and write
\[
  \bar{\partial}^{(s)} g = \sum \alpha^i g_i,
\]
where $\bar{\partial}$ is the differential of $\bar{P}^\bullet$ and $\{g_i\}$ is a set of generators of $P^{(s - 1)}$. We can then write the $\tau$ component of the equation as
\[
  \bar{\partial}^{(s - 2)} h_\tau^{(s)} g = \sum \left(\alpha^i h^{(s - 1)}_\tau (g_i) + A\left(\alpha^i, \partial_0^{(s - 2)} \partial_0^{(s - 1)}g_i\right) \right) \equiv t_g.
\]
Thus, we want to choose $h^{(s)}_\tau g$ to be a lift of $t_g$ along $\bar{\partial}^{(s - 2)}$.

\begin{lemma}
  For any valid choice of $h^{(k)}$ for $k < s$, the equation can be solved.
\end{lemma}
Thus, we can choose the homotopies iteratively.
\begin{proof}
  By exactness, we need to check that $\bar{\partial}^{(s - 3)} t_g = 0$. Informally, this follows from the Toda bracket manipulation
  \[
    \partial^{(s - 3)} \langle \partial^{(s - 2)}, \partial^{(s - 1)}, \partial^{(s)}\rangle = \langle \partial^{(s - 3)}, \partial^{(s - 3)}, \partial^{(s - 1)} \rangle \partial^{(s)} = 0.
  \]

  In more detail, we can identify $\bar{\partial}^{(s - 3)} t_g$ as the $\tau$ component of $\partial_1^{(s - 3)} h^{(s - 1)} \partial_0^{(s)}$. It is convenient to temporarily pick an arbitrary (and immaterial) value for $h^{(s)}_\tau$, so that $h^{(s)}$ is a null-homotopy of $\partial^{(s)} \partial^{(s - 1)}$. Then we have
  \[
    \partial_1^{(s - 3)} h^{(s - 1)} \partial_0^{(s)} = h^{(s - 2)} \partial_0^{(s - 1)}\partial_0^{(s)} = h^{(s - 2)} d^{(s)} h^{(s)} = \partial_1^{(s - 3)} \partial_1^{(s - 2)} h^{(s)}.
  \]
  The $\tau$ component of the right-hand side is $\bar\partial^{(s - 3)} \bar\partial^{(s - 2)} h^{(s)}_\tau$, which vanishes for any choice of $h^{(s)}_\tau$. So we are done.
\end{proof}

In practice, we are rarely provided with a description of $M$ itself, but just the Steenrod module $\bar{M} \equiv \pi_0 M$. In this case, the strategy is to lift the chain complex $\bar{P}^{(\bullet)}$ itself, without the augmentation to $\bar{M}$. The above argument applies for most of the chain complex, except at the very beginning, where the chain complex is no longer exact.

Thus, we need to carefully choose $h^{(2)}_\tau$ such that the lift $h^{(3)}_\tau$ can be made. This is not always possible, and one checks that the obstruction to doing so lives in $\Ext^{3, 1}_\A(\bar{M}, \bar{M})$. If we manage to choose such an $h^{(2)}_\tau$, we can perform the rest of the inductive procedure and obtain a secondary chain complex $P^\bullet$. Finally, the spectral sequence implies that $\Tot(P^\bullet)$ is a lift of $\bar M$ to a periodic $\B$-module.

In general, different choices of $h^{(2)}_\tau$ will lead to different lifts of $\bar M$, and one can check directly that the set of choices is a torsor over $\Ext^{2, 1}_\A(\bar{M}, \bar{M})$. In practice, we often work in the case where there is a unique lift, and further, for degree reasons, \emph{any} choice of $h^{(2)}_\tau$ allows for a lift. Then we simply choose $h^{(2)}_\tau = 0$.

\begin{remark}
  When implementing this algorithm, the most expensive steps are evaluating $\partial_0^{(s - 2)} \partial_0^{(s - 1)}$, and then applying $A$ to it. Both of these steps are fully parallelizable with no data dependencies, so can be computed in a scalable and distributed fashion. Note that the composite $\partial_0^{(s - 2)} \partial_0^{(s - 1)}$ will be used again when computing products, and should not be discarded after applying $A$.
\end{remark}

\section{Computing the Adams \texorpdfstring{$E_3$}{E3} page}\label{section:e3}
Now that we can compute cofibrant replacements, we can, in theory, calculate anything we want in $\Mod_{C\tau^2}^{ft}$. In this \namecref{section:e3}, we will describe in detail the procedure for computing various useful sets of data. To simplify the presentation, we shall focus on obtaining data useful for computing the Adams spectral sequence for $\pi_{*, *} X$. The general case of computing $[X, Y]$ is not too much more difficult, but involves slightly more linear algebra\footnote{While it is true that $[X, Y] = \pi_* DX \otimes Y$ when $X$ is finite, this trick is not so useful in practice, since we lose the composition product structure under this transformation.}.

In this \namecref{section:e3}, we will assume that the free $\A$-resolution $\bar{P}^\bullet \to H^*X$ is a minimal resolution, i.e.\ $\bar{P}^\bullet \otimes_\A k$ has trivial differential. In this case, $\Hom_\A(\bar{P}^\bullet, k)$ also has trivial differential, so an element in $\Ext^{s, t}_\A(H^*X, k)$ is exactly a map $\bar{P}^{(s)} \to k[t]$. Thus, a choice of generators of each $\bar{P}^{(s)}$ also grants $\Ext^{s, t}_\A(H^*X, k)$ an $\F_p$-vector space basis.

\begin{remark}
  This \namecref{section:e3} is better seen as a technical documentation of the algorithm rather than a piece of mathematical exposition. Most of the content involves explicitly writing out large formulas that have to be satisfied and then observing that we can indeed iteratively make choices to satisfy these equations.
\end{remark}

\subsection{Computing \texorpdfstring{$d_2$}{d2}}\label{section:d2}
We begin by computing the $d_2$ differential in the Adams spectral sequence of $X$.

\begin{lemma}
  We can read off the $d_2$ differential of $X$ from a minimal free resolution of $\H^{(2)} X$.
\end{lemma}

\begin{proof}
  Recall that the $d_2$ differential is given by the connecting homomorphism of the cofiber sequence
  \[
    \Sigma^{0, -1} C\tau \otimes \nu X \to C\tau^2 \otimes \nu X \to C\tau \otimes \nu X
  \]
  in $\Syn_\HFp$. More precisely, it is obtained by applying $\pi_{*, *}$ to the connecting homomorphism of this cofiber sequence.

  By \Cref{lemma:unique-lift}, the cofiber sequence lifts uniquely to a sequence in $\Mod_{C\tau^2}$, and then the $d_2$ differential is induced by applying $[C\tau^2, -]^{*, *}_{\Mod_{C\tau^2}}$ to the connecting homomorphism.

  By \Cref{lemma:shift-ctau}, this is equivalent to applying $[\H^{(2)}X, -]^{*, *}$ to the sequence
  \[
    \Sigma k[1] \to \kk \to k.
  \]
  This connecting homomorphism can be computed as in the proof of the snake lemma, i.e.\ as the obstruction to lifting a map $\H^{(2)}X \to k$ along the projection $\kk \to k$.

  Explicitly, let $P^{\bullet} \to \H^{(2)}X$ be a minimal free secondary resolution lifting $\bar{P}^\bullet \to H^*X$. An element in $\Ext^{s, t}_\A(H^*X, k)$ is represented by a map $x\colon \bar{P}^{(s)} \to k[t]$, which lifts to a map $\tilde{x} \colon P^{(s)} \to \kk[t]$ using the section $\sigma^\kk$ we have chosen.

  We can try to make this into a map of chain complexes
  \[
    \begin{tikzcd}[column sep=huge, row sep = huge, ampersand replacement=\&]
      P_0^{(s + 2)} \oplus P_1^{(s + 1)} \ar[d, "{\begin{pmatrix}\partial^{(s + 2)}_0 & d^{(s + 1)} \\ -h^{(s + 2)} & -\partial^{(s + 1)}_1 \end{pmatrix}}"'] \ar[r] \& 0 \ar[d]\\
      P_0^{(s + 1)} \oplus P_1^{(s)} \ar[d, "{\begin{pmatrix}\partial^{(s + 1)}_0 & d^{(s)} \\ -h^{(s + 1)} & -\partial^{(s)}_1 \end{pmatrix}}"'] \ar[r, "{\begin{pmatrix} 0 & x_1 \end{pmatrix}}"] \& \kk_1[t] \ar[d, "d^\kk"] \\
      P_0^{(s)} \oplus P^{(s - 1)}_1 \ar[r, "{\begin{pmatrix} x_0 & 0 \end{pmatrix}}"] \& \kk_0[t]
    \end{tikzcd}
  \]
  where $x_0$ and $x_1$ are the components of $\tilde{x}$. The minimality of the resolution ensures that $x_0 \partial_0^{(s + 1)}$ is literally zero as a map of modules, so the bottom square commutes.\footnote{If the resolution is not minimal, then we get a non-trivial map $P_0^{(s + 1)} \to \kk_1[t]$ representing a null-homotopy of $x_0 \partial_0^{(s + 1)}$, and we have to adjust the upcoming argument accordingly.}

  On the other hand, the top square need not commute. Again, minimality ensures the $P_1^{(s + 1)}$ component of the map is automatically fine, and the $P_0^{(s + 2)}$ component is the map
  \[
    -x_1 h^{(s + 2)} \colon P_0^{(s + 2)} \to \kk_1[t].
  \]
  The commutativity of the bottom square ensures this is in the kernel of $d^\kk$, so it factors through $\pi_1 \kk[t] = k[t + 1]$. Concretely, this map is given by $-x h^{(s + 2)}_\tau \colon P_0^{(s + 2)} \to k[t + 1]$, which represents an element in $\Ext^{s + 2, t + 1}_\A(H^*X, k)$.

  This obstruction is exactly the connecting homomorphism of the cofiber sequence $\Sigma k[1] \to \kk \to k$. Hence $d_2(x) = x h_\tau^{(s + 2)}$ in the Adams spectral sequence.
\end{proof}

When working with synthetic spectra, computing the $E_3$ page is more than computing the $d_2$ differential, which is the connecting homomorphism of
\[
  \Sigma^{0, -1} C\tau \otimes \nu X \to C\tau^2 \otimes \nu X \to C\tau \otimes \nu X.
\]
Instead, we want to compute $\pi_{*, *} C\tau^2 \otimes \nu X$. Since we have already computed the connecting homomorphism, it remains to solve the extension problem. Fortunately, this is reasonably straightforward --- multiplication by $p$ is detected by $h_0$.

To state our result, we recall the \emph{carrying cocycle} \cite{arithmetic}:
\begin{notation}
  The \emph{carrying cocycle} $x \mathbin{\tilde{+}} y$ of $\sigma^\kk$ is defined by
  \[
    \sigma^\kk(x) + \sigma^\kk(y) = \sigma^\kk(x + y) + p (x \mathbin{\tilde{+}} y)
  \]
  for $x, y \in \F_p$. This naturally extends to a function $V \times V \to V$ for any $\F_p$-vector space $V$ with a basis.
\end{notation}

\begin{thm}\label{thm:lift}
  Fix a minimal free resolution, hence a basis of $\Ext^{s, t}(H^*X, k)$ for every $s, t$. For every element $x \in \Ext^{s, t}(H^* X, k)$ such that $d_2(x) = 0$, there is a standard lift $[x] \in \pi_{*, *} C\tau^2 \otimes \nu X$ with the property that
  \[
    [x + y] = [x] + [y] + \tau h_0 (x \mathbin{\tilde{+}} y),
  \]
\end{thm}
This completely specifies the additive structure of $\pi_{*, *} C\tau^2 \otimes \nu X$.

\begin{proof}
  Let $x$ be an element that survives the $E_2$ page. Then the argument above describes a lift of $x$ to an element in $[\H^{(2)}X, \kk]^{s, t}$, or equivalently, $\pi_{t - s, s} C\tau^2 \otimes \nu X \otimes$. We call this element $[x]$, which is a well-defined lift once we fixed every choice we have made so far.

  This is, of course, not the only lift. If $y \in \Ext^{s + 1, t + 1}_\A(H^*X, k)$, then we can add it to the $\tau$ component of the map $P_0^{(s + 1)} \to \kk_1[t]$ that we originally picked to be zero. It would still be a chain map, and this represents $[x] + \tau y$.

  The failure of $[-]$ to be additive arises from the fact that the sum of two standard lifts need not be a standard lift, since $\sigma^\kk$ is not linear. To prove the additive relation claimed, we have to show that if $x \colon \bar{P}^{(s)} \to k[t]$ is any map, then the chain map $P^{\bullet} \to \kk[t]$ given by $p \tilde{x}$ is homotopic to $\tau h_0 x$. This is a straightforward computation, with the homotopy being given by $\mu_0 x_0$. Ultimately, this boils down to the relation $\beta \mu_0 = \tau$.
\end{proof}

\subsection{Computing products}
We now turn to computing composition products in $\Mod_{C\tau^2}$. To simplify matters, we shall only consider the case of computing the $\pi_{*, *} C\tau^2$ action on $\pi_{*, *} C\tau^2 \otimes \nu X$.

Let $P^\bullet \to \H^{(2)} X$ be a minimal free resolution of $\H^{(2)} X$, and let $Q^\bullet \to \kk$ be a minimal free resolution of $\kk$. From the previous \namecref{section:d2}, we know that an element $f\in \pi_{t - s, s} C\tau^2 \otimes \nu X$ is represented by a chain map $\Tot(P^\bullet) \to \Sigma^s \kk[t]$, which we can lift to a map $\Tot(P^\bullet) \to \Sigma^s \Tot(Q^\bullet)[t]$ since the source is cofibrant.

Now an element in $\pi_{t' - s', s'} C\tau^2$ is represented by a chain map $\Tot(Q^\bullet) \to \Sigma^{s'} \kk[t']$. To compute the product with $f$, we simply compose this with the lifted chain map $\Tot(P^\bullet) \to \Sigma^s \Tot(Q^\bullet)[t]$ and read off the composite.

Most of the hard work is in actually writing down the lift of $f$ to $\Sigma^s \Tot(Q^\bullet)[t]$. To simplify notation slightly, we shift $X$ so that $t = 0$. We begin by computing a lift of $f$ over the ordinary Steenrod algebra, i.e.\ a lift to a map $\bar{f} \colon \bar{P}^\bullet \to \Sigma^s \bar{Q}^\bullet$. Our standard splittings then let us lift this to a diagram
\[
  \begin{tikzcd}
    \vdots \ar[d, "\partial^{(s + 2)}"] \ar[r] & \vdots \ar[d, "\partial^{(2)}"] \\
    P^{(s + 1)} \ar[d, "\partial^{(s + 1)}"] \ar[r, "f^{(s + 1)}"] & Q^{(1)} \ar[d, "\partial^{(1)}"] \\
    P^{(s)} \ar[r, "f^{(s)}"] \ar[d]  & Q^{(0)} \ar[d] \\
    P^{(s - 1)} \ar[r] \ar[d] & 0 \ar[d] \\
    \vdots \ar[r] & \vdots
  \end{tikzcd}
\]
By construction, each of these squares commute up to homotopy, and again our job is to find a suitable homotopy $H$ that induces a chain map on $\Tot(-)$. To do so, we write down the induced map on the total chain complex:
\[
  \begin{tikzcd}[column sep=2.8cm, row sep = huge, ampersand replacement=\&]
    \vdots \ar[d] \ar[r] \& \vdots \ar[d] \\
    P_0^{(s + k)} \oplus P_1^{(s + k - 1)} \ar[d, "{\begin{pmatrix}\partial^{(s + k)}_0 & d^{(s + k - 1)} \\ -h^{(s + k)} & -\partial^{(s + k - 1)}_1 \end{pmatrix}}"'] \ar[r, "{\begin{pmatrix} f^{(s + k)}_0 & 0 \\ H^{(s + k)} & f^{(s + k - 1)}_1 \end{pmatrix}}"] \& Q_0^{(k)} \oplus Q_1^{(k - 1)} \ar[d, "{\begin{pmatrix}\partial^{(k)}_0 & d^{(k - 1)} \\ -h^{(k)} & -\partial^{(k - 1)}_1 \end{pmatrix}}"]\\
    P_0^{(s + k - 1)} \oplus P^{(s + k - 2)}_1 \ar[d] \ar[r, "{\begin{pmatrix} f^{(s + k - 1)}_0 & 0 \\ H^{(s + k - 1)} & f^{(s + k - 2)}_1 \end{pmatrix}}"] \& Q_0^{(k - 1)} \oplus Q^{(k - 2)}_1 \ar[d]\\
    \vdots \ar[r] \& \vdots
  \end{tikzcd}
\]
Expanding the matrices gives the equations
\[
  \begin{aligned}
    d^{(k)} f_1^{(s + k)} &= f_0^{(s + k)} d^{(s + k)} \\
    H^{(s + k)} d^{(s + k)} &= f_1^{(s + k - 1)} \partial_1^{(s + k)} - \partial_1^{(k)} f_1^{(s + k)} \\
    d^{(k - 1)} H^{(s + k)} &= f_0^{(s + k - 1)} \partial^{(s + k)}_0 - \partial_0^{(k)} f_0^{(s + k)} \\
    \partial_1^{(k - 1)} H^{(s + k)} + H^{(s + k - 1)} \partial^{(s + k)}_0 &= f_1^{(s + k - 2)} h^{(s + k)} - h^{(k)} f_0^{(s + k)}
  \end{aligned}
\]
As usual, the first two are implied by the remaining via $\B$-linearity, and the third equation simply says $H^{(k)}$ is a homotopy between $\partial^{(k)} f^{(k)}$ and $f^{(k - 1)} \partial^{(k)}$. The main equation of content is the last one. We can interpret this as follows --- there are two natural null-homotopies of $f^{(s + k - 2)} \partial^{(s + k - 1)} \partial^{(s + k)}$:
\begin{itemize}
  \item we can use the null-homotopy of $\partial^{(s + k - 1)} \partial^{(s + k)}$ and compose it with $f^{(s + k - 2)}$; or
  \item we can homotope to $\partial^{(k - 1)} \partial^{(k)} f^{(s + k)}$ and apply the null-homotopy of $\partial^{(k - 1)} \partial^{(k)}$.
\end{itemize}
The equation asserts that these two null-homotopies agree.

Regardless of how we are supposed to interpret this equation, our job is to choose the $\tau$ part of $H$ so that the equation is satisfied.

Let $g$ be a generator of $P^{(s + k)}$. Write
\[
  \bar{\partial} (g) = \sum \alpha^i g_i,\quad \bar{f} (g) = \sum \beta^j g_j.
\]
Then the $\tau$ part of the last equation says
\begin{multline*}
\bar\partial^{(k - 1)} H^{(s + k)}_\tau g \\
+ \sum \left(\alpha^i H^{(s + k - 1)}_\tau g_i + A\left(\alpha^i, \left(f_0^{(s + k - 1)} \partial^{(s + k)}_0 - \partial_0^{(k)} f_0^{(s + k)}\right)g_i\right)\right) \\
= \bar{f}^{(s + k - 2)} h^{(s + k)}_\tau g - \sum \left(\beta^j h_\tau^{(k)} g_j + A\left( \beta^j, \partial_0^{(k - 1)} \partial^{(k)}_0 g_j\right)\right).
\end{multline*}

We again solve this inductively. The first homotopy is $H^{(s + 1)}$. Its equation takes values in the zero group, so it always holds. Changing the $\tau$ part modifies $f$ by a $\tau$-multiple. If we want to set $f = [\bar{f}]$, then we choose the $\tau$ part of $H^{(s + 1)}$ to vanish.

As for the second homotopy $H^{(s + 2)}$, we have chosen the target of the chain map to be a resolution of $\kk$, and we always choose minimal resolutions. So $\bar\partial^{(1)}$ is surjective in positive internal degree. The only term in the equation that is in internal degree zero is the $\bar{f}^{(s + k -2)} h_\tau^{(s + k)}g$ term. Requiring this to vanish is exactly the requirement that $\bar{f}$ survives the Adams $d_2$. (Note that in $\beta^j h_\tau^{(2)} g_j$, the $h^{(2)}$ is the homotopy of the free resolution of $\kk$, which is zero for degree reasons)

Afterwards, exactness implies that the lift can always be performed. To see this, as in the case of constructing a free resolution, we have to check that $\partial_1^{(k - 2)}$ applied to the last equation is always satisfied. Instead of painstakingly tracking through each of the terms, it suffices to observe that by induction, the equations desired always hold after applying the differential in $\Tot(Q^\bullet)$. Since the first three equations always hold, we know that $\partial_1^{(k - 2)}$ must kill the last equation.

Once we have performed this lift, given any class $x \in \Ext^{s', t'}(k, k)$ that survives the $d_2$ differential, we can compute the product $[x]f$ by composing the chain maps. We then end up with
\[
  \begin{tikzcd}[column sep=4cm, ampersand replacement=\&]
    P_0^{(s + s' + 1)} \oplus P_1^{(s + s')} \ar[d] \ar[r, "{\begin{pmatrix} x_1 H^{(s + s' + 1)} & x_1 f_1^{(s + s')}\end{pmatrix}}"] \& \kk_1[t'] \ar[d, "d^\kk"]\\
    P_0^{(s + s')} \oplus P^{(s + s' - 1)}_1 \ar[r, "{\begin{pmatrix} x_0 f_0^{(s + s')} & 0 \end{pmatrix}}"] \& \kk_0[t']
  \end{tikzcd}
\]
One has to be extremely careful here, as $x_0 f_0^{(s + s')}$ need not be the standard lift of $x\bar{f}$; one has to use the formula $p = \tau h_0$ to express the result in terms of $[x\bar{f}]$.

The interesting part is, of course, the $x_1 H^{(s + s' + 1)}$ term. Note that this term lives in the kernel of $d^\kk$. Indeed, by the commutativity of the diagram, $d^\kk x_1 H^{(s + s' + 1)} = x_0 f_0^{(s + s')} \partial_0^{(s + 1)}$, but $\partial_0^{(s + 1)}$ only hits decomposables by minimality, which is then killed by $x_0$. The $\tau$ component of $x_1 H^{(s + s' + 1)}$ is $x H^{(s + s' + 1)}_\tau$, so the product $[x] f$ is given by whatever $x_0 f_0$ represents plus $\tau (x H^{(s + s' + 1)}_\tau)$.

In a very imprecise manner, this suggests that $H^{(s + s' + 1)}_\tau$ encodes the hidden extension part of the product. Of course, ``the hidden extension part'' only makes sense after choosing the standard lifts $[-]$; it is not a homotopically meaningful concept.

\subsection{Massey products}
We finally turn to Massey products, which requires lifting chain homotopies. This actually contains two kinds of information. Firstly, we get to compute hidden extensions in Massey products that jump by one filtration, which is always useful. But we also get to compute Massey products of the form $\langle x, \tau y, z \rangle$ when $xy$ and $yz$ are not zero in the $E_2$ page but are hit by a differential. This is exactly the $E_3$ page Massey product as described in Moss' convergence theorem \cite{moss}. While these are extremely easy to compute in the $E_3$ page, our synthetic approach gives us the answer up to mod $\tau^2$ instead of mod $\tau$, which is extra information one can capitalize on. (It also provides a very useful test case for the algorithm, since we can verify the answers by hand)

This is more involved than the previous cases, and is largely unpleasant. In particular, we apologize in advance for the overwhelming number of objects called $h$.

We begin with a general remark on chain homotopies, which is rather important. Suppose we have a chain map between (ordinary) chain complexes over $\A$ of the form
\[
  \begin{tikzcd}
    \vdots \ar[d, "\partial^{(s + 2)}"] \ar[r] & \vdots \ar[d, "\partial^{(2)}"] \\
    P^{(s + 1)} \ar[d, "\partial^{(s + 1)}"] \ar[r, "f^{(s + 1)}"] & Q^{(1)} \ar[d, "\partial^{(1)}"] \\
    P^{(s)} \ar[r, "f^{(s)}"] \ar[d]  & Q^{(0)} \ar[d] \\
    P^{(s - 1)} \ar[r] \ar[d] & 0 \ar[d] \\
    \vdots \ar[r] & \vdots
  \end{tikzcd}
\]
Null-homotopies of this chain map consist of maps $\H^{(s + k)} \colon P^{(s + k)} \to Q^{(k + 1)}$ satisfying the equation
\[
  \partial^{(k + 1)} \H^{(s + k)} + \H^{(s + k - 1)} \partial^{(s + k)} = f^{(s + k)}.
\]
The first possible non-zero map is $\H^{(s - 1)}$. A useful trick is that if these are minimal resolutions and $Q^\bullet$ resolves $k$, then we can always choose $\H^{(s - 1)} = 0$. Indeed, since $f$ is null and the resolution is minimal, $f^{(s)}$ can only hit positive degree elements. Since $\partial^{(1)}$ is surjective in positive internal degree, the lifting problem for $H^{(s)}$ can always be solved. Subsequent lifting problems can then be solved by exactness.

On the other hand, this is not true in general, and in particular not true when working in the secondary setting.

Indeed, over the ordinary Steenrod algebra, the set of null-homotopies is a torsor over the appropriate $\Ext$ group, which acts by modifying $\H^{(s - 1)}$. While we can always choose $\H^{(s - 1)} = 0$, this should not be thought of as having a special status.

In the secondary setting, it is no longer a torsor over the full $\Ext$ group; we can only add elements that survive the Adams $d_2$. Then the homotopy with $\H^{(s - 1)} = 0$ need not be a valid basepoint of this torsor. Part of our job when lifting chain homotopies is to find a choice of $\H^{(s - 1)}$.

Let $P^\bullet \to \H^{(2)} X$ and $Q^\bullet \to \kk$ be minimal free resolutions. Suppose we have a chain map between them of the form
\[
  \begin{tikzcd}[column sep=3cm, row sep = huge, ampersand replacement=\&]
    P_0^{(s + k + 1)} \oplus P_1^{(s + k)} \ar[d, "{\begin{pmatrix}\partial^{(s + k + 1)}_0 & d^{(s + k)} \\ -h^{(s + k + 1)} & -\partial^{(s + k)}_1 \end{pmatrix}}"'] \ar[r, "{\begin{pmatrix} f^{(s + k + 1)}_0 & 0 \\ H^{(s + k + 1)} & f^{(s + k)}_1 \end{pmatrix}}"] \& Q_0^{(k + 1)} \oplus Q_1^{(k)} \ar[d, "{\begin{pmatrix}\partial^{(k + 1)}_0 & d^{(k)} \\ -h^{(k + 1)} & -\partial^{(k)}_1 \end{pmatrix}}"]\\
    P_0^{(s + k)} \oplus P^{(s + k - 1)}_1 \ar[r, "{\begin{pmatrix} f^{(s + k)}_0 & 0 \\ H^{(s + k)} & f^{(s + k - 1)}_1 \end{pmatrix}}"] \& Q_0^{(k)} \oplus Q^{(k - 1)}_1
  \end{tikzcd}
\]
that is null-homotopic. Then the induced map on $\A \otimes_\B(-)$ is also null-homotopic and admits null-homotopies, say $\H^{(s + k)} \colon P^{(s + k)} \to Q^{(k + 1)}$. These lift to maps
\[
  \begin{pmatrix}
    \H^{(s + k)}_0 & 0\\
    -\eta^{(s + k)} & -\H^{(s + k - 1)}_1
    \end{pmatrix}
    \colon
    P_0^{(s + k)} \oplus P_1^{(s + k - 1)} \to Q_0^{(k + 1)} \oplus Q_1^{(k)}.
\]
Expanding the definition of a chain homotopy gives the equations
\[
  \begin{aligned}
    d^{(k)} \eta^{(s + k)} &= \partial_0^{(k + 1)} \H^{(s + k)}_0 + \H_0^{(s + k - 1)} \partial_0^{(s + k)} - f_0^{(s + k)}\\
    \partial_1^{(k)} \eta^{(s + k)} &= \eta^{(s + k - 1)} \partial_0^{(s + k)} + h^{(k + 1)} \H_0^{(s + k)} - \H_1^{(s + k - 2)} h^{(s + k)} - H^{(s + k)}.
  \end{aligned}
\]
The first equation simply states that $\eta$ witnesses the equation $\partial \H + \H \partial \simeq f$, and the second is some complex compatibility condition that we have to solve iteratively.

Note that the equation for $\eta^{(s + k)}$ takes values in $Q_1^{(k - 1)}$. Thus, it is automatically satisfied for \emph{both} $k = -1$ and $k = 0$. The $k = 1$ step is the only step where it might be impossible, and all higher $k$ can be solved by exactness.

When $k = 1$, since $\bar{\partial}^{(1)}$ is surjective in positive internal degrees, we can only fail to lift in internal degree $0$. We shall now analyze which terms can contribute to the $\tau$ part.
\begin{itemize}
  \item In $\eta^{(s)} \partial_0^{(s + 1)}$, the only possible contribution comes from the $f_0$ part of $\eta$. The idea is that minimality ensures $\partial_0$ introduces an algebra element of positive degree. On the other hand, when multiplying with a homotopy, the $A$ function shows up, and lowers degree by 1. So as long as two $\partial_0$'s show up, we are clear.
  \item The $h^{(2)} \H_0^{(s + 1)}$ term cannot contribute. Indeed, $\H_0^{(s + 1)}$ takes values in $Q_0^{(2)}$, whose generator of lowest degree is $2$. While $h^{(2)}$ can lower degree by $1$, we still have $2 - 1 > 0$.
  \item The $\H_1^{(s - 1)} h^{(s + 1)}$ and $H^{(s + 1)}$ terms can contribute.
\end{itemize}

Our goal is then to choose $\H_1^{(s - 1)}$ appropriately so that the right-hand side vanishes in degree $0$. This has a very natural interpretation. In degree $0$, The map $f_0$ is necessarily $p$ times some $\Ext$ class, and one checks that the term $\eta^{(s)} \partial_0^{(s + 1)}$ picks out $h_0$ times said $\Ext$ class. The term $H^{(s + 1)}$ is the $\tau$ part of $f$ itself. So these two terms combined give us the $\tau$ part of $f$ after normalizing the degree $0$ part to exactly $0$. The equation then tells us $\H_1^{(s - 1)}$ should be the class whose $d_2$ kills the $\tau$ part of $f$, witnessing the fact that $f$ is indeed null.

In general, finding this $\H_1^{(s - 1)}$ is low-dimensional linear algebra, and is relatively easy. While the equation it has to satisfy comes from the $k = 1$ case, it ends up not depending on the $k = -1$ and $k = 0$ data, so it can be computed before we start the lifting process.

Once we manage to lift chain homotopies, computing Massey products becomes relatively straightforward. Unfortunately, the sign conventions surrounding Massey products are rather confusing and not well-documented in the literature. Instead of tackling this problem, we are content with computing Massey products up to a sign.

\section{The computer implementation}\label{section:data}
\subsection{The generated data}
We ran our algorithm on the sphere at the prime $2$ up to the $140$\textsuperscript{th} stem. The output of the algorithm is available at \cite[\texttt{d2-data.zip}]{raw_data}, and the contents of each file are described in \Cref{table:data}. In the table, the first three sets of files are the data generated by the new algorithm, while the rest are auxiliary data to assist the user in interpreting the data.

\begin{table}[ht]
  \centering
  \caption{List of data files}\label{table:data}
  \vspace{0.4cm}
  \begin{tabularx}{\textwidth}{lX}
    \toprule
    Filename & Description\\
    \midrule
    \texttt{d2} & $d_2$ differentials. \\
    \texttt{product\_a} & The product of all elements with \texttt{[a]}. We have computed products with all indecomposables up to the $39$\textsuperscript{th} stem, and the names are listed in \Cref{table:prod-name}. Note that these products include the twist of $(-1)^{s't}$.  \\
    \texttt{massey\_a\_b} & The Massey product \texttt{<-, [b], [a]>} up to a sign. The only exception is \texttt{massey\_P}, which contains the Adams periodicity operator $\langle -, [h_0^4], [h_3]\rangle$. \\
    \midrule
    \texttt{change\_of\_basis} & Translation between our basis and the basis of the Bruner--Rognes dataset \cite{bruner-rognes}. \\
    \texttt{filtration\_one} & All ($E_2$ page) filtration one products. This is useful for identifying elements by hand. \\
    \texttt{charts.pdf} & Adams charts displaying the $E_2$ and $E_3$ pages. When a bidegree has more than one basis element, they are ordered bottom-to-top, left-to-right.\\
    \texttt{clean\_charts.pdf} & The same charts as above but without $h_2$ products. \\
    \texttt{differentials.gz} & Differentials in our minimal resolution. This contains information to uniquely identify all of our basis elements and lifts, but is most likely not of much use to humans. \\
    \bottomrule
  \end{tabularx}
\end{table}

\begin{table}[ht]
  \centering
  \caption{Names of products}\label{table:prod-name}
  \vspace{0.4cm}
  \begin{tabular}{cccc}
    \toprule
    $n$ & $s$ & class & name \\
    \midrule
    0 & 1 & $[1]$ & \verb|h_0| \\
    1 & 1 & $[1]$ & \verb|h_1| \\
    3 & 1 & $[1]$ & \verb|h_2| \\
    7 & 1 & $[1]$ & \verb|h_3| \\
    8 & 3 & $[1]$ & \verb|c_0| \\
    9 & 5 & $[1]$ & \verb|Ph_1| \\
    11 & 5 & $[1]$ & \verb|Ph_2| \\
    14 & 4 & $[1]$ & \verb|d_0| \\
    15 & 2 & $[1]$ & \verb|h_0h_4| \\
    16 & 2 & $[1]$ & \verb|h_1h_4| \\
    16 & 7 & $[1]$ & \verb|Pc_0| \\
    17 & 9 & $[1]$ & \verb|P^2h_1| \\
    18 & 2 & $[1]$ & \verb|h_2h_4| \\
    19 & 3 & $[1]$ & \verb|c_1| \\
    19 & 9 & $[1]$ & \verb|P^2h_2| \\
    20 & 4 & $[1]$ & \verb|g| \\
    22 & 8 & $[1]$ & \verb|Pd_0| \\
    23 & 4 & $[1]$ & \verb|h_4c_0| \\
    23 & 9 & $[1,1]$ & \verb|h_0^2i| \\
    24 & 11 & $[1]$ & \verb|P^2c_0| \\
    25 & 13 & $[1]$ & \verb|P^3h_1| \\
    27 & 13 & $[1]$ & \verb|P^3h_2| \\
    30 & 2 & $[1]$ & \verb|h_4^2| \\
    30 & 6 & $[1]$ & \verb|Dh_2^2| \\
    \bottomrule
  \end{tabular}
  \hspace{1em}
  \begin{tabular}{cccc}
    \toprule
    $n$ & $s$ & class & name \\
    \midrule
    30 & 12 & $[1]$ & \verb|P^2d_0| \\
    31 & 4 & $[1]$ & \verb|h_0^3h_5| \\
    31 & 5 & $[0,1]$ & \verb|n| \\
    31 & 8 & $[1,1]$ & \verb|d_0e_0| \\
    32 & 2 & $[1]$ & \verb|h_1h_5| \\
    32 & 4 & $[1]$ & \verb|d_1| \\
    32 & 6 & $[1]$ & \verb|Dh_1h_3| \\
    32 & 15 & $[1]$ & \verb|P^3c_0| \\
    33 & 4 & $[1]$ & \verb|p| \\
    33 & 17 & $[1]$ & \verb|P^4h_1| \\
    34 & 2 & $[1]$ & \verb|h_2h_5| \\
    34 & 8 & $[1]$ & \verb|e_0^2| \\
    35 & 17 & $[1]$ & \verb|P^4h_2| \\
    36 & 6 & $[1]$ & \verb|t| \\
    37 & 5 & $[1]$ & \verb|x| \\
    37 & 8 & $[0,1]$ & \verb|e_0g| \\
    38 & 2 & $[1]$ & \verb|h_3h_5| \\
    38 & 4 & $[1,0]$ & \verb|e_1| \\
    38 & 16 & $[1]$ & \verb|P^3d_0| \\
    39 & 4 & $[1]$ & \verb|h_5c_0| \\
    39 & 9 & $[1]$ & \verb|Dh_1d_0| \\
    39 & 12 & $[1]$ & \verb|Pd_0e_0| \\
    39 & 17 & $[1,1]$ & \verb|h_0^2P^2i| \\
    \\
    \bottomrule
  \end{tabular}
\end{table}

In all files, the results are expressed in terms of our $E_2$ page basis. We adopt the following naming conventions:
\begin{enumerate}
  \item \verb|x_(n, s, i)| is the $i$\textsuperscript{th} basis element in filtration $s$ of the $n$\textsuperscript{th} stem.
  \item We use $[-]$ to denote the standard lift to $\Mod_{C\tau^2}$ as in \Cref{thm:lift}. Note again that $[a + b] \not= [a] + [b]$ in general.
  \item If an element is in a known degree (e.g.\ it is the value of a product), we will write an element in vector form under our basis, e.g.\ as $[1, 0]$. We shall not put an extra pair of brackets around the vector to denote the secondary lift. It should be clear from context whether we mean the $E_2$ page element or its secondary lift.
  \item We use $\tau$ to denote multiplication by $\tau$ (our files are UTF-8 encoded).
\end{enumerate}

\subsection{Generating the data}
The code used for the calculation is available at \cite{ext_rs}, and the latest version of this software is available at \url{https://github.com/SpectralSequences/sseq}. This is a monorepo, and we will work in the \texttt{ext/} subdirectory throughout. This repository comes with a reasonable amount of documentation, and the \texttt{README} in \texttt{ext/} contains instructions for accessing said documentation.

The commands used to generate the data are packaged into a script, which is available at \cite[\texttt{script.sh}]{raw_data}. This should be run in the \texttt{ext/} directory of the repository. The save files for the computations are at \cite[\texttt{S\_2\_milnor.tar}]{raw_data}.

To assist the reader in further exploring the resolution, we illustrate the full interactive session that generates the data we need in \Cref{figure:session}. Assuming Rust is installed, running the commands as indicated in any subdirectory of \texttt{ext/} will compute the $d_2$ differentials for $S_2$, the product with $g$ as well as the Adams periodicity operator.\footnote{As mentioned in the documentation, when resolving to larger stems, one ought to supply the \texttt{-{}-release}, \texttt{-{}-features concurrent} and \texttt{-{}-no-default-features} flags after \texttt{cargo run} for much improved performance.} This guide is written for the version in \cite{ext_rs} but should work with future versions with little modifications.

\begin{figure}[ht]
\begin{shell}
 $ cargo run --example secondary > d2
\prompt{Module (default: S_2):} S_2
\prompt{Module save directory (optional):} S_2_milnor
\prompt{Max n (default: 30):} 40
\prompt{Max s (default: 15):} 20

 $ cargo run --example secondary_product > product_g
\prompt{Module (default: S_2):} S_2
\prompt{Module save directory (optional):} S_2_milnor
\prompt{Max n (default: 30):} 40
\prompt{Max s (default: 15):} 20
\prompt{Name of product:} g
\prompt{n of Ext class g:} 20
\prompt{s of Ext class g:} 4
\prompt{Input ext class:} [1]

 $ cargo run --example secondary_massey > massey_P
\prompt{We are going to compute <-, b, a> for all (-), where a is an}
\prompt{element in Ext(M, k) and b and (-) are elements in Ext(k, k).}
\prompt{Module (default: S_2):} S_2
\prompt{Module save directory (optional):} S_2_milnor
\prompt{Max n (default: 30):} 40
\prompt{Max s (default: 15):} 20
\prompt{n of a:} 7
\prompt{s of a:} 1
\prompt{Name of Ext part of a:} h_3
\prompt{Input Ext class h_3:} [1]
\prompt{Name of τ part of a:}
\prompt{n of b:} 0
\prompt{s of b:} 4
\prompt{Name of Ext part of b:} h_0^4
\prompt{Input Ext class h_0^4:} [1]
\prompt{Name of τ part of b:}
\end{shell}
\caption{Interactive session to generate the dataset. The grey text is the computer's prompt and the black text is the user's input.}\label{figure:session}
\end{figure}

\subsection{Runtime performance}
We ran the program on a computational server of the Harvard University Mathematics Department. It has two \texttt{Xeon E5-2690 v2} CPUs (10 cores/20 threads each) and 125 GiB of memory. Using all 40 threads of the machine, computing the secondary resolution up to the 140\textsuperscript{th} stem took 3.3 hours and 7.8 GiB of memory.\footnote{We will mostly focus on analyzing the performance of computing the secondary resolution itself. For the products, computing the $h_0$ product generally takes twice as long as computing the resolution, while computing $\langle -, [h_1], [h_0]\rangle$ took $2.6\times$ as long. The time taken decreases rapidly as the stem of the multiplicand increases. For example, if we want to compute the product with $[g]$, to stay within the range, we would only multiply $[g]$ with elements up to the 120\textsuperscript{th} stem, and we only need to lift the chain map for 120 stems, as opposed to 140 for $h_0$.}

To understand the asymptotic complexity of the algorithm, recall that for each generator, to compute $h_\tau^{(s)}g$, we have to solve the equation
\[
  \bar{\partial}^{(s - 2)} h_\tau^{(s)} g = \sum \left(\alpha^i h_\tau^{(s - 1)}g_i + A\left(\alpha^i, \partial_0^{(s - 1)} \partial_0^{(s - 2)} g_i\right)\right).
\]

We break this up into two steps --- we first compute $A\left(\alpha^i, \partial_0^{(s - 1)} \partial_1^{(s - 2)} g_i \right)$, and then solve the rest of the lifting problem. The first step is fully parallelizable, as there are no dependencies between different generators, while the second step requires the value of $h_\tau^{(s - 1)} g_i$, so must be computed in some specific order.

In practice, the second step is much faster than the first step even after parallelization. Moreover, the cost of the second step is exactly the cost of computing a single product; if it takes too long, we have bigger problems to deal with.

Thus, we shall focus on the cost of the first step. Our objective is to understand how this grows with the stem. There are two separate questions we can ask:
\begin{enumerate}
  \item What is the maximum time it takes to perform the first step for a single generator?
  \item How much time does it take in total to compute up to a stem?
\end{enumerate}

The first question is relevant in a situation where we have an extremely large number of cores/machines that can parallelize the computation, in which case the bottleneck is the slowest generator. The second question is relevant where we have a fixed, small(ish) number of cores that will be saturated throughout the process.

\begin{figure}[t]
  \centering
  \input{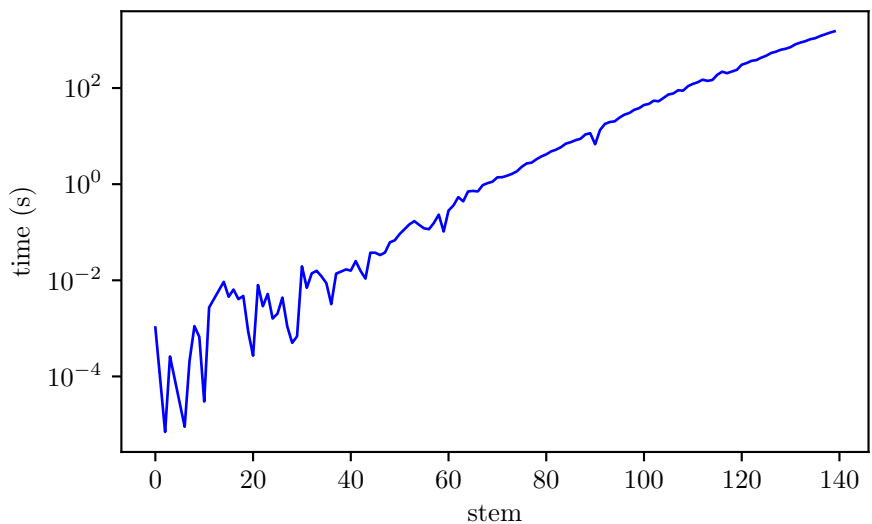}
  \caption{Time taken by the slowest generator in each stem}\label{figure:max}
\end{figure}

\begin{figure}[t]
  \centering
  \input{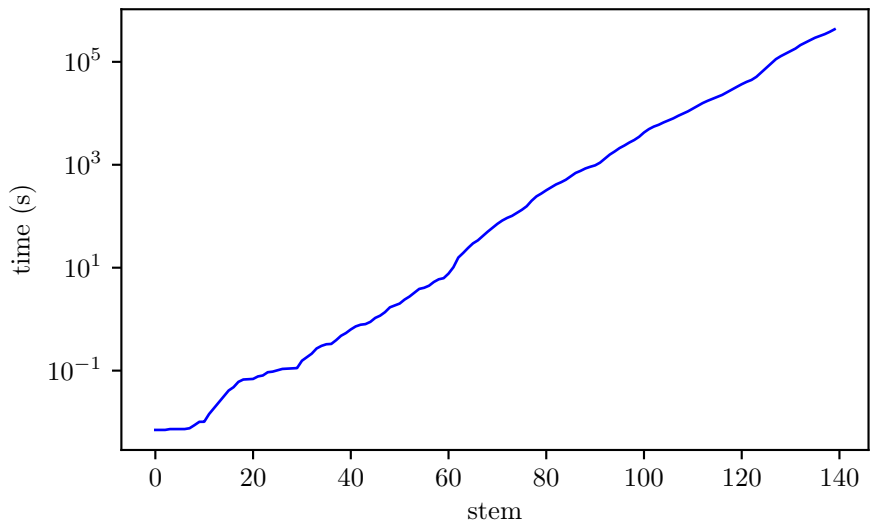}
  \caption{Time taken to compute up to each stem}\label{figure:cumulative}
\end{figure}

\begin{figure}[t]
  \centering
  \input{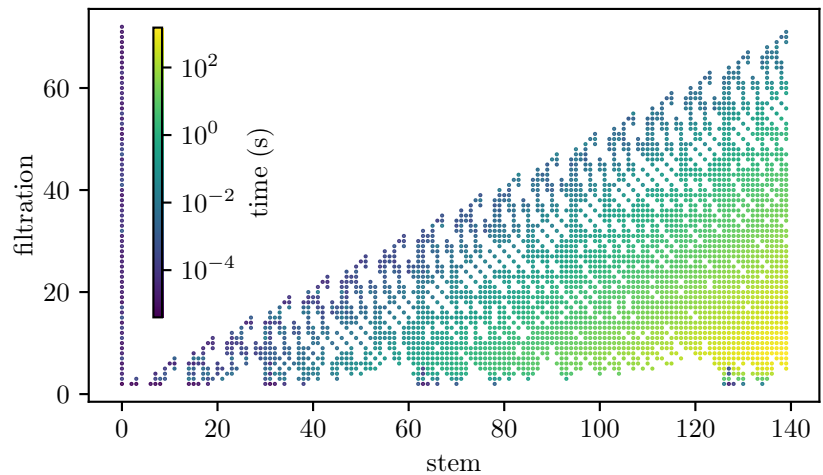}
  \caption{Time taken by the slowest generator in each bidegree}\label{figure:scatter}
\end{figure}

To answer these questions, we timed the computation for each generator and generated three charts:
\begin{enumerate}
  \item \Cref{figure:max} shows the time taken by the slowest generator in each stem.
  \item \Cref{figure:cumulative} shows the time taken to compute up to each stem.
  \item \Cref{figure:scatter} shows the time taken by the slowest generator in each bidegree.
\end{enumerate}
Again, these figures only include the time taken by the first part, and measure CPU time as opposed to wall time (so the actual time needed to compute up to a stem is around $\frac{1}{40}$ of the time indicated).

The most obvious feature that stands out is that the time tends to grow exponentially in stem (the time axis uses a log scale). Fitting a simple linear regression on the datapoints beyond the 50\textsuperscript{th} stem, we see that the maximum time increases by a factor of 3 every 10 stems, while the cumulative time increases by a factor of $~3.85$ every 10 stems. For comparison, the cumulative time of Nassau's algorithm for computing the minimal resolution increases by a factor of 2.8 every 10 stems \cite[Abbildung 2.13]{nassau}.

From the scatter plot, we see that for each stem, the slowest bidegrees are the ones with the lowest filtrations, except for a few exceptions in very low filtrations. This is expected, since the Adams vanishing line suggests there will be more and lower degree generators in low filtrations, hence the resolution is larger in these bidegrees. This also explains the irregularities one observes in the graphs. In \Cref{figure:max}, the dips coincide with the stems where there are no low filtration elements. Conversely, the jumps in the cumulative time occur in stems near $2^n$, where there is a higher density of low filtration elements.

\subsection{Future work}
There are a few obvious improvements one can make to the dataset:
\begin{enumerate}
  \item Compute further into more stems. This mostly requires more computational power. The ``secondary'' part of the process is fully parallelizable, and the code supports distributing the work across multiple machines.
  \item Compute \emph{all} products by indecomposables, not just the indecomposables up to the $39$\textsuperscript{th} stem. The remaining products are extremely fast to compute, since the cost depends on the stem of the multiplicand, which now only goes up to at most 100. The main bottleneck is enumerating the indecomposables, which has been done manually so far. To push the product computation further, we ought to automate the process of finding all indecomposables.
  \item Compute more Massey products. However, computing \emph{all} potential Massey products seems prohibitively expensive.
  \item Have a dataset expressed in terms of ``human names'' of the classes. The main blocker is in coming up with a reasonable database of names.
  \item Have a program to propagate differentials back and forth using the Leibniz rules and all available products.
\end{enumerate}

\part{Computation of Adams differentials}\label{part:computation}
\section{Overview}
The goal of this \namecref{part:computation} is to use the computer calculations of \Cref{part:secondary} to compute differentials in the Adams spectral sequence.

The computer algorithm automatically gives us all $d_2$ differentials. To compute longer differentials, we introduce the notion of a hidden extension on the $E_k$ page. Essentially by definition, hidden extensions on the $E_3$ page can be read off from the computer calculated $\Mod_{C\tau^2}$ composition products. Equipped with these hidden extensions, a generalized Leibniz rule then lets us relate differentials of different lengths.

After introducing this machinery in \Cref{section:diff-hidden}, we proceed to perform two sets of computations.

In \Cref{section:old-diff}, we compute the first 35 stems of the Adams spectral sequence. Of course, all of these results are well-known; the goal is to illustrate the techniques in more familiar territory.

In \Cref{section:new-diff}, we resolve previously unknown differentials in the Adams spectral sequence. In particular, we compute all unknown $d_2$, $d_3$, $d_4$ and $d_5$ differentials up to the 95\textsuperscript{th} stem listed in \cite{more-stable-stems}. Since this \namecref{section:new-diff} builds upon the results of \cite{more-stable-stems}, we assume the reader is already familiar with \cite{more-stable-stems}.

\section{Differentials and hidden extensions}\label{section:diff-hidden}
The arguments of this \namecref{section:diff-hidden} are quite generally applicable, and we shall present them in more generality than are needed for our calculations. In this \namecref{section:diff-hidden}, we work in $\Syn_E$ for some fixed Adams type spectrum $E$, and $X$ and $Y$ will be arbitrary synthetic spectra, not necessarily of the form $\nu (-)$. To streamline the presentation, we shall adopt the following conventions:
\begin{itemize}
  \item ``The Adams spectral sequence of $X$'' will mean ``the $\tau$-Bockstein spectral sequence of $X$ with the change of sign'' (recall that for any spectrum $X$, the $\tau$-Bockstein spectral sequence of $\nu X$ agrees with the Adams spectral sequence of $X$ up to a sign \cite[Theorem A.1]{manifold-synthetic}).
  \item We will write $X/\tau^k$ for $C\tau^k \otimes X$.
  \item We will omit all suspensions $\Sigma^{a, b}$ in $\Syn_E$; they can be inferred from context if necessary.
\end{itemize}

\begin{notation}
  Define maps $r_m, r_{n, m}, \delta_m, \delta_{n, m}$ by the cofiber sequences
  \[
    \begin{tikzcd}[row sep=tiny]
      X \ar[r, "\tau^n"] & X \ar[r, "r_m"] & X / \tau^m \ar[r, "\delta_m"] & X \\
      X/\tau^{n - m} \ar[r, "\tau^m"] & X/\tau^n \ar[r, "r_{n,m}"] & X / \tau^m \ar[r, "\delta_{n, m}"] & X/\tau^{n - m}.
    \end{tikzcd}
  \]
  Note that $\tau^m$ will always denote a map $X/\tau^{n - m} \to X/\tau^n$, as opposed to the endomorphism of $X/\tau^n$ of the same name. In particular, $\tau^m$ is non-zero on $X/\tau$.
\end{notation}

\begin{notation}
  If $x \in \pi_{*, *}X/\tau^m$ and $y \in \pi_{*, *}X/\tau^n$ are such that $r_{m, k} x = r_{n, k} y$, we say $x \equiv y \mod \tau^k$. Note in particular that $x$ and $y$ may live in different groups.
\end{notation}

One immediately sees that
\begin{lemma}\pushQED{\qed}
  For any $n, k > m$, we have a commutative diagram
  \[
    \begin{tikzcd}
      X \ar[r, "\tau^n"] \ar[d, "\tau^{n - m}"] & X \ar[r, "r_n"] \ar[d, equals] & X/\tau^n \ar[r, "\delta_n"] \ar[d, "r_{n, m}"] & X \ar[d, "\tau^{n - m}"] \\
      X \ar[r, "\tau^m"] \ar[d, equals] & X \ar[r, "r_m"] \ar[d, "\tau^{k - m}"] & X/\tau^m \ar[r, "\delta_m"] \ar[d, "\tau^{k - m}"] & X \ar[d, equals] \\
      X \ar[r, "\tau^k"] \ar[d, "r_{\ell - k}"] & X \ar[r, "r_k"] \ar[d, "r_\ell"] & X/\tau^k \ar[r, "\delta_k"] \ar[d, equals] & X  \ar[d, "r_{\ell - k}"]\\
      X/ \tau^{\ell - k} \ar[r, "\tau^k"] & X / \tau^\ell \ar[r, "r_{\ell, k}"] & X/\tau^k \ar[r, "\delta_{\ell, k}"] & X/ \tau^{\ell - k}.
    \end{tikzcd}\qedhere
  \]
\end{lemma}

The following are standard properties of Bockstein spectral sequences, whose proofs are left to the reader.
\begin{lemma}\pushQED{\qed}
  Let $x \in \pi_{*, *} X/\tau$.
  \begin{enumerate}
    \item For any representative of $d_{k + 1}(x)$ on the $E_2$ page, there is a lift of $x$ to $[x] \in \pi_{*, *} X/\tau^k$ such that $\delta_k [x] \equiv -d_{k + 1}(x) \mod \tau$.
    \item If $\tau^k x = 0$, then $x$ is the target of a $d_{k + 1}$ differential.
    \item If $\delta x = \tau^{k - 2} y$ for some $y \in \pi_{*, *} X$, then $x$ survives to the $E_k$ page, and $y \equiv -d_k(x) \mod \tau$.\qedhere
  \end{enumerate}
\end{lemma}

We now define hidden extensions. Classically, they are defined for classes on the $E_\infty$ page in terms of multiplication in homotopy groups. For our purposes, we need to generalize this to potentially non-surviving classes. Such a notion was first introduced by Cooley in his thesis \cite[pp.\ 18--21]{cooley-thesis}, together with a version of \Cref{thm:hidden-ext} \cite[Theorem 1.24]{cooley-thesis}. While we believe our definition agrees with Cooley's, we shall make no attempts to compare them.

Fix a map of synthetic spectra $\alpha\colon X \to Y$.
\begin{defi}\label{defi:hidden-ext}
  Let $x\in \pi_{*, *} X/\tau$ and $y \in \pi_{*, *} Y / \tau$. Suppose $x$ survives to the $E_r$ page and $s < r - 1$. We say there is a hidden $\alpha$-extension by $s$ from $x$ to $y$ on the $E_r$ page if there is a lift $\{x\}$ of $x$ to $\pi_{*, *} X / \tau^{r - 1}$ and $\{y\}$ of $y$ to $\pi_{*, *} Y / \tau^{r - 1 - s}$ such that
  \[
    \alpha \{x\} = \tau^s \{y\}.
  \]
  Alternatively, this says $\alpha \{x\}$ is $\tau^s$-divisible, and a $\tau^s$ division of $\alpha\{x\}$ is equal to $y$ mod $\tau$.

  We say this hidden extension is maximal if $\alpha\{x\}$ is not $\tau^{s + 1}$ divisible. This is automatic if $r = s + 2$. In case $r = \infty$ and $\alpha\{x\}$ is $\tau^s$ divisible for all $s$ (e.g.\ it is zero), we say there is a maximal hidden extension by $\infty$ to $0$.
\end{defi}
In particular, a hidden extension by $0$ is a regular, non-hidden extension.

\begin{remark}
  The jump $s$ is redundant information given $x$, $y$ and $\alpha$, and we omit it when no confusion can arise.
\end{remark}

\begin{remark}
  After fixing an $\{x\}$, the value of $y$ is well-defined up to images of $d_2, \ldots, d_{s + 1}$, and we shall consider $y$ as an element in this quotient. It is, however, inaccurate to say it is well-defined on the $E_{s + 2}$ page; it may not survive that long.

  Of course, different lifts $\{x\}$ give different values of $y$, and in general they can belong to different filtrations. However, this is not an issue when $s = 1$; there is a hidden extension by $1$ iff $\alpha x = 0$ on the $E_2$ page, and the indeterminacy in $y$ is exactly $\alpha$-multiples of classes in the bidegree right above $x$ on the $E_2$ page.
\end{remark}

\begin{remark}
  Let $s + 1 < q < r$. If there is a hidden $\alpha$-extension by $s$ from $x$ to $y$ on the $E_r$ page, then there is a hidden $\alpha$-extension from $x$ to $y$ on the $E_q$ page. The converse holds if there is no indeterminacy (and $x$ survives long enough).
\end{remark}

\begin{thm}[Generalized Leibniz rule]\label{thm:extend-diff}
  Let $x \in \pi_{*, *} X/\tau$ survive to the $E_r$ page. Fix a representative of $d_r(x)$ on the $E_2$ page. Then there is a differential from a maximal $\alpha$-extension of $x$ on the $E_r$ page to a maximal $\alpha$-extension of $d_r(x)$ on the $E_\infty$ page.
\end{thm}

\begin{proof}
  Pick a lift $\{x\}$ of $x$ to $\pi_{*, *}X/\tau^{r - 1}$ such that $\delta_{r - 1}\{x\}$ is a lift of $-d_r(x)$. Then we have
  \[
    \alpha \{x\} = \tau^{?} y,\quad \alpha \delta_{r - 1}\{x\} = \tau^\iq z
  \]
  for some $y$ and $z$ whose reduction mod $\tau$ are maximal hidden $\alpha$-extensions of $x$ and $-d_r(x)$ respectively. Then
  \[
    \delta_{r - ? - 1} y = \delta_{r - 1} \tau^? y = \delta_{r - 1} \alpha \{x\} = \alpha \delta_{r - 1} \{x\} = \tau^\iq z.
  \]
  So
  \[
    \delta (r_{r - ? - 1, 1} y) = \tau^{r - ? - 2} \delta_{r - ? - 1} y = \tau^{r - ? - \iq - 2} z.\qedhere
  \]
\end{proof}

\begin{remark}
  There are also differentials between non-maximal extensions, but they all vanish since they are pre-empted by shorter differentials.
\end{remark}

We end the \namecref{section:diff-hidden} with a result that identifies hidden extensions with differentials in the cofiber, which can be useful if we want to compute longer hidden extensions by hand. Define $C\alpha, \iota_\alpha, \delta_\alpha$ by the cofiber sequence
\[
  \begin{tikzcd}
    X \ar[r, "\alpha"] & Y \ar[r, "\iota_{\alpha}"] & C\alpha \ar[r, "\delta_{\alpha}"] & X.
  \end{tikzcd}
\]

\begin{thm}\label{thm:hidden-ext}
  Let $x \in \pi_{*,*} X/\tau$ be such that $d_{k + 1}(x) = 0$. Suppose $\bar{x} \in \pi_{*, *} C\alpha / \tau$ is such that $\delta_\alpha \bar{x} = x$, and suppose $y \in \pi_{*, *} Y/\tau$ is such that $\iota_\alpha y = d_{k + 1} \bar{x}$ on the $E_k$ page. Then there is a hidden $\alpha$ extension from $x$ to $y$ on the $E_{k + 2}$ page.
\end{thm}

\begin{proof}
  Consider the cofiber sequences
  \[
    \begin{tikzcd}[row sep = tiny]
      X \ar[r, "\alpha"] & Y \ar[r, "\iota_{\alpha}"] & C\alpha \ar[r, "\delta_{\alpha}"] & X \\
      \S/\tau \ar[r, "\tau^k"] & \S/\tau^{k + 1} \ar[r, "r_{k + 1, k}"] & \S/\tau^k \ar[r, "\delta_{k + 1, k}"] & \S/\tau
    \end{tikzcd}
  \]
  Taking the tensor product of these cofiber sequences gives
  \[
    \begin{tikzcd}
      Y/\tau^{k + 1} \ar[d, "r_{k + 1, k}"] \ar[r, "\iota_\alpha"] & C\alpha/\tau^{k + 1} \ar[d, "r_{k + 1, k}"] \ar[r, "\delta_\alpha"] & X/\tau^{k + 1} \ar[d, "r_{k + 1, k}"] \\
      Y/\tau^k \ar[r, "\iota_\alpha"] \ar[d, "\delta_{k + 1, k}"] & C\alpha/\tau^k \ar[d, "\delta_{k + 1, k}"] \ar[r, "\delta_\alpha"] & X/\tau^k \ar[d, "\delta_{k + 1, k}"] \\
      Y/\tau \ar[r, "\iota_\alpha"] & C\alpha/\tau \ar[r, "\delta_\alpha"] & X/\tau.
    \end{tikzcd}
  \]
  Since $d_{k + 1} \bar{x} = \iota_\alpha y$ on the $E_k$ page, we can pick a lift $\{\bar{x}\} \in \pi_{*, *} C \alpha$ of $\bar{x}$ such that
  \[
    \delta_{k + 1, k} \{\bar{x}\} = -\iota_\alpha y.
  \]
  By \cite[Section 6]{may-additivity} (see also \cite[Lemma 9.3.2]{inverting-hopf}), there is an $\{x\} \in X/\tau^{k + 1}$ such that
  \[
    r_{k + 1, k} \{x\} = \delta_\alpha \{\bar{x}\},\quad \alpha \{x\} = \tau^k y.
  \]
  The first condition tells us
  \[
    \{x\} \equiv \delta_\alpha \{\bar{x}\} \equiv \delta_\alpha \bar{x} = x \mod \tau.
  \]
  So $\{x\}$ is a lift of $x$ to $X/\tau^{k + 1}$, and the result follows.
\end{proof}

\section{Computation of old differentials}\label{section:old-diff}
\DeclareSseqCommand\htwo {} {
  \class (\lastx + 3, \lasty + 1) \structline
}
\DeclareSseqCommand\hone {} {
  \class (\lastx + 1, \lasty + 1) \structline
}
\DeclareSseqCommand\hzero {} {
  \class (\lastx, \lasty + 1) \structline
}
\DeclareSseqCommand\hzeroi {} {
  \class (\lastx, \lasty - 1) \structline
}

\DeclareSseqGroup\tower {m} {
  \class(0, 0)
  \savestack
  \foreach \n in {1, ..., #1} { \hzero }
  \restorestack
}

\DeclareSseqGroup\zigzag {O{}} {
  \class(0, 2)
  \hzeroi \hzeroi
  \savestack
  \hone
  \hzeroi \hzeroi
  \restorestack
}
\begin{sseqdata}[
  name = s2 ass,
  Adams grading,
  grid = go, tick step = 4, grid step = 2,
  classes = { fill, minimum size = 0.25em },
  class labels = { below },
  class label handler={\def\result{\scalebox{0.5}{$#1$}}},
  label distance = -1pt,
  differentials = { blue, - },
  x range = {0}{35}, y range = {0}{20}, scale = 0.462]
  \tower{21}
  \savestack
  \hone \classoptions["h_1"] \hone \hone \structline(0, 2)
  \restorestack
  \htwo \classoptions["h_2"]
  \savestack
  \hzero \structline(3, 3) \structline(0, 1)
  \restorestack
  \htwo \htwo

  \tower(7, 1){3}
  \classoptions["h_3"]
  \hone \structline(9, 3)

  \class["c_0" left](8, 3) \hone

  \class["Ph_1"](9, 5)\hone \hone \hzeroi \hzeroi \classoptions["Ph_2"]
  \foreach \n in {2, 3, 4} {
    \class["P^{\n} h_1"](1 + 8 * \n, 1 + 4 * \n)\hone \hone \hzeroi \hzeroi \classoptions["P^{\n} h_2"]
  }

  \class["Pc_0" left](16, 7) \hone
  \class(24, 11) \hone
  \class(32, 15) \hone

  \class["h_3^2"](14, 2) \hzero
  \class["d_0"] (14, 4) \hzero \hzero \structline(11, 5)

  \tower(15, 1){7} \classoptions["h_4"]
  \hone \hone \hone \structline(15, 3)
  \class(18, 2) \structline(15, 1) \hzero \structline(15, 2) \structline(18, 4, 1)
  \class(21, 3) \structline(18, 2)

  \class(15, 5) \structline(14, 4)
  \hone
  \hone \structline(14, 6)
  \hzeroi \structline(14, 5)
  \hzeroi \structline(14, 4)
  \hzeroi \classoptions["e_0" left]

  \hone \hzeroi \classoptions["f_0" right]

  \class(20, 6) \structline(17, 5)
  \hzeroi \structline(17, 4)
  \hzeroi \classoptions["g"]
  \hone \structline(18, 4, 2)
  \class(23, 5) \structline(20, 4) \htwo
  \class(23, 6) \structline(23, 5) \structline(20, 5)

  \class["c_1"](19, 3)\htwo
  \class["h_4 c_0"](23, 4) \hone

  \zigzag(25, 8)

  \tower(23, 9){3} \hone \hone \structline(25, 10)
  \zigzag(22, 8)
  \classoptions["Pd_0"]
  \zigzag(28, 8)
  \classoptions["d_0^2" left]

  \structline(23, 8)(23, 9, 1)
  \structline(23, 9, 2)(24, 10)

  \classoptions["k"](29, 7)
  \tower(31, 1){14}
  \classoptions["h_5"]
  \hone \hone \hone \structline(31, 3)
  \hzeroi \structline(31, 2)
  \hzeroi \structline(31, 1)
  \htwo

  \tower(30, 2){3} \classoptions["h_4^2"]
  \hone

  \tower(30, 6){5}
  \classoptions["\Delta h_2^2"]
  \htwo
  \class["\Delta h_1 h_3"](32, 6) \structline(33, 7)

  \tower(30, 12){2}
  \classoptions["P^2 d_0" left]
  \hone \hone \hone
  \zigzag(33, 12)
  \structline(33, 14)(33, 15)

  \zigzag(31, 8)\structline(31, 9, 1)

  \zigzag(34, 8)

  \class["n"](31, 5) \htwo \htwo

  \class["d_1"](32, 4) \htwo \htwo
  \class["p"](33, 4) \hzero \structline(32, 4)

  % these only exist for the h_2 structlines to be drawn
  \zigzag(37, 8)
  \zigzag(36, 12)

  \foreach \x/\y in {22/8, 22/9, 22/10, 23/7, 23/8, 25/8, 25/9, 26/7, 26/8, 29/7, 29/8, 32/7, 32/8, 34/8, 34/9, 35/7, 35/8, 30/12, 30/13, 30/14, 33/12, 33/13, 34/11, 34/12} {
    \structline(\x, \y)(\x + 3, \y + 1)
  }
  \structline(28, 8)(31, 9, 2)
  \structline(28, 9)(31, 10, 2)
  \structline(31, 8, 2)(34, 9)
  \structline(31, 9, 2)(34, 10)

  \structline(19,9)(22,10)
  \structline(27,13)(30,14)
  \class(38, 18)
  \structline(35,17)(38,18)

  % differentials
  \d2(15, 1)
  \d[red]3(15, 2)
  \d[red]3(15, 3)

  \d2(17, 4)
  \d2(18, 4, 2)
  \d2(18, 5)

  \d2(23, 7)
  \d2(23, 8)
  \d2(25, 8)
  \d2(26, 7)
  \d2(26, 8)
  \d2(26, 9)
  \d2(29, 7)
  \d2(29, 8)

  \d2(31, 1)
  \d2(31, 2)
  \d2(31, 3)
  \d[red]3(31, 4)
  \d[red]3(31, 5)
  \d[red]3(31, 6)
  \d[red]3(31, 7)
  \d[red]3(31, 8, 1)
  \d[green!75!black]4(31, 8, 2)
  \d[green!75!black]4(31, 9, 1)
  \d[green!75!black]4(31, 10, 1)
  \d[green!75!black]4(32, 9)(31, 13, 2)

  \d[red]3(30, 6)

  \d2(32, 7)(31, 9, 2)
  \d2(32, 8)(31, 10, 2)

  \d2(35, 7)
  \d2(35, 8)
  \d2(33, 12)
  \d2(34, 11)
  \d2(34, 12)
  \d2(34, 13)
  \d[red]3(34, 2)

  \d2(38, 7)

  \class(38, 4)
  \d[red]3(38, 4)

  % hidden extensions
  \structline[red, dashed, page = \infty](15, 4)(16, 7)
  \structline[red, bend left = 20, page = \infty](20, 6)(23, 9, 2)
  \structline[red, page = \infty](23, 6)(23, 9, 2)
  \structline[red, page = \infty](21, 5)(22, 8)

  \structline[blue, page = \infty](23, 9, 1)(24, 11)
  \structline[blue, page = \infty](30, 2)(33, 4)
  \structline[blue, dashed, bend left = 20, page = \infty](32, 6)(35, 9)
  \structline[green!75!black, dashed, page = \infty](31, 11)(32, 15)
\end{sseqdata}

\begin{figure}[p]
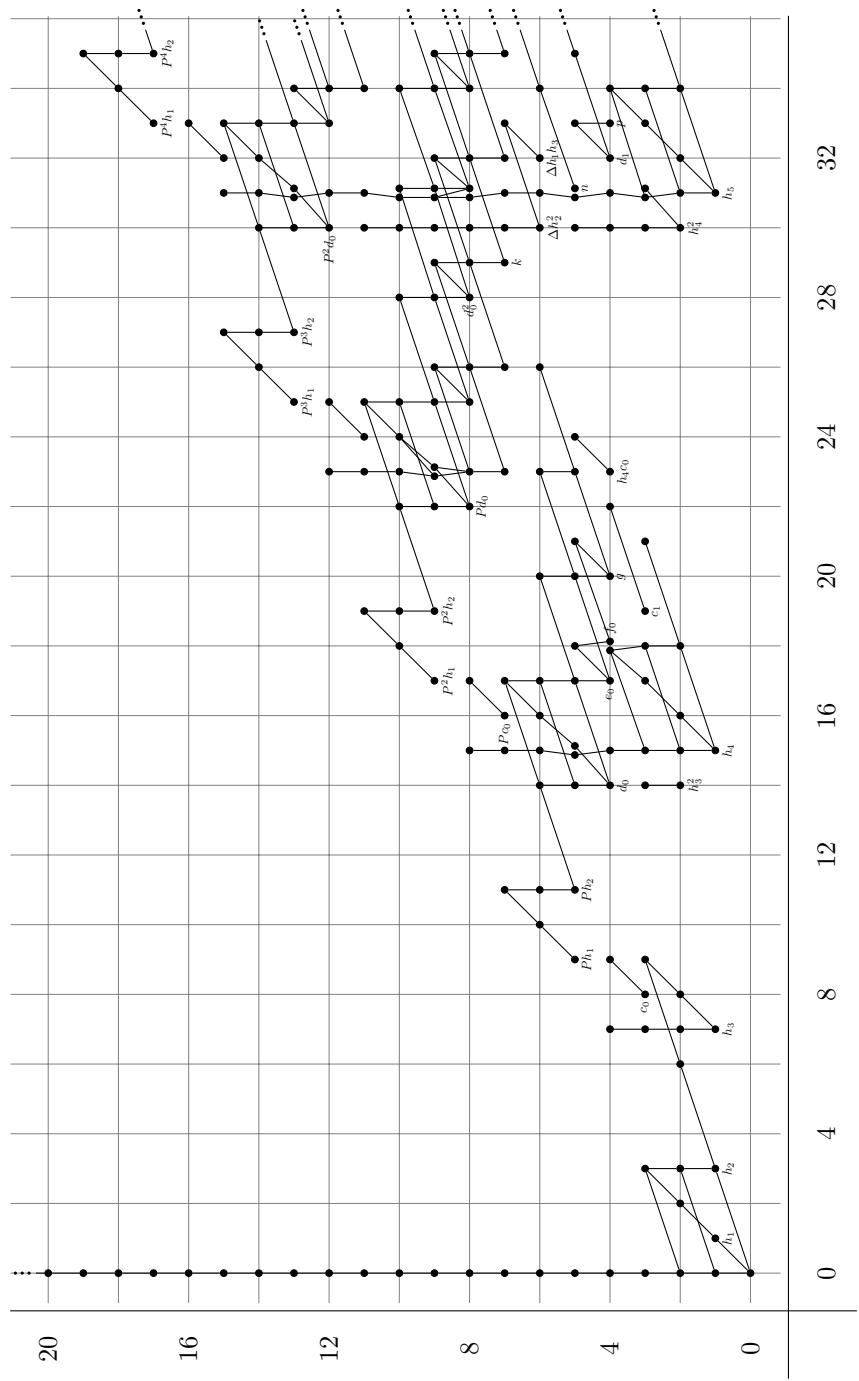

  \centering
  \begin{sideways}
    \printpage[name = s2 ass, page = 1]
  \end{sideways}
  \caption{The $E_2$ page of the classical Adams spectral sequence}\label{fig:e2}
\end{figure}

\begin{figure}[p]
  \centering
  \begin{sideways}
    \printpage[name = s2 ass, page = 2--\infty]
  \end{sideways}
  \caption{The differentials of the classical Adams spectral sequence}\label{fig:diff}
\end{figure}

\begin{figure}[p]
  \centering
  \begin{sideways}
    \printpage[name = s2 ass, page = \infty]
  \end{sideways}
  \caption{The $E_\infty$ page of the classical Adams spectral sequence}\label{fig:einfty}
\end{figure}

To illustrate how one can make use of the computer-generated data, we compute all differentials in the Adams spectral sequence up to the 35\textsuperscript{th} stem at the prime $2$ and resolve most hidden extension. The resulting Adams charts are displayed in \Cref{fig:e2,fig:diff,fig:einfty} (the dashed hidden extensions in the $E_\infty$ page are those we do not prove). Many of the arguments can be simplified if we are willing to use other tools, but we restrict ourselves to ``straightforward'' manipulations using the computer data.
\afterpage{\clearpage}
\subsection*{Conventions}
We assume the reader is familiar with the names of classes in the homotopy groups of spheres and the classical Adam $E_2$ page, as well as the translation between the two. For convenience, we label the relevant $E_2$ page names in the Adams charts as well.

We adopt the following naming conventions:
\begin{enumerate}
  \item If $\alpha \in \pi_* \S$ is an element in the classical homotopy groups of the sphere, we use the same name to denote the corresponding element in the homotopy groups of the \emph{synthetic} sphere. By this we mean an element in $\pi_{*, *} \S$ whose $\tau$ inversion gives the original class $\alpha$, and has maximum Adams filtration amongst such elements. While this is potentially ambiguous, the ambiguity is irrelevant in all cases of interest in this \namecref{section:old-diff}.

    To avoid any confusion, we shall never refer to the classical homotopy groups of the sphere in this \namecref{section:old-diff}. All such names always refer to the synthetic version.

  \item If $a \in \Ext_\A(\F_2, \F_2)$ is a permanent cycle, we use $\{a\}$ to denote any lift of $a$ to $\pi_{*, *} \S$. Again the ambiguities end up being irrelevant.
  \item If $a \in \Ext_\A(\F_2, \F_2)$ survives the $d_2$ differential, we use $[a]$ to denote the specific lift to $\pi_{*, *} C\tau^2$ constructed in \Cref{thm:lift} using the minimal resolution generated by our program. Of course, the precise choice of lift is irrelevant; what matters is that $[a]$ refers to the same lift throughout the whole dataset.

    Note that in general, $[a + b] \not= [a] + [b]$. Instead, there is a correction term as specified in \Cref{thm:lift}.
\end{enumerate}

\subsection{Differentials in stems \texorpdfstring{$0$}{0} to \texorpdfstring{$28$}{28}}
We first look at the differentials in the first $28$ stems, which are relatively straightforward.

\begin{lemma}
  We have
  \[
    d_3(h_0 h_4) = h_0 d_0.
  \]
\end{lemma}

\begin{proof}
  In the computer data, we see that
  \[
    [d_0][h_0 h_4] = \tau k.
  \]
  Since $[d_0]$ detects $\kappa$, this means there is a hidden $\kappa$-extension from $h_0 h_4$ to $k$. So this follows from dividing the differential $d_2(k) = h_0 d_0^2$.
\end{proof}

\begin{cor}
  We have $\delta h_4 = \tilde{2} \sigma^2$.
\end{cor}

\begin{proof}
  We know that $\pi_{14, 3} \S$ is spanned by $\tilde{2} \sigma^2$ and $\tau \kappa$ as an $\F_2$-module. Since $\delta h_4 = \tilde{2} \sigma^2 \mod \tau$, we know that $\delta h_4 = \tilde{2} \sigma^2 + a \tau \kappa$ for some $a \in \F_2$. By computer calculation, $\tilde{2}^2 \sigma^2 = \tau \tilde{2} \kappa$. So we get
  \[
    \delta h_0 h_4 = \tilde{2}^2 \sigma^2 + a \tau \tilde{2} \kappa = (a + 1) \tau \tilde{2} \kappa,
  \]
  Since $d_3 (h_0 h_4) = h_0 d_0$, this expression must also equal $\tau \tilde{2} \kappa$. So we must have $a = 0$, as desired.
\end{proof}

\begin{remark}
  One can prove these two results in the opposite order. First, using the fact that $\tau \tilde{2} \sigma^2 = 2 \sigma^2 = 0$, we learn that we must have $\delta h_4 = \tilde{2} \sigma^2$. Then $\delta h_0 h_4 = \tilde{2}^2 \sigma^2 = \tau \tilde{2}\kappa$. So $d_3(h_0 h_4) = h_0 d_0$.
\end{remark}

\begin{cor}
  $h_1 h_4$, $h_2 h_4$ and $c_0 h_4$ are permanent.
\end{cor}

\begin{proof}
  We have $\delta (h_1h_4) = \eta \delta h_4 = \eta \tilde{2} \sigma^2 = 0$ since $\tilde{2} \eta = 0$. The others follow similarly with $\nu \sigma = 0$ and $\tilde{2} \epsilon = 0$.
\end{proof}

\begin{remark}
  This is the same proof as the Moss' convergence theorem proof, constructing $[h_1 h_4]$ as $\langle \sigma^2, 2, \eta \rangle$. When applying Moss' convergence theorem, one has to verify that the product vanishes in homotopy and work out indeterminacies. In the synthetic proof, this translates to keeping track of higher $\tau$-divisible terms that can show up in the products and $\delta$. In particular, knowing the full value of $\delta$ instead of just the Adams differential is often extremely useful for future computations.
\end{remark}

\begin{lemma}
  $g$ is permanent.
\end{lemma}

\begin{proof}
  If $g$ supported a differential, then so would $Pg = d_0^2$, since $P$ acts injectively on all potential targets with no indeterminacy. But $d_0^2$ is permanent.
\end{proof}

\subsection{Hidden extensions in stems \texorpdfstring{$0$}{0} to \texorpdfstring{$28$}{28}}

\begin{lemma}
  We have
  \[
    \nu^3 + \eta^2 \sigma = \tau \eta \epsilon.
  \]
\end{lemma}

\begin{proof}
  The computer data gives
  \[
    [h_2]^3 = [h_1^2 h_3],\quad [h_1]^2 [h_3] = [h_1^2 h_3] + \tau h_1 c_0.\qedhere
  \]
\end{proof}

\begin{lemma}
  $\delta e_0 = \eta^2 \kappa$. Thus, $\tau \eta^2 \kappa = 0$.
\end{lemma}

\begin{proof}
  Since $d_2(e_0) = h_1^2 d_0$, the only other possibility is $\delta e_0 = \eta^2 \kappa + \tau \{Pc_0\}$. This would imply that
  \[
    \tau \eta^2 \kappa = \tau^2 \{Pc_0\}.
  \]
  Multiplying by $\eta$ gives
  \[
    \tau^2 \eta \{Pc_0\} = \tau \eta^3 \kappa = \tau \tilde{2}^2 \nu \kappa = 0,
  \]
  which is a contradiction.
\end{proof}

\begin{remark}
  One can similarly show that $\delta f_0 = \tilde{2} \nu \kappa$.
\end{remark}

\begin{lemma}
  We have
  \[
    \nu^3 \kappa = \tau^2 \eta \{Pd_0\}.
  \]
\end{lemma}
Note that this hidden extension jumps by $2$ filtrations, and we are able to compute this by iterating hidden extensions by $1$.

\begin{proof}
  Since $\nu^3 = \eta^2 \sigma + \tau \eta\epsilon$, multiplying by $\kappa$ gives
  \[
    \nu^3 \kappa = \eta^2 \sigma \kappa + \tau \eta \epsilon \kappa.
  \]
  By computer calculation, we know that
  \[
    \epsilon \kappa = \tau \{Pd_0\}.
  \]
  So it remains to show that the first term vanishes. Since the 23\textsuperscript{rd} stem is non-$\tau$-torsion, it suffices to show that $\tau \eta^2 \sigma \kappa = 0$. But we have already seen that $\tau \eta^2 \kappa = \tau \delta e_0 = 0$. So we are done.
\end{proof}

\begin{cor}
  We have
  \[
    \tilde{2}^2 \nu \kappabar = \tau^2 \eta \{Pd_0\},\quad \eta \kappabar = \tau^2 \{Pd_0\}.
  \]
\end{cor}

\begin{proof}
  The first follows from the identity
  \[
    \tilde{2}^2 \kappabar = \nu^2 \kappa.
  \]
  The second follows from $\eta^3 = \tilde{2}^2 \nu$.
\end{proof}

\begin{lemma}
  We have
  \[
    \sigma \{Ph_1\} = \eta^2 \kappa + \tau \{Pc_0\}.
  \]
  Thus, in the (classical) stable homotopy groups of spheres, there is a hidden $\sigma$ extension from $Ph_1$ to $Pc_0$.
\end{lemma}

\begin{proof}
  This follows from the computer data, since there are no higher filtration terms.
\end{proof}
\subsection{Stems \texorpdfstring{$29$}{29} to \texorpdfstring{$35$}{35}}
\begin{lemma}\label{lemma:h03-h5}
  $d_3(h_0^3 h_5) = h_0 \Delta h_2^2$ and $d_4(h_0^8 h_5) = h_0 P^2 d_0$.
\end{lemma}

\begin{proof}
  To compute $d_3(h_0^3 h_5)$, the obvious approach of starting with $d_2(h_0^2 h_5) = h_0^3 h_4^2$ and then computing a hidden $\tilde{2}$ extension does not work, since $h_0 \Delta h_2^2$ is in the indeterminacy. Instead, we start with $d_2(h_0 h_5) = h_0^2 h_4^2$ and compute a hidden $\tilde{2}^2$ extension. Indeed, computer calculation gives
  \[
    \begin{aligned}
      \relax[h_0] [h_0^2 h_4^2] &= [h_0^3 h_4^2]\\
      \relax[h_0] [h_0^3 h_4^2] &= \tau [h_0 \Delta h_2^2].
    \end{aligned}
  \]
  So there is a hidden $\tilde{2}^2$ extension from $h_0^2 h_4^2$ to $h_0 \Delta h_2^2$ with no indeterminacy, and the $d_3$ follows. The $d_4$ is similar.
\end{proof}

\begin{lemma}
  $d_3(d_0 e_0) = h_0^5 \Delta h_2^2$ and $d_4(d_0e_0 + h_0^7 h_5) = P^2 d_0$.
\end{lemma}
% d_0 e_0 is [1, 1] in our charts
\begin{proof}
  We have
  \[
    \delta (d_0 e_0) = \kappa \delta e_0 = \eta^2 \kappa^2,
  \]
  which computer calculation tells us is $\tau h_0^5 \Delta h_2^2$ mod $\tau^2$. The next differential follows from $h_0$-division in a purely classical manner, since $h_0 d_0 e_0 = 0$ on the $E_3$ page.
\end{proof}

\begin{remark}
  Here it is important for us to precisely identify the value of $\delta e_0$. A simple hidden extension argument would not work since $h_0^5 \Delta h_2^2 = d_0 Pc_0$ is in the indeterminacy.
\end{remark}

\begin{cor}
  $d_3(\Delta h_2^2) = h_1 d_0^2$.
\end{cor}

\begin{proof}
  The source and target have hidden $\eta$ extensions to $d_0e_0$ and $h_0^5 \Delta h_2^2$ respectively.
\end{proof}

\begin{cor}
  We have
  \[
    \delta h_5 = \tilde{2} \{h_4^2\}.
  \]
\end{cor}

\begin{proof}
  From our calculations, $\pi_{30, 3} \S = \F_2$ and is generated by $\tilde{2} \{h_4^2\}$.
\end{proof}

\begin{cor}
  $h_1 h_5$ and $p$ are permanent and $d_3(h_2 h_5) = h_0 p$.
\end{cor}

\begin{proof}
  The first follows from $\tilde{2} \eta = 0$. The rest follow from the hidden $\nu$ extension from $h_4^2$ to $p$.
\end{proof}

To prove that the remaining elements are permanent, we have to look beyond what our charts in \Cref{fig:e2,fig:diff,fig:einfty} cover. The reader can instead refer to Isaksen's charts at \cite{charts}.
\begin{lemma}
  $d_1$ and $\Delta h_1 h_3$ are permanent.
\end{lemma}

\begin{proof}
  These elements can only hit $h_0^{15} h_5$, and it is easy to check that $Ph_0^{15} h_5$ cannot be hit by an element of filtration at least $8$.
\end{proof}

\begin{lemma}
  $h_0 h_2 h_5$ is permanent.
\end{lemma}

\begin{proof}
  Since $[h_0] [h_0 p] = 0$, we know $h_0 h_2 h_5$ does not hit $h_1 \Delta h_1 h_3$. So the only potential targets are $h_1 P^3 c_0$ and $P^4 h_1$. We again rule this out by Adams periodicity. We have
  \[
    P h_0 h_2 h_5 = h_0 P h_2 h_5.
  \]
  Since $P [h_0 p] = 0$ with no indeterminacy in $\Mod_{C\tau^2}$, we know that $d_4(P h_2 h_5) = 0$. So the shortest differential $P h_2 h_5$ can support hits $h_1 P^4 c_0$ or $P^5 h_1$, both of which preclude a differential on $h_0 P h_2 h_5$.
\end{proof}

\section{Computation of new differentials}\label{section:new-diff}
We now turn to the computation of new differentials. These differentials are listed in \Cref{table:new-diff}, with the proofs indicated in the last column. Many of these new differentials are easy consequences of the generalized Leibniz rule using hidden $\tilde{2}$, $\eta$, $\nu$ and $\sigma$ extensions, which are listed in \Cref{table:hidden-two,table:hidden-eta,table:hidden-nu,table:hidden-sigma}. We shall not provide further explanation for these differentials. The remainder of the differentials require extra arguments, and are explained in \Cref{section:hard-diff}.

Throughout the \namecref{section:new-diff}, we shall use the names of \cite{more-stable-stems}. The identifications of their classes in our basis are listed in \Cref{table:class-ident} with brief justifications. We encourage the reader to refer to the charts at \cite{charts} when reading this \namecref{section:new-diff}.

\begin{remark}
  Some of these new differentials have been independently computed in unpublished work of Burklund--Isaksen--Xu. Specifically, they have computed the differentials on $\Delta^2 g_2$, $h_1 \Delta^2 g_2$, $h_0^3 \Delta h_2^2 h_6$, $x_{95, 7}$, $\Delta^2 Mh_1$ and $\Delta^2 M h_1^2$. Their arguments are similar to ours, except they had to compute hidden extensions by hand.
\end{remark}
\subsection{Computation of new differentials}\label{section:hard-diff}
\begin{lemma}\label{lemma:h07-h6}
  $d_3(h_0^7 h_6) = \Delta^2 h_0 h_3^2$.
\end{lemma}

\begin{proof}
  As in \Cref{lemma:h03-h5}, we start with $d_2(h_0^5h_6) = h_0^6 h_5^2$, and observe that there is a hidden $\tilde{2}^2$ extension from $h_0^6 h_5^2$ to $\Delta^2 h_0 h_3^2$. Indeed, we have
  \begin{align*}
    [h_0] [h_0^6 h_5^2] &= [h_0^7 h_5^2] + \tau (h_1 \Delta x + \Delta^2 h_3^2)\\
    [h_0] [h_0^7 h_5^2] &= 0.\qedhere
  \end{align*}
\end{proof}

\begin{lemma}\label{lemma:h0-g-b4}
  $d_3(h_0 g B_4) = h_0 \Delta^2 d_0 e_0 + \Delta h_1 e_0^2 g$.
\end{lemma}

\begin{proof}
  First, $h_0 g B_4$ is $h_2$-divisible with
  \[
    h_0 g B_4 = h_2 e_0 B_4,\quad d_2 (e_0 B_4) = h_0 M d_0 e_0.
  \]

  There is a hidden $\nu$-extension from $h_0 M d_0 e_0$ to $\Delta h_1 e_0^2 g$ with indeterminacy $h_0 \Delta^2 d_0 e_0$. So we know that
  \[
    d_3(h_0 g B_4) = \Delta h_1 e_0^2g + ?h_0 \Delta^2 d_0 e_0.
  \]

  To determine the indeterminacy, we consider further $\nu$-multiplication. There is a hidden $\nu$-extension from $h_0 g B_4$ to $h_0 \Delta^2 m$. Thus, we find that
  \[
    h_2 d_3(h_0 g B_4) = d_2(h_0 \Delta^2 m) = h_2 (h_0 \Delta^2 d_0 e_0).
  \]
  Since $h_2 \Delta h_1 e_0^2 g = 0$, the result follows.
\end{proof}

\begin{remark}\label{remark:err}
  Our calculations have uncovered two incorrect $d_3$'s in \cite{more-stable-stems}. Write $\tau_m$ for the $\tau$ in \cite{more-stable-stems}, which is $\tau^2$ in $\Syn_{BP}$.

  \cite[Lemma 5.21]{more-stable-stems} claims that $d_3(h_0 g B_4) = \Delta^2 h_0 d_0 e_0$ (note that $d_0 B_5 = h_0 g B_4$ in the classical Adams spectral sequence). Their argument neglects the possibility that in the motivic Adams spectral sequence, $d_3(\tau_\mot^2 d_0 B_5) = \Delta^2 h_0 d_0 e_0 + \tau_\mot^3 \Delta h_1 e_0^2 g$, which would make $\Delta^2 h_0 d_0 e_0$ $\tau_\mot$-divisible in the $E_\infty$-page. Indeed, our argument shows this is exactly what happens.

  \cite[Lemma 5.26]{more-stable-stems} claims that $d_3(M h_0 d_0 k) = P \Delta^2 h_0 d_0 e_0$. This is a clerical error; in $\mathrm{mmf}$, there is a $d_2$ hitting $\Delta h_1 d_0^2 e_0^2 + \tau_\mot^3 P \Delta h_1 dg^2$, so in the $E_3$ page we have $\Delta h_1 d_0^2 e_0^2 = \tau_\mot^3 P \Delta h_1 dg^2$. Since $\tau_\mot^2 h_0 M d_0 k$ has trivial image in $\mathrm{mmf}$, its $d_3$ must be $\Delta h_1 d_0^2 e_0^2 + \tau_\mot^3 P \Delta h_1 dg^2$.
\end{remark}

\begin{lemma}\label{lemma:delta3-h1-h3}
  $d_3(\Delta^3 h_1 h_3) = \Delta h_1 e_0^2 g$.
\end{lemma}

\begin{proof}
  We first show that $d_3(\Delta^3 h_1 h_3)$ is non-zero. This is the argument of \cite[Lemma 5.20]{more-stable-stems}. Since $Ph_1 \cdot \Delta^3 h_1 h_3$ supports a $d_4$, we know $\Delta^3 h_1 h_3$ supports a  differential of length at most $d_4$. The target bidegree of the $d_4$ is zero and computer calculation shows it does not support a $d_2$, so it must support a $d_3$.

  Next, we observe that there is a hidden $\nu$-extension by $1$ from $\Delta^3 h_1 h_3$ to $0$. So $d_3(\Delta^3 h_1 h_3)$ must be killed by $h_2$. This leaves $h_0 \Delta^2 d_0 e_0$ as the only possibility.
\end{proof}

\begin{lemma}\label{lemma:h1-d2-g2}
  $d_6(h_1 \Delta^2 g_2) = 0$.
\end{lemma}

\begin{proof}
  We have to show that $\delta (h_1 \Delta^2 g_2) = 0$ mod $\tau^5$. To do so, we use the $E_2$ page relation $h_1 \Delta^2 g_2 = h_3 \Delta^2 e_1$.

  Since $d_3(\Delta^2 e_1) = h_2^2 \Delta^2 n$, we can write
  \[
    \delta (\Delta^2 e_1) = \tau \nu^2 \{\Delta^2 n\} + {?} \tau^2 \eta \{M \Delta h_1 d_0\} + {?} \tau^3 \{\Delta h_1 g^3\}.
  \]
  We shall show that all terms are trivial mod $\tau^5$ after multiplication by $\sigma$.

  \begin{enumerate}
    \item Since $\nu \sigma = 0$, the first term vanishes completely.
    \item By computer calculation, we know
      \[
        [h_3] [M \Delta h_1 d_0] = \tau h_0^6 x_{91, 11}.
      \]
      Note that $h_0^6 x_{91, 11}$ is itself a $h_3$ multiple, hence is in the indeterminacy. The only term in the bidegree above $h_0^6 x_{91, 11}$ is $h_0^7 x_{91, 11}$. So we get
      \[
        \sigma \{M \Delta h_1 d_0\} = {?} \tau h_0^6 x_{91, 11} + {?} \tau^2 h_0^7 x_{91, 11} \mod {\tau^3},
      \]
      where the coefficients are potentially different. However, the left-hand side is permanent, while the terms on the right support differentials. So the coefficients must in fact vanish.

    \item Since hidden $\sigma$-extensions by $1$ vanish identically at bidegree $(85, 17)$, we know that
      \[
        \sigma \{\Delta h_1 g^3\} = 0 \mod \tau^2.\qedhere
      \]
  \end{enumerate}
\end{proof}

\begin{lemma}\label{lemma:x-94-8}
  $d_3(x_{94, 8}) = h_1 x_{92, 10}$.
\end{lemma}

\begin{proof}
  This follows from \cite[Remark 5.2]{more-stable-stems} since $d_2(x_{94, 8}) = 0$.
\end{proof}

\subsection{Tables}\label{section:tables}
This \namecref{section:tables} contains the following tables:
\begin{itemize}
  \item \Cref{table:new-diff} contains all the newly computed differentials and their proofs. The grey rows consist of old differentials we include for reference.
  \item \Cref{table:hidden-two,table:hidden-eta,table:hidden-nu,table:hidden-sigma} contain the hidden extensions that we use to compute the differentials. The first four columns are lifted straight out of computer-generated data, while the last two columns identify the names of the classes.
  \item \Cref{table:class-ident} gives the identification between the names in \cite{more-stable-stems} and our basis, together with a brief justification for each. We omit cases where the group is one-dimensional.
\end{itemize}

\begin{longtabu}{cccccc}
  \caption{Newly computed differentials}\label{table:new-diff} \\
  \toprule
  $n$ & $s$ & $r$ & source & target & proof \\
  \midrule
  \endfirsthead
  \caption{Newly computed differentials (continued)} \\
  \toprule
  $n$ & $s$ & $r$ & source & target & proof \\
  \midrule
  \endhead
  \bottomrule
  \endfoot
  63 &  8 & 3 & $h_0^7 h_6$ & $\Delta^2 h_0 h_3^2$ & \Cref{lemma:h07-h6} \\
  69 &  8 & 2 & $D_3'$ & $0$ & - \\
  69 &  8 & 3 & $D_3'$ & $h_2 M g$ & $\tilde{2}$ division \\
  \rowfont{\color{gray!70!black}} 69 & 10 & 2 & $P(A + A')$ & $h_0 h_2 M g$ & - \\
  80 & 14 & 3 & $h_0 g B_4$ & $h_0 \Delta^2 d_0 e_0 + \Delta h_1 e_0^2 g$ & \Cref{lemma:h0-g-b4} \\
  80 & 14 & 3 & $\Delta^3 h_1 h_3$ & $\Delta h_1 e_0^2 g$ & \Cref{lemma:delta3-h1-h3} \\
  \rowfont{\color{gray!70!black}} 85 & 17 & 2 & $M d_0 j$ & $h_0 MPd_0 e_0$ & - \\
  \rowfont{\color{gray!70!black}} 87 & 17 & 2 & $\Delta^3 h_1 d_0$ & $e_0^3 m$ & - \\
  88 & 18 & 3 & $\Delta^3 h_1^2 d_0$ & $\Delta h_1 d_0^2 e_0^2$ & $\eta$ multiplication\\
  88 & 18 & 3 & $h_2 M d_0 j$ & $\Delta h_1 d_0^2 e_0^2 + h_0 P \Delta^2 d_0 e_0$ & $\nu$ multiplication\\
  92 & 12 & 5 & $\Delta^2 g_2$ & $0$ & $\eta$ division \\
  93 &  9 & 5 & $h_0^2 \Delta h_2^2 h_6$ & $h_0^2 \Delta^2 g_2$ & \cite{more-stable-stems} \\
  93 & 10 & 6 & $h_0^3 \Delta h_2^2 h_6$ & $M\Delta h_2^2 e_0$ & $\tilde{2}$ multiplication \\
  93 & 13 & 6 & $h_1 \Delta^2 g_2$ & $0$ & \Cref{lemma:h1-d2-g2} \\
  94 &  8 & 2 & $x_{94, 8}$ & $0$ & - \\
  94 &  8 & 3 & $x_{94, 8}$ & $h_1 x_{92, 10}$ & \Cref{lemma:x-94-8} \\
  94 & 15 & 3 & $\Delta^2 M h_1$ & $M d_0 e_0^2$ & $\tilde{2}$ division \\
  \rowfont{\color{gray!70!black}} 94 & 17 & 3 & $M d_0 m$ & $MP\Delta h_1^2 d_0$ & $\eta$ division \\
  95 &  7 & 2 & $x_{95, 7}$ & $h_0 x_{94, 8}$ & - \\
  95 & 16 & 4 & $\Delta^2 M h_1^2$ & $M P \Delta h_0^2 e_0$ & $\eta$ multiplication \\
  \rowfont{\color{gray!70!black}} 95 & 19 & 2 & $x_{95, 19, 0}$ & $MP\Delta h_1^3 d_0$ & - \\
\end{longtabu}

\begin{longtable}{cccccc}
  \caption{Selected hidden $\tilde{2}$-extensions}\label{table:hidden-two} \\
  \toprule
  $n$ & $s$ & source & target & name of source & name of target \\
  \midrule
  \endfirsthead
  \caption{Selected hidden $\tilde{2}$-extensions (continued)} \\
  \toprule
  $n$ & $s$ & source & target & name of source & name of target \\
  \midrule
  \endhead
  \bottomrule
  \endfoot
  69 & 8 & $[1, 0]$ & $[1]$ & $D_3'$ & $P(A + A')$ \\
  92 & 14 & $[0, 1, 0]$ & $[1]$ & $h_0^2 \Delta^2 g_2$ & $M\Delta h_2^2 e_0$ \\
  93 & 18 & $[1]$ & $[1, 0]$ & $Md_0 e_0^2$ & $MP\Delta h_1^2 d_0$ \\
  94 & 15 & $[1]$ & $[1]$ & $\Delta^2 M h_1$ & $Md_0 m$ \\
\end{longtable}

\begin{longtable}{cccccc}
  \caption{Selected hidden $\eta$ extensions}\label{table:hidden-eta} \\
  \toprule
  $n$ & $s$ & source & target & name of source & name of target\\
  \midrule
  \endfirsthead
  \caption{Selected hidden $\eta$ extensions (continued)} \\
  \toprule
  $n$ & $s$ & source & target & name of source & name of target\\
  \midrule
  \endhead
  \bottomrule
  \endfoot
  86 & 19 & $[1]$ & $[0, 1, 1]$ & $e_0^3 m$ & $\Delta h_1 d_0^2 e_0^2$ \\
  91 & 17 & $[0, 0, 1]$ & $[0, 1]$ & $Md_0 \ell + h_0^6 x_{91, 11}$ & $MP\Delta h_1 d_0$ \\
  93 & 18 & $[1]$ & $[0, 1]$ & $M d_0 e_0^2$ & $MP\Delta h_0^2 e_0$\\
  94 & 17 & $[1]$ & $[1, 0]$ & $M d_0 m$ & $x_{95, 19, 0}$ \\
\end{longtable}

\begin{longtable}{cccccc}
  \caption{Selected hidden $\nu$ extensions}\label{table:hidden-nu} \\
  \toprule
  $n$ & $s$ & source & target & name of source & name of target \\
  \midrule
  \endfirsthead
  \caption{Selected hidden $\nu$ extensions (continued)} \\
  \toprule
  $n$ & $s$ & source & target & name of source & name of target \\
  \midrule
  \endhead
  \bottomrule
  \endfoot
  76 & 15 & $[1]$ & $[1, 1, 1]$ & $h_0 M d_0 e_0$ & $h_0 \Delta^2 d_0 e_0 + \Delta h_1 e_0^2 g$ \\
  80 & 14 & $[1, 0]$ & $[0, 1]$ & $h_0 g B_4$ & $h_0 \Delta^2 m$\\
  80 & 14 & $[0, 1]$ & $[0, 0]$ & $\Delta^3 h_1 h_3$ & $0$\\
  84 & 19 & $[1, 1]$ & $[0, 0, 1]$ & $h_0 MPd_0 e_0$ & $\Delta h_1 d_0^2 e_0^2 + h_0 P \Delta^2 d_0 e_0$ \\
\end{longtable}

\begin{longtable}{cccccc}
  \caption{Selected hidden $\sigma$ extensions}\label{table:hidden-sigma} \\
  \toprule
  $n$ & $s$ & source & target & name of source & name of target \\
  \midrule
  \endfirsthead
  \caption{Selected hidden $\nu$ extensions (continued)} \\
  \toprule
  $n$ & $s$ & source & target & name of source & name of target \\
  \midrule
  \endhead
  \bottomrule
  \endfoot
  84 & 15 & $[1, 1]$ & $[1, 0, 0]$ & $M \Delta h_1 d_0$ & $h_0^6 x_{91, 11}$ \\
  85 & 17 & $[1, 0, 0]$ & $[0, 0]$ & $x_{85, 17, 0}$ & $0$ \\
  85 & 17 & $[0, 1, 1]$ & $[0, 0]$ & $x_{85, 17, 1} + x_{85, 17, 2}$ & $0$ \\
\end{longtable}

\begin{longtable}{ccccc}
  \caption{Identification of classes}\label{table:class-ident} \\
  \toprule
  $n$ & $s$ & class & name & identification\\
  \midrule
  \endfirsthead
  \caption{Identification of classes (continued)} \\
  \toprule
  $n$ & $s$ & class & name & identification\\
  \midrule
  \endhead
  \bottomrule
  \endfoot
  62 & 10 & $[1, 0, 0]$ & $h_1 \Delta x$ & $h_1$-divisible \\
  62 & 10 & $[0, 1, 0]$ & $\Delta^2 h_3^2$ & $h_1$-torsion \\
  68 & 12 & $[1, 0]$ & $h_0 h_2 Mg$ & $h_2$-divisible \\
  69 & 8 & $[1, 0]$ & $D_3'$ & $h_0$-torsion \\
  79 & 17 & $[1, 0, 0]$ & $h_0 \Delta^2 d_0 e_0$ & $h_0$-divisible \\
  79 & 17 & $[0, 1, 1]$ & $\Delta h_1 e_0^2 g$ & $h_0$-torsion \\
  80 & 14 & $[0, 1]$ & $\Delta^3 h_1 h_3$ & $h_0$-torsion \\
  80 & 14 & $[1, 0]$ & $h_0 g B_4$ & $h_0$-divisible \\
  83 & 16 & $[0, 1]$ & $h_0 \Delta^2 m$ & $h_0$-divisible \\
  84 & 15 & $[1, 1]$ & $M\Delta h_1 d_0$ & $h_0$-torsion \\
  84 & 19 & $[1, 1]$ & $h_0 MPd_0 e_0$ & $h_2$-divisible \\
  85 & 17 & $[?, 1, 0]$ & $Md_0 j$ & $d_0$-divisible \\
  87 & 21 & $[0, 1, 0]$ & $h_0 P \Delta^2 d_0 e_0$ & $h_0$-divisible \\
  87 & 21 & $[0, 1, 1]$ & $\Delta h_1 d_0^2 e_0^2$ & $h_0$-torsion \\
  91 & 17 & $[1, 0, 0]$ & $h_0^6 x_{91, 11}$ & $h_0$-divisible \\
  91 & 17 & $[1, 0, 1]$ & $Md_0 \ell$ & $d_0$-divisible \\
  92 & 14 & $[0, 1, 0]$ & $h_0^2 \Delta^2 g_2$ & $h_0$-divisible \\
  92 & 19 & $[1, 1]$ & $e_0 g^2 m$ & $g$-divisible \\
  93 & 20 & $[1, 0]$ & $MP\Delta h_1^2 d_0$ & $h_1$-divisible \\
  94 &  9 & $[1, 0, 0]$ & $h_0 x_{94, 8}$ & $h_0$-divisible \\
  94 & 20 & $[0, 1]$ & $MP\Delta h_0^2 e_0 + ? e_0^3 g$ & $h_0$-non-torsion \\
  95 &  7 & $[?, 1]$ & $x_{95, 7}$ & non-$h_6$-divisible \\
\end{longtable}

\printbibliography
\end{document}